\documentclass[11pt,a4]{article}
\usepackage[T2A]{fontenc}
\usepackage[cp866]{inputenc}
\usepackage[mathscr]{eucal}
\usepackage{amssymb}
\usepackage{amsmath}
\usepackage{graphicx}
\usepackage[russian]{babel}  
\usepackage{array}  

\renewcommand{\leq}{\leqslant}

\newcommand{\R}{{\mathbb R}}

\newcommand{\K}{{\mathbb K}}
\renewcommand{\P}{{\mathbb P}}

\renewcommand{\k}{\rule{0.7em}{0.7em}}
\usepackage{fancyhdr}
\begin{document}
\sloppy

\lhead
    [\scriptsize V.N. Gorbuzov]
    {\scriptsize V.N. Gorbuzov}
\rhead
    [\it \scriptsize Stereographically conjugate differential systems]
    {\it \scriptsize Stereographically conjugate differential systems}

{\normalsize

\thispagestyle{empty}

\mbox{}
\\[-0.15ex]
\centerline{
{\large
\bf
STEREOGRAPHICALLY CONJUGATE DIFFERENTIAL SYSTEMS 
}
}
\\[2.5ex]
\centerline{
\bf 
V.N. Gorbuzov
}
\\[2ex]
\centerline{
\it
Department of Mathematics and Computer Science,
}
\\[1ex]
\centerline{
\it
Yanka Kupala Grodno State University,
}
\\[1ex]
\centerline{
\it
Ozeshko {\rm 22}, Grodno, {\rm 230023}, Belarus
}
\\[1.5ex]
\centerline{
E-mail: gorbuzov@grsu.by
}
\\[5.5ex]
\centerline{{\large\bf Abstract}}
\\[1ex]
\indent
Topological bases of behaviour of trajectories for autonomous differential systems of the second order 
on sphere are stated. Stereographic atlas of trajectories is constructed. 
Differential connections between trajectories of stereographically conjugate differential systems are established. 
The behaviour of trajectories in an neighbourhood of infinitely remote point of the phase plane is investigated. 
Examples of global qualitative research of trajectories of stereographically conjugate differential systems are given.
\\[1.5ex]
\\[1.25ex]
\\[7.5ex]
\centerline{{\large\bf Contents}}
\\[1.5ex]
{\bf  Introduction}                   \dotfill\ 2
\\[1ex]
{\bf \S 1. Stereographic atlas of sphere}
                                                 \dotfill \ 2
\\[0.5ex]
\mbox{}\hspace{1.35em}
1. Stereographic projection of plane
                                                 \dotfill \ 2
\\[0.5ex]
\mbox{}\hspace{1.35em}
2. Stereographic atlas of sphere
                                                 \dotfill \ 4
\\[1ex]
\noindent
{\bf \S 2. 
Stereographically conjugate differential system
}
                                                 \dotfill \ 6
\\[0.75ex]
\mbox{}\hspace{1.35em}
3. Bendixon's transformation
                                                 \dotfill \ 6
\\[0.5ex]
\mbox{}\hspace{1.35em}
4. Form of stereographically conjugate differential system 
                                                 \dotfill \ 7
\\[0.5ex]
\mbox{}\hspace{1.35em}
5. Stereographic atlas of trajectories for differential system
                                                 \dotfill \ 9
\\[1ex]
\noindent
{\bf \S 3. Trajectories of stereographically conjugate differential systems
}
                                                 \dotfill \ 16
\\[0.75ex]
\mbox{}\hspace{1.35em}
6. Regular points and equilibrium states of stereographically conjugate
\\
\mbox{}\hspace{2.6em}
differential systems
                                                 \dotfill \ 16
\\[0.5ex]
\mbox{}\hspace{1.35em}
7. Stereographic cycles
                                                 \dotfill \ 20
\\[0.5ex]
\mbox{}\hspace{1.35em}
8. Symmetry of phase directional field for stereographically conjugate
\\ 
\mbox{}\hspace{2.6em}
differential systems
                                                 \dotfill \ 24
\\[0.5ex]
\mbox{}\hspace{1.35em}
9. Infinitely remote equilibrium state
                                                 \dotfill \ 25
\\[1ex]
{\bf References}
                                              \dotfill \ 31

\newpage

\mbox{}
\\[-2.5ex]
\centerline{\large\bf  Introduction}
\\[1ex]
\indent
Research object is ordinary autonomous polynomial differential system of the second order
\\[1.5ex]
\mbox{}\hfill                                   
$
\displaystyle 
\dfrac{dx}{dt} =
\sum \limits_{k=0}^{n}\, X_k^{}(x,y)\equiv
 X(x,y), 
\qquad 
\dfrac{dy}{dt} =
\sum \limits_{k=0}^{n}\,Y_k^{}(x,y)\equiv
Y(x,y),
$
\hfill (D)
\\[1.75ex]
where $X_k^{}$ and $Y_k^{}$ 
\vspace{0.5ex}
are homogeneous polynomials of degree $k,\ k=0,1,\ldots,n,$ on variables $x$ and $y$
such that 
\vspace{0.35ex}
$|X_n^{}(x,y)|+|Y_n^{}(x,y)|\not\equiv0$ on $\R^2,$ and polynomials $X$ and $Y$ are relatively prime, 
\vspace{0.35ex}
i.e. they have no the common divisors which are distinct from real numbers.

This paper is a continuation of researches stated in [1] and 
\vspace{0.15ex}
the main results of this paper were originally published by the author in 
the articles [2] and [3]. 
\\[3.5ex]
\centerline{
{\bf\large \S\;\!1. Stereographic atlas of sphere}}
\\[1.5ex]
\centerline{
{\bf  1. 
Stereographic projection of plane
}
}
\\[1ex]
\indent
Let's introduce three-dimensional rectangular Cartesian coordinate system 
\vspace{0.15ex}
$O ^ {\star} x ^ {\star} y ^ {\star} z ^ {\star}, $ 
combined with the right rectangular Cartesian coordinate system $Oxy, $ 
\vspace{0.15ex}
meeting conditions: 
the straight line $OO ^ {\star} $ is orthogonal to plane $Oxy,$ 
\vspace{0.15ex}
the length of the segment $OO ^ {\star} $ is equal to one unit of the scale of the system of coordinate $Oxy;$ 
\vspace{0.15ex}
the axis $O ^ {\star} x ^ {\star} $ is codirected with the axis $Ox,$ 
the axis $O ^ {\star} y ^ {\star} $ is codirected with the axis $Oy,$ 
\vspace{0.15ex}
and the axis $O ^ {\star} z ^ {\star} $ is directed so that the system of coordinate $O ^ {\star} x ^ {\star} y ^ {\star} z ^ {\star}$ will be right; 
\vspace{0.15ex}
a scale in the system of coordinate $O ^ {\star} x ^ {\star} y ^ {\star} z ^ {\star} $ same, as well as in the system of coordinate $Oxy.$
\vspace{0.25ex}
We will construct the sphere with the centre $O^ {\star}$ of unit radius: 
\\[2.25ex]
\mbox{}\hfill                                   
$
S^2=\bigl\{
(x^{\star},y^{\star},z^{\star})
\colon
x^{\star}{}^{\,^{\scriptstyle 2}}+
y^{\star}{}^{\,^{\scriptstyle 2}}+
z^{\star}{}^{\,^{\scriptstyle 2}}=1\bigr\}.
$
\hfill (1.1)
\\[2.25ex]
\noindent
Points $N (0,0,1) $ and $S (0,0, {}-1) $ 
\vspace{0.25ex}
are according to northern and southern poles of this sphere. 
\vspace{0.25ex}
Thus the southern pole $\!S(0,0, \!{}-1)\!$ coincides with the beginning $\!O(0,0)\!$ of the system of coordinate $Oxy. $
\vspace{0.25ex}
The equation $\!z ^ {\star}\! = -\;\!1\!$ is the equation in the system of coordinate 
$\!O ^ {\star} x ^ {\star} y ^ {\star} z ^ {\star}\!$ of the plane $\!Oxy.\!$
The plane $Oxy$ 
\vspace{0.35ex}
concerns by sphere (1.1) in the southern pole $S (0,0, {}-1). $

On the plane $Oxy $ 
\vspace{0.15ex}
arbitrarily we will choose the point $M(x,y)$ and we will spend the ray with the beginning 
$M$ through the northern pole $N.$ 
\vspace{0.15ex}
The ray $MN$ intersects the sphere (1.1) in some point $P.$
\vspace{0.15ex}
Thereby to each point of the plane $Oxy $ there correspond one point of the sphere (1.1),
\vspace{0.15ex}
and to each point of  the sphere (1.1), except the northern pole,  
there correspond  one point of the plane $Oxy.$ 
\vspace{0.15ex}
Such projection is called
{\it stereographic projection} of plane on sphere [4, pp. 83 --- 84].
\vspace{0.25ex}
The point $N$ is said to be the {\it centre} of stereographic projection (Fig. 1.1).
\vspace{0.5ex}

{\bf Lemma 1.1.}
{\it
Stereographic projection of plane is 
bijection between plane and sphere without northern pole {\rm(}centre of this projection{\rm)}.
}
\vspace{0.15ex}

To spread correspondence to all sphere (1.1), we will arrive as follows.
On the plane $Oxy$ we introduce conditional
{\it infinitely remote point} $M^{}_{\infty},$ which is 
pre-image of the northern pole $N(0,0,1)$ by the stereographic projection.
\vspace{0.15ex}

The plane $Oxy$ plas the infinitely remote point $M^{}_{\infty}$
\vspace{0.15ex}
(the image of this point under the stereographic projection of the plane $Oxy$ on the sphere (1.1) 
\vspace{0.35ex}
is the northern pole  $N(0,0,1))$ 
is called the {\it extended plane} $Oxy$ and is denoted by $\overline{Oxy},$ i.e. $\overline{Oxy} = Oxy\sqcup M^{}_{\infty}\,.$
\vspace{0.5ex}

The extended plane $\overline{Oxy}$ consists of two parts: final (the plane $Oxy)$ and
\vspace{0.25ex}
infinite (the point $M^{}_{\infty}).$ 
\vspace{0.15ex}
Then, 
by the stereographic projection, 
for any sequence  $\{M_k^{}\}$ of points $M_k^{}$ of the plane $Oxy$ 
such that these points are leaving in infinity 
\vspace{0.25ex}
(tend to the point  $M^{}_{\infty}),$ we get 
the unique correspondence sequence 
\vspace{0.35ex}
$\{P_k^{}\}$ of points $P_k^{}$ of the sphere (1.1) 
such that these points are tend to the northern pole  $N(0,0,1).$ 

Suppose the point $M$  is lying on the plane $Oxy$ and has coordinates $M(x,y).$ 
Then, this point in space coordinate system $O^{\star}x^{\star}y^{\star}z^{\star}$ has the coordinates $M(x,y,{}-1).$ 
The straight line $MN$ in the coordinate system $O^{\star}x^{\star}y^{\star}z^{\star}$ is the system of equations
\\[2ex]
\mbox{}\hfill
$
\dfrac{x^{\star}}{x}=\dfrac{y^{\star}}{y}=\dfrac{z^{\star}-1}{{}-2}\,.
\hfill
$
\\[2ex]
Using the stereographic projection with the centre in the northern pole $N(0,0,1),$ we obtain 
to the point $M(x,y,{}-1)$ corresponds the point $P(x^{\star},y^{\star},z^{\star}),$ which is
the point of intersection of the straight line $MN$ and the sphere (1.1). 
Therefore the coordinates $x^{\star},\  y^{\star},\  z^{\star}$ of the point $P$ are the solutions 
of the algebraic system of equations 
\\[2ex]
\mbox{}\hfill                 
$
\dfrac{x^{\star}}{x} =\dfrac{y^{\star}}{y}=\dfrac{1-z^{\star}}{2}\,,
\qquad 
x^{\star}{}^{\,^{\scriptstyle 2}}+
y^{\star}{}^{\,^{\scriptstyle 2}}+
z^{\star}{}^{\,^{\scriptstyle 2}}=1.
$
\hfill (1.2)
\\[2.25ex]
\indent
Assume that the point 
$P(x^{\star},y^{\star},z^{\star})$ lies on the sphere (1.1) and is not the northern pole  $N(0,0,1).$
Then the applicate of this point is $z^{\star}\in [{}-1;1).$
Further, we solve the system of equations (1.2) for 
$x^{\star},\ y^{\star},\ z^{\star}$ at ${}-1\leq z^{\star}<1,$
and have the bijective reflexion 
\\[2ex]
\mbox{}\hfill                                   
$
\psi_{_N} \colon (x,y)\to \bigl(x^{\star}(x,y), y^{\star}(x,y), z^{\star}(x,y)\bigr),
\hfill                                  
$
\\[-0.75ex]
\mbox{}\hfill {\rm (1.3)}
\\[0.75ex]
\mbox{}\hfill 
$
x^{\star}(x,y)=\dfrac{4x}{x^2+y^2+4}\,, \ 
y^{\star}(x,y)=\dfrac{4y}{x^2+y^2+4}\,, \ 
z^{\star}(x,y)=\dfrac{x^2+y^2-4}{x^2+y^2+4}
$
for all 
$
(x,y)\in \R^2
\hfill 
$
\\[2.75ex]
of the plane $Oxy$  on the sphere (1.1) without the northern pole  $N(0,0,1).$

The coordinate functions of map (1.3) are continuously differentiable. The Jacobians 
\\[2.25ex]
\mbox{}\hfill                                  
$
\dfrac{{\sf D} (x^{\star},y^{\star})}
{{\sf D} (x,y)}=
{}-16\, \dfrac{x^2+y^2-4}{(x^2+y^2+4)^3}\,,
\qquad 
\dfrac{{\sf D} (x^{\star},z^{\star})}{{\sf D} (x,y)}=
64\, \dfrac{y}{(x^2+y^2+4)^3}\,,
\hfill                                  
$
\\[2.75ex]
\mbox{}\hfill                                  
$
\dfrac{{\sf D} (y^{\star},z^{\star})}{{\sf D} (x,y)}=
{}-64\, \dfrac{x}{(x^2+y^2+4)^3}
$ 
\ \ for all 
$
(x,y) \in \R^2
\hfill
$
\\[2.5ex]
at the same time are nonvanishing in any point of the plane $Oxy.$

Hence, taking into account Lemma 1.1, we obtain the next statement.

{\bf Theorem 1.1.}
{\it
The stereographic map {\rm (1.3)}  of the plane $Oxy$ on the sphere {\rm (1.1)} without the northern pile  $N(0,0,1)$
is a diffeomorphism.
}

The basic property of the diffeomorphism (1.3) [5, columns 222 -- 223]: 
an angle between curves on the plane equals to the angle between stereographic images of these curves on the sphere. 
\\[5ex]
\mbox{}\hfill
{\unitlength=1mm
\begin{picture}(52,56)
\put(-7,0){\includegraphics[width=52mm,height=56mm]{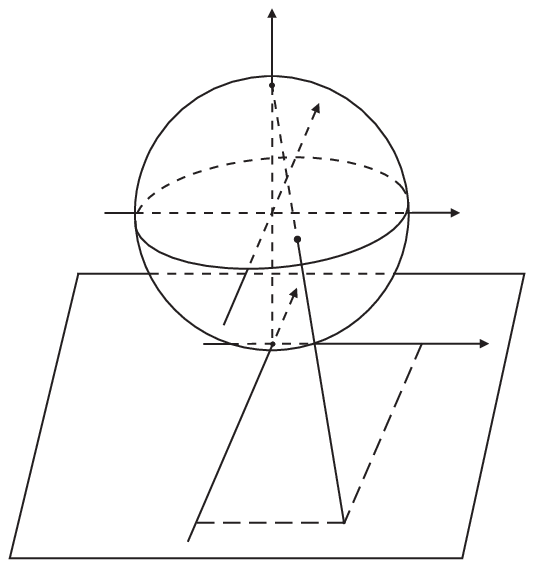}}

\put(20.9,19.7){\makebox(0,0)[cc]{\scriptsize $O$}}
\put(17.7,23.7){\makebox(0,0)[cc]{\scriptsize $S$}}
\put(35.5,20.5){\makebox(0,0)[cc]{\scriptsize $x$}}
\put(40,20){\makebox(0,0)[cc]{\scriptsize $x$}}
\put(10,3.9){\makebox(0,0)[cc]{\scriptsize $y$}}
\put(22,25){\makebox(0,0)[cc]{\scriptsize $y$}}
\put(28,2.6){\makebox(0,0)[cc]{\scriptsize $M$}}

\put(17.2,37){\makebox(0,0)[cc]{\scriptsize $O^{\star}$}}
\put(24,32.5){\makebox(0,0)[cc]{\scriptsize $P$}}
\put(37.6,33.3){\makebox(0,0)[cc]{\scriptsize $x^{\star}$}}
\put(26,44.3){\makebox(0,0)[cc]{\scriptsize $y^{\star}$}}

\put(17.6,47){\makebox(0,0)[cc]{\scriptsize $N$}}
\put(21.5,54.5){\makebox(0,0)[cc]{\scriptsize $z^{\star}$}}

\put(20,-5){\makebox(0,0)[cc]{\rm Fig. 1.1}}
\end{picture}}
\qquad\qquad
{\unitlength=1mm
\begin{picture}(50,55)
\put(0,0){\includegraphics[width=49.79mm,height=55.84mm]{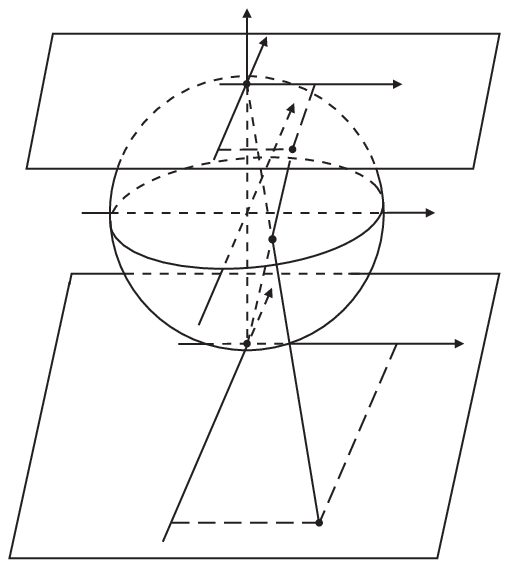}}

\put(25.9,19.7){\makebox(0,0)[cc]{\scriptsize $O$}}
\put(22.7,23.7){\makebox(0,0)[cc]{\scriptsize $S$}}
\put(40,20.5){\makebox(0,0)[cc]{\scriptsize $x$}}
\put(44.7,20){\makebox(0,0)[cc]{\scriptsize $x$}}
\put(15,3.9){\makebox(0,0)[cc]{\scriptsize $y$}}
\put(27,25){\makebox(0,0)[cc]{\scriptsize $y$}}
\put(34,2.6){\makebox(0,0)[cc]{\scriptsize $M$}}

\put(21.7,33.5){\makebox(0,0)[cc]{\scriptsize $O^{\star}$}}
\put(28.9,32.5){\makebox(0,0)[cc]{\scriptsize $P$}}
\put(42.6,33.3){\makebox(0,0)[cc]{\scriptsize $x^{\star}$}}
\put(27,46){\makebox(0,0)[cc]{\scriptsize $y^{\star}$}}

\put(31.9,42){\makebox(0,0)[cc]{\scriptsize $M^{\ast}$}}
\put(21.3,46.3){\makebox(0,0)[cc]{\scriptsize $O^{\ast}$}}
\put(22.7,50.3){\makebox(0,0)[cc]{\scriptsize $N$}}
\put(38.3,46){\makebox(0,0)[cc]{\scriptsize $u$}}
\put(31.5,49.8){\makebox(0,0)[cc]{\scriptsize $u$}}
\put(27.9,51.5){\makebox(0,0)[cc]{\scriptsize $v$}}
\put(20,41.9){\makebox(0,0)[cc]{\scriptsize $v$}}
\put(26.5,55.5){\makebox(0,0)[cc]{\scriptsize $z^{\star}$}}

\put(25,-5){\makebox(0,0)[cc]{\rm Fig. 2.1}}
\end{picture}}
\hfill\mbox{}
\\[-7ex]

\newpage

\mbox{}
\\[-1.5ex]
\centerline{\bf  2. 
Stereographic atlas of sphere }
\\[1ex]
\indent
On the plane, that concerning the sphere (1.1) in the northern pole $N(0,0,1),$
\vspace{0.15ex}
we introduce the right rectangular Cartesian coordinate system
\vspace{0.25ex}
$O^{\ast}uv$ such that its origin $O^{\ast}(0,0)$ coincides with northern pole
$N(0,0,1)$ of sphere (1.1), the axis $O^{\ast}u$
\vspace{0.25ex}
is codirected with the axis $O^{\star}x^{\star},$ 
the axis $O^{\ast}v$ is codirected with the axis $O^{\star}y^{\star}$ (Fig. 2.1).
\vspace{0.25ex}
Scale in  the coordinate systems $Oxy, \ O^{\star}x^{\star}y^{\star}z^{\star},$ and $O^{\ast}uv$  is identical.
\vspace{0.35ex}

If the point $M^{\ast}$ in the plane $O^{\ast}uv$ has coordinates $M^{\ast}(u,v),$ 
\vspace{0.25ex}
then this point in the space coordinate system $O^{\star}x^{\star}y^{\star}z^{\star}$ has 
coordinates $M^{\ast}(u,v,1).$ 
\vspace{0.25ex}
The straight line $M^{\ast}S$ in the coordinate system $O^{\star}x^{\star}y^{\star}z^{\star}$ 
is defined by the system of equations
\\[2ex]
\mbox{}\hfill
$
\dfrac{x^{\star}}{u}=\dfrac{y^{\star}}{v}=\dfrac{z^{\star}+1}{2}\,.
\hfill
$
\\[2.25ex]
\indent
The image of the point $M^{\ast}(u,v,1)$ 
\vspace{0.25ex}
by the stereographic projection with the centre in southern pole $S(0,0,{}-1)$ of the sphere (1.1) 
is the point  $P(x^{\star},y^{\star},z^{\star}),$ 
\vspace{0.25ex}
which is the point of intersection of the straight line $M^{\ast}S$ and the sphere (1.1).
\vspace{0.25ex}
Therefore the coordinates $x^{\star}, y^{\star}, z^{\star}$ of the point $P$ are solutions of algebraic system of equations
\\[2ex]
\mbox{}\hfill                 
$
\dfrac{x^{\star}}{u} =\dfrac{y^{\star}}{v}=\dfrac{z^{\star}+1}{2}\,,
\quad \ 
x^{\star}{}^{\,^{\scriptstyle 2}}+
y^{\star}{}^{\,^{\scriptstyle 2}}+
z^{\star}{}^{\,^{\scriptstyle 2}}=1.
$
\hfill (2.1)
\\[2.25ex]
\indent
Let us consider the point $P(x^{\star},y^{\star},z^{\star})$ 
\vspace{0.25ex}
such that this point lies on the sphere (1.1) and 
this point is not southern pole  $S(0,0,{}-1).$
Then its applicate $z^{\star}\in ({}-1;1].$
\vspace{0.35ex}

Having resolved the system of equations (2.1) rather $x ^ {\star}, \, y ^ {\star}, \, z ^ {\star} $ at $ {}-1< z ^ {\star} \leq 1, $ 
we receive the bijective reflexion 
\\[2ex]
\mbox{}\hfill                                   
$
\psi_{_S} \colon (u,v)\to\ \bigl(x^{\star}(u,v), y^{\star}(u,v), z^{\star}(u,v)\bigr),
\hfill                                  
$
\\[-0.75ex]
\mbox{}\hfill {\rm (2.2)}
\\[0.75ex]
\mbox{}\hfill 
$
x^{\star}(u,v)\!=\dfrac{4u}{u^2+v^2+4}\,, \ 
y^{\star}(u,v)\!=\dfrac{4v}{u^2+v^2+4}\,, \ 
z^{\star}(u,v)\!=\!{}-\dfrac{u^2+v^2-4}{u^2+v^2+4}
$\,
for all 
$\!(u,v)\!\in\! \R^2
\hfill 
$
\\[2.75ex]
of the plane $O^{\ast}uv $ on the sphere (1.1) without the southern pole  $S(0,0,{}-1).$
\vspace{0.35ex}

The coordinate functions of map (2.2) are continuously differentiable.
The Jacobians 
\\[2.25ex]
\mbox{}\hfill                                  
$
\dfrac{{\sf D} (x^{\star},y^{\star})}
{{\sf D} (u,v)}=
{}-16\, \dfrac{u^2+v^2-4}{(u^2+v^2+4)^3}\,,
\qquad 
\dfrac{{\sf D} (x^{\star},z^{\star})}{{\sf D} (u,v)}=
{}-64\, \dfrac{v}{(u^2+v^2+4)^3}\,,
\hfill                                  
$
\\[2.75ex]
\mbox{}\hfill                                  
$
\dfrac{{\sf D} (y^{\star},z^{\star})}{{\sf D} (u,v)}=
64\, \dfrac{u}{(u^2+v^2+4)^3}
$ 
\ \ for all 
$
(u,v) \in \R^2
\hfill
$
\\[2.5ex]
at the same time are nonvanishing in any point of the plane $O^{\ast}uv.$
\vspace{0.35ex}

Therefore 
\vspace{0.25ex}
the stereographic map (2.2) of the plane $O^{\ast}uv$ 
on the sphere (1.1) without the southern pole  $S(0,0,{}-1)$
is a diffeomorphism.
\vspace{0.35ex}

The plane $\!O^{\ast}uv\!$ plas the infinitely removed point $\!\!M^{\ast}_{\infty}\!\!$ 
\vspace{0.35ex}
(the image of this point under the stereographic projection of the plane $O^{\ast}uv$ on the sphere (1.1) 
\vspace{0.5ex}
is the southern pole $\!S(0,0,\!-1))$ 
is called the {\it extended plane} $O^{\ast}uv$  
and is denoted by $\overline{O^{\ast}uv},$ i.e., 
$\overline{O^{\ast}uv} = O^{\ast}uv\sqcup M^{\ast}_{\infty}\,.$
\vspace{0.35ex}

Let us cover  the sphere (1.1) by two spheres $U_1^{}$ and $U_2^{},$ 
\vspace{0.35ex}
without accordingly the northern pole $N(0,0,1)$ and the southern pole $S(0,0,{}-1)\colon$
\\[2ex]
\mbox{}\hfill
$
\displaystyle
U_1^{}=\Bigl\{
(x^{\star}, y^{\star}, z^{\star})\colon\,   
{x^{\star}}^{\,^{\scriptstyle 2}}+{y^{\star}}^{\,^{\scriptstyle 2}}+{z^{\star}}^{\,^{\scriptstyle 2}}=1, \ \,
{}-1\leq z^{\star}<1\Bigr\}
\hfill
$
\\[1ex]
and
\\[1ex]
\mbox{}\hfill
$
\displaystyle
U_2^{}=\Bigl\{
(x^{\star}, y^{\star}, z^{\star})\colon\,   
{x^{\star}}^{\,^{\scriptstyle 2}}+{y^{\star}}^{\,^{\scriptstyle 2}}+{z^{\star}}^{\,^{\scriptstyle 2}}=1, \ \,
{}-1< z^{\star}\leq1\Bigr\}.
\hfill
$
\\[2.25ex]
\indent
We introduce the diffeomorphic maps 
\\[2ex]
\mbox{}\hfill                                                              
$
\displaystyle
\varphi_1^{}\colon (x^{\star}, y^{\star}, z^{\star})\to \
\bigl(
x(x^{\star}, y^{\star}, z^{\star}),\, y(x^{\star}, y^{\star}, z^{\star})
\bigr),
\hfill
$
\\[0ex]
\mbox{}\hfill (2.3)
\\[0ex]
\mbox{}\hfill
$
\displaystyle
x(x^{\star}, y^{\star}, z^{\star})=\dfrac{2x^{\star}}{1-z^{\star}}\,, 
\qquad 
y(x^{\star}, y^{\star}, z^{\star})=\dfrac{2y^{\star}}{1-z^{\star}}
$
\ for all 
$
(x^{\star}, y^{\star}, z^{\star})\in U_1^{},
\hfill
$
\\[3ex]
\mbox{}\hfill                                                              
$
\displaystyle
\varphi_2^{}\colon (x^{\star}, y^{\star}, z^{\star})\to\ 
\bigl(
u(x^{\star}, y^{\star}, z^{\star}),\, v(x^{\star}, y^{\star}, z^{\star})
\bigr),
\hfill
$
\\[0ex]
\mbox{}\hfill (2.4)
\\[0ex]
\mbox{}\hfill
$
\displaystyle
u(x^{\star}, y^{\star}, z^{\star})=\dfrac{2x^{\star}}{1+z^{\star}}\,, 
\qquad 
v(x^{\star}, y^{\star}, z^{\star})=\dfrac{2y^{\star}}{1+z^{\star}}
$
\ for all 
$
(x^{\star}, y^{\star}, z^{\star})\in U_2^{}.
\hfill
$
\\[2.5ex]
\indent
The map (2.3) is inverse to the map (1.3). 
\vspace{0.25ex}
Therefore the map (2.3) is diffeomorphic map of the sphere $U_1^{}$ without the northern pole $N(0,0,1)$ on the plane $Oxy,$
\vspace{0.25ex}
i.e, it is stereographic map of the sphere (1.1) on the plane $Oxy$ 
\vspace{0.25ex}
from the centre in the northern pole $N(0,0,1)$
\linebreak
\text{[5, columns 222 -- 223; 6, p. 37].}
The map (2.4) is inverse to the map (2.2). 
\vspace{0.25ex}
Therefore the  map (2.4)  is diffeomorphic map of the sphere $U_2^{}$ 
\vspace{0.25ex}
without the southern pole $S(0,0,{}-1)$ on the plane $O^{\ast}uv,$
\vspace{0.25ex}
i.e, it is stereographic map of the sphere (1.1) on the plane $O^{*}uv$ from the centre in the southern pole $S(0,0,{}-1).$
\vspace{0.35ex}

Thus two charts $(U_1^{}, \varphi_1^{}) $ and $(U_2^{}, \varphi_2^{}) $ of sphere (1.1) are constructed. 
\vspace{0.35ex}
The set of charts $(U_1^{}, \varphi_1^{}) $ and $(U_2^{}, \varphi_2^{}) $
are  a {\it stereographic atlas} of sphere  (1.1) [6, p. 103].
\vspace{0.5ex}

Let us establish connections between 
\vspace{0.15ex}
the local coordinate systems $Oxy$ and $O^{*}uv$ of
stereographic atlas of sphere  (1.1). 
\vspace{0.15ex}
The stereographic maps (1.3) and (2.2) are inverse to the maps (2.3) and (2.4) of the stereographic atlas of sphere  (1.1),
respectively, i.e.,
\\[2ex]
\mbox{}\hfill                                                      
$
\varphi_1^{}=\psi_{_N}^{{}-1},
\ \ 
\psi_{_N}=\varphi_1^{{}-1}
$
\ \ 
and 
\ \ 
$
\varphi_2^{}=\psi_{_S}^{{}-1},
\ \ 
\psi_{_S}=\varphi_2^{{}-1}.
\hfill
$
\\[2.5ex]
\indent
Combining sequentially two diffeomorphic maps  $\psi_{_N}$ and $\psi_{_S}^{{}-1},$ 
we obtain the diffeomorphic map (Fig. 2.1)
\\[1.25ex]
\mbox{}\hfill                                                     
$
\varphi_{21}^{}=\psi_{_S}^{{}-1}\circ\psi_{_N}= \varphi_2^{}\circ\psi_{_N}
$
\hfill (2.5)
\\[2.25ex]
of the plane $Oxy$ without the origin $O(0,0)$ 
\vspace{0.5ex}
on the plane $O^{\ast}uv$ without the origin $O^{\ast}(0,0).$  

Using the analytical representations (1.3) and (2.4) of the maps  $\psi_{_N}$ and $\varphi_{2}^{},$
\vspace{0.25ex}
we get the analytical representation of the diffeomorphic map (2.5) in the form
\\[2ex]
\mbox{}\hfill                           
$
\varphi_{21}^{} \colon (x,y) \to\
\Bigl(\;\!\dfrac{4x}{x^{2}+y^{2}}\,,\
\dfrac{4y}{x^{2}+y^{2}}\Bigr)
$
\ for all 
$
(x,y)\in\R^2\backslash\{(0,0)\}.
$
\hfill (2.6)
\\[2ex]
\indent
The coordinate functions 
\\[2ex]
\mbox{}\hfill                           
$
u \colon (x,y) \to\ \dfrac{4x}{x^{2}+y^{2}}\,,
\quad\
v \colon (x,y) \to\ \dfrac{4y}{x^{2}+y^{2}} 
$
\ \ for all 
$
(x,y)\in\R^2\backslash\{(0,0)\}
$
\hfill (2.7)
\\[2ex]
of the diffeomorphic map  (2.6) are 
\vspace{0.25ex}
functions of transition  [6, p. 99] from the local coordinates  $(u,v)$ to the local coordinates $(x,y)$ 
of the stereographic atlas of sphere  (1.1).
\vspace{0.35ex}

The diffeomorphic map 
\\[2ex]
\mbox{}\hfill                           
$
\varphi_{12}^{} \colon (u,v) \to\
\Bigl(\;\!\dfrac{4u}{u^{2}+v^{2}}\,,\
\dfrac{4v}{u^{2}+v^{2}}\Bigr)
$
\ for all 
$ 
(u,v)\in\R^2\backslash\{(0,0)\}
$
\hfill (2.8)
\\[2ex]
has the coordinate functions 
\\[2ex]
\mbox{}\hfill                           
$
x \colon (u,v) \to\ \dfrac{4u}{u^{2}+v^{2}}\,,
\quad \ 
y \colon (u,v) \to\ \dfrac{4v}{u^{2}+v^{2}} 
$
\ \ 
for all 
$
(u,v)\in\R^2\backslash\{(0,0)\}.
$
\hfill (2.9)
\\[2ex]
\indent
The coordinate functions (2.9) are 
\vspace{0.25ex}
functions of transition  from the local coordinates  $(x,y)$ to the local coordinates $(u,v)$ 
of the stereographic atlas of sphere  (1.1).

The diffeomorphic map (2.8), which is received as inverse map to the map (2.6), 
is an analytical representation of the map  
\\[1.5ex]
\mbox{}\hfill
$
\varphi_{12}^{}=\psi_{_N}^{{}-1}\circ\psi_{_S}= \varphi_1^{}\circ\psi_{_S}
\hfill
$
\\[2ex]
of the plane $O^{\ast}uv$ without the origin $O^{\ast}(0,0)$ 
\vspace{0.75ex}
on the plane $Oxy$ without the origin $O(0,0).$

{\bf Theorem 2.1.}
\vspace{0.15ex}
{\it
The diffeomorphic map {\rm (2.6)} and the identity map 
of the plane $Oxy$ on itself are a group of the second order.}
\vspace{0.25ex}

{\sl Indeed}, the transformation (2.6) is bijective
\\[1.5ex]
\mbox{}\hfill
$
(x,y)\ \stackrel{\stackrel{\scriptstyle \varphi_{21}^{}}{\mbox{}}}
{\longrightarrow}\ 
\Bigl(\;\!\dfrac{4x}{x^2+y^2}\,,\
\dfrac{4y}{x^2+y^2}\Bigr)\ 
\stackrel{\stackrel{\scriptstyle \varphi_{21}^{}}{\mbox{}}}
{\longrightarrow}\ 
\left(\dfrac{4\cdot\dfrac{4x}{x^2+y^2}}
{\Bigl(\dfrac{4x}{x^2+y^2}\Bigr)^2
+\Bigl(\dfrac{4y}{x^2+y^2}\Bigr)^2}\,,
\right.
\hfill
$
\\[2.75ex]
\mbox{}\hfill
$
\left.
\dfrac{4\cdot\dfrac{4y}{x^2+y^2}}
{\Bigl(\dfrac{4x}{x^2+y^2}\Bigr)^2
+\Bigl(\dfrac{4y}{x^2+y^2}\Bigr)^2}\right)=(x,y)
$
\ for all 
$
(x,y)\in\R^2\backslash \{(0,0)\},
\hfill
$
\\[2ex]
and also, 
$
\varphi_{21}^{}\circ I=I\circ \varphi_{21}^{}=\varphi_{21}^{}$ and $I\circ I=I,
$
where $I\colon (x,y)\to (x,y)$ for all $(x,y)\in\R^2.\k$
\\[5ex]
\centerline{\bf\large 
\S\;\!2. Stereographically conjugate differential system}
\\[1.5ex]
\centerline{
\bf  3. Bendixon's transformation 
}
\\[1.15ex]
\indent
{\it Bendixon's transformation} of the phase plane $Oxy$ of the differential system (D) is the transformation 
\\[1.5ex]
\mbox{}\hfill                           
$
x = \dfrac{4u}{u^{2}+v^{2}}\,,
\qquad 
y =\dfrac{4v}{u^{2}+v^{2}} \,.
$
\hfill (3.1)
\\[2ex]
This transformation is constructed by the function of transition (2.9) from the local coordinates $x,\ y$ to the local 
co\-or\-di\-na\-tes $u,\ v$ 
\vspace{0.35ex}
of the stereographic atlas of sphere (1.1).

By Bendixon's transformation (3.1), we reduces the differential system (D) to the differential system [7, p. 239]
\\[2ex]
\mbox{}\hfill                           
$
\dfrac{du}{dt} ={}-\dfrac{\,u^{2}
-v^{2}}{4}\,  X\Bigl(\dfrac{4u}{u^{2}+v^{2}}\,,\
\dfrac{4v}{u^{2}+v^{2}}\Bigr)-
\dfrac{uv}{2}\,   Y\Bigl(\dfrac{4u}{u^{2}+v^{2}}\,,\
\dfrac{4v}{u^{2}+v^{2}}\Bigr)\equiv
\ \stackrel{\ast}{U}(u,v),
\hfill
$
\\[.5ex]
\mbox{}\hfill (3.2)
\\[.5ex]
\mbox{}\hfill
$
\dfrac{dv}{dt}=
{}-\dfrac{u v}{2}\,   X\Bigl(\dfrac{4u}{u^{2}+v^{2}}\,,\
\dfrac{4v}{u^{2}+v^{2}}\Bigr)+
\dfrac{u^{2}-v^{2}}{4}\,
Y\Bigl(\dfrac{4u}{u^{2}+v^{2}}\,,\
\dfrac{4v}{u^{2}+v^{2}}\Bigr)\equiv \ \stackrel{\ast}{V}(u,v).
\hfill
$
\\[2.5ex]
\indent
Since $X$ and $Y$ are polynomials, we see that the system (3.2) has the form
\\[1.5ex]
\mbox{}\hfill                         
$
\dfrac{d u}{dt}=\dfrac{U (u,v )}{(u^2+v^2)^{m}}\,,
\qquad
\dfrac{dv}{dt}=\dfrac{V(u,v)}{(u^2+v^2)^{m}}\,,
\hfill
$
\\[1.5ex]
where $U$ and $V$ are polynomials such that they  are not dividing simultaneously on $u^2+v^2,$ 
and $m$ is an nonnegative integer.
\vspace{0.35ex}

The autonomous polynomial differential system
\\[1.5ex]
\mbox{}\hfill                           
$
\dfrac{d u}{d\tau}=U(u,v),
\qquad
\dfrac{d v}{d\tau}=V(u,v),
$
\hfill (3.3)
\\[2ex]
where $(u^2+v^2)^{m} d\tau=dt,$ 
\vspace{0.35ex}
and the polynomials $U$ and $V$ are relatively prime, is called
 {\it stereographically conjugate}
to the differential system (D). 

Taking into account Theorem 2.1, we obtain  the system (D) 
\vspace{0.15ex}
is stereographically conjugate to the system (3.3). And, using Bendixon's transformation 
\\[1.5ex]
\mbox{}\hfill
$
u = \dfrac{4x}{x^{2}+y^{2}}\,,
\qquad 
v =\dfrac{4y}{x^{2}+y^{2}},
\hfill
$
\\[1.5ex]
we get the system (3.3) is reduced to the system (D).
The differential systems  (D) and (3.3) are stereographically mutually conjugate.
\vspace{0.35ex}

The phase plane $Oxy$ 
\vspace{0.5ex}
(the extended phase plane $\overline{Oxy}\,)$ of the differential system (D) and 
the phase plane $O^{\ast}uv$ (the extended phase plane $\overline{O^{\ast}uv}\,)$ of the differential system (3.3) are 
said to be {\it stereographically conjugate}.
\\[2.25ex]
\centerline{
\bf  4. Form of stereographically conjugate differential system}
\\[1.15ex]
\indent
The form of the differential system (3.3) 
\vspace{0.35ex}
(this system is stereographically conjugate to the system  (D))
depends on whether divides on $x^2+y^2$ or not divides on $x^2+y^2$ the polynomial
\\[1.25ex]
\mbox{}\hfill
$
\displaystyle
W_n^{}\colon (x,y) \to\  x\;\!Y_n^{}(x,y) - yX_n^{}(x,y) 
$
\ for all $(x,y)\in \R^2.
\hfill
$
\\[1.75ex]
\indent
If 
\vspace{0.5ex}
$W_n^{}(x,y)\not\equiv (x^2+y^2)P(x,y)$ on $\R^2,$ 
where $P$ is some polynomial, then the differential system (3.3) has the form  
\\[1.5ex]
$
\displaystyle
\dfrac{du}{d\tau} =
\dfrac{v^2 - u^2}{4}\ \sum\limits_{j = 0}^{n}\,
(u^2 + v^2)^{n - j}\  X_j^{}(4u,4v)  -
 \dfrac{uv}{2}\  \sum\limits_{j = 0}^{n}\,
(u^2 + v^2)^{n - j}\ Y_j^{}(4u,4v) \equiv U_{_0} (u,v),
\hfill
$
\\[0.15ex]
\mbox{}\hfill                   (4.1)
\\[0.15ex]
$
\displaystyle
\dfrac{dv}{d\tau} =
{} - \dfrac{uv}{2}\  \sum\limits_{j = 0}^{n}\,
(u^2+ v^2)^{n - j}\  X_j^{}(4u,4v)  +
\dfrac{u^2 - v^2}{4} \ \sum\limits_{j = 0}^{n}\,
(u^2 + v^2)^{n - j}\ Y_j^{}(4u,4v)\equiv V_{_0} (u,v),
\hfill
$
\\[1.5ex]
where  $(u^2 + v^2)^n\, d\tau = dt.$
\vspace{0.5ex}

Suppose 
\vspace{0.5ex}
$
W_n^{}(x,y)=(x^2+y^2)\;\!P(x,y)
$ 
for all $(x,y)\in \R^2,
$ 
where $P$ is some polynomial and  
the case $P(x,y)=0$ for all $(x,y)\in \R^2$ is not excluded. 
\vspace{0.25ex}

Now if the identities hold
\\[1ex]
\mbox{}\hfill  
$
\displaystyle
{}-2y\bigl(xY_{n-r+1}^{}(x,y) - yX_{n-r+1}^{}(x,y)\bigr)  -
 (x^2 + y^2) X_{n-r+1}^{}(x,y) =
(x^2 + y^2)^{k - r+ 1}\ K_r^{}(x,y),
\hfill
$
\\[2.25ex]
\mbox{}\hfill
$
\displaystyle
2x\bigl(xY_{n-r+1}^{}(x,y) - yX_{n-r+1}^{}(x,y)\bigr) -
 (x^2 + y^2) Y_{n-r+1}^{}(x,y) 
=
(x^2 + y^2)^{k - r + 1}\ Q_r^{}(x,y)
$
\hfill (4.2)
\\[2.25ex]
\mbox{}\hfill
for all 
$
(x,y)\in\R^2, 
\quad r=1,\ldots, k,
\hfill
$
\\[1.75ex]
where the natural number $k$ such that $2k\leq n + 2,$ 
\vspace{0.35ex}
$K_r^{}$ and $Q_r^{},\ r=1,\ldots, k,$ are some polynomials,
then the differential system (3.3) has the form 
\\[1.5ex]
\mbox{}\qquad
$
\displaystyle
\dfrac{du}{d\theta_k^{}} =
\dfrac{v^2 - u^2}{4} \  \sum\limits_{j = 0}^{n-k}\,
(u^2 +  v^2)^{n - j - k} \, X_j^{}(4u,4v)
 -
\dfrac{uv}{2}\  \sum\limits_{j = 0}^{n - k}\,
(u^2 + v^2)^{n - j - k}\ Y_j^{}(4u,4v)  \,+
\hfill
$
\\[1.25ex]
\mbox{}\qquad\qquad\quad
$
\displaystyle
+\,
\sum \limits_{r = 1}^{k}\, 4^{2k-2r-1}\ K_{r}^{}(4u,4v)\equiv U_k^{} (u,v),
\hfill
$
\\[0.25ex]
\mbox{}\hfill (4.3)
\\[0.25ex]
\mbox{}\qquad
$
\displaystyle
\dfrac{dv}{d\theta_k^{}} = {}-\dfrac{uv}{2} \
\sum\limits_{j = 0}^{n - k}\,
(u^2 + v^2)^{n - j - k}\ X_j^{}(4u,4v)  
+
 \dfrac{u^2 - v^2}{4}\ 
\sum\limits_{j = 0}^{n - k}\,
(u^2 + v^2)^{n - j - k}\ Y_j^{}(4u,4v)  \,+
\hfill
$
\\[1.25ex]
\mbox{}\qquad\qquad\quad
$
\displaystyle
+\,
\sum\limits_{r=1}^{k} \, 4^{2k - 2r - 1}\ Q_{r}^{}(4u,4v)\equiv V_k^{} (u,v),
$
\ where $(u^2 + v^2)^{n - k}\, d\theta_k^{}  = dt.
\hfill
$
\\[-2ex]

\newpage

{\bf Example 4.1.} 
Consider the autonomous differential system
\\[2ex]
\mbox{}\hfill
$
\displaystyle
\dfrac{dx}{dt}  = a_{_0}\equiv X(x,y),
\qquad 
\dfrac{dy}{dt}  = b_{_0}\equiv Y(x,y),
\qquad 
| a_{_0}|+| b_{_0}| \ne 0,
$
\hfill (4.4)
\\[2.5ex]
with $W_{_0}(x,y)=b_{_0} x-a_{_0} y\not\equiv (x^2+y^2)P(x,y)$ on $\R^2,$ where $P$ is some polynomial.
\vspace{1ex}

The stereographically conjugate system to the differential system (4.4) is the system 
\\[2.25ex]
\mbox{}\hfill
$
\displaystyle
\dfrac{du}{dt} =
{}-\dfrac{a_{_0}}{4}\; u^2-\dfrac{b_{_0}}{2}\;  uv+\dfrac{a_{_0}}{4}\; v^2,
 \qquad
\dfrac{dv}{dt} =
\dfrac{b_{_0}}{4}\; u^2-\dfrac{a_{_0}}{2}\;  uv-\dfrac{b_{_0}}{4}\; v^2,
\qquad
| a_{_0}|+| b_{_0}| \ne 0.
$
\hfill (4.5)
\\[3.25ex]
\indent
{\bf Example 4.2.} 
Let us consider the autonomous linear system 
\\[2ex]
\mbox{}\hfill
$
\displaystyle
\dfrac{dx}{dt}  =  a_{_0}  +  a_{1}^{}x  +  a_{2}^{}y
\equiv X(x,y),
\qquad
\dfrac{dy}{dt} =  b_{_0}  +  b_{1}^{}x  +  b_{2}^{}y
\equiv Y(x,y),
$
\hfill (4.6)
\\[2ex]
where $|a_{1}^{}| + |a_{2}^{}| +|b_1^{}| + |b_2^{}|\ne 0.$
\vspace{1ex}

If 
\vspace{0.75ex}
$
|a_1^{} - b_2^{}|  +| a_2^{} + b_1^{}|  \ne  0,
\
|a_{1}^{}| + |a_{2}^{}| +|b_1^{}| + |b_2^{}|\ne 0,
$
then
the stereographically conjugate system to the autonomous linear system (4.6) is the differential system
\\[2ex]
$
\displaystyle
\dfrac{du}{d\tau}  =
{} - a_1^{}u^3 -  (a_2^{} + 2b_1^{})u^2\;\!v
+  (a_1^{} - 2b_2^{})uv^2 +  a_2^{}v^3 \ -
\hfill 
$
\\[2ex]
\mbox{}\hspace{13.35em} 
$
-\ \dfrac{1}{4}\; a_{_0}u^4  -  2b_{_0}u^3\;\!v  -
2b_{_0}uv^3  +  \dfrac{1}{4}\; a_{_0}v^4\equiv U_{_0} (u,v),
$
\hfill (4.7)
\\[1.5ex]
$
\displaystyle
\dfrac{dv}{d\tau}  =
b_1^{}u^3  -  (2a_1^{} - b_2^{})u^2\;\!v  - (2a_2^{} + b_1^{})uv^2  -  b_2^{}v^3   +
b_{_0}u^4  -  2a_{_0}u^3\;\!v  - 2a_{_0}uv^3 +  b_{_0}v^4\equiv V_{_0} (u,v),
\hfill
$
\\[2.25ex]
where  $(u^2 + v^2)\,d\tau  =  dt.$
\vspace{0.75ex}

If 
\vspace{0.75ex}
$
b_2^{} = a_1^{},
\ 
b_1^{} ={} - a_2^{},
\
|a_{1}^{}| + |a_{2}^{}|\ne 0,
$
then
the stereographically conjugate system to the autonomous linear system (4.6) is the differential system
\\[2ex]
\mbox{}\hfill
$
\displaystyle
\dfrac{du}{dt}  =
{} - a_1^{}u  +  a_2^{}v  -  \dfrac{a_{_0}}{4}\; u^2  -
\dfrac{b_{_0}}{2}\; uv   +  \dfrac{a_{_0}}{4}\; v^2\equiv U_{1}^{} (u,v),
\hfill
$
\\[0.25ex]
\mbox{}\hfill                   (4.8)
\\[0.25ex]
\mbox{}\hfill
$
\displaystyle
\dfrac{dv}{dt}  =
{} -  a_2^{}u  -  a_1^{}v  +   \dfrac{b_{_0}}{4}\; u^2  -
\dfrac{a_{_0}}{2}\; uv   -  \dfrac{b_{_0}}{4}\; v^2\equiv V_{1}^{} (u,v).
\hfill 
$
\\[3.5ex]
\indent
{\bf Example 4.3.} Consider the autonomous quadratic system 
\\[2ex]
\mbox{}\hfill
$
\displaystyle
\dfrac{dx}{dt}  =  a_{_0} +  a_1^{}x  +  a_2^{}y  +  a_3^{}x^2
 +  a_4^{}xy  +  a_5^{}y^2\equiv X(x,y),
\hfill
$
\\
\mbox{}\hfill                   (4.9)
\\
\mbox{}\hfill
$
\displaystyle
\dfrac{dy}{dt}  =  b_{_0} +  b_1^{}x  +  b_2^{}y  +
b_3^{}x^2  +  b_4^{}xy  +  b_5^{}y^2\equiv Y(x,y),
\;\ \hfill 
$
\\[2.75ex]
where  $|a_{3}^{}| + |a_{4}^{}| + |a_{5}^{}|+|b_3^{}| + |b_4^{}| + |b_{5}^{}|\ne 0.$
\vspace{1.25ex}

If
\vspace{1ex}
$|a_5^{} - a_3^{} + b_4^{}|  + |a_4^{} + b_3^{} - b_5^{}| \ne  0,\ 
|a_{3}^{}| + |a_{4}^{}| + |a_{5}^{}|+|b_3^{}| + |b_4^{}| + |b_{5}^{}|\ne 0,$
then the stereographically conjugate system to the system (4.9) is the differential system  
\\[2.25ex]
\mbox{}\qquad
$
\dfrac{du}{d\tau}  =
{}- 4a_3^{}u^4 -  4(a_4^{} + 2b_3^{})u^3\;\!v  + 4(a_3^{} - 2b_4^{} - a_5^{})u^2\;\!v^2   +
 4(a_4^{} - 2b_5^{})uv^3 +4a_5^{}v^4
\, -
\hfill 
$
\\[2ex]
\mbox{}\qquad\qquad\quad
$
-\, 
a_1^{}u^5  +  (2b_1^{} - a_2^{})u^4\;\!v  -
2b_2^{}u^3\;\!v^2  -
2b_1^{}u^2\;\!v^3  +  (a_1^{} - 2b_2^{})uv^4  +
 a_2^{}v^5
\,- 
$
\hfill (4.10)
\\[2.5ex]
\mbox{}\qquad\qquad\quad
$
-\,  \dfrac{a_{_0}}{4}\; u^6  -
 \dfrac{b_{_0}}{2}\; u^5\;\!v  -
\dfrac{a_{_0}}{4}\; u^4\;\!v^2  -
b_{_0}u^3\;\!v^3  +  \dfrac{a_{_0}}{4}\; u^2\;\!v^4  -
\dfrac{b_{_0}}{2}\; uv^5  +  \dfrac{a_{_0}}{4}\; v^6 \equiv U_{_0} (u,v),
\hfill
$
\\[1.5ex]
\mbox{}\qquad
$
\dfrac{dv}{d\tau}  =
4b_3^{}u^4  +  4(b_4^{} - 2a_3^{})u^3\;\!v  -
4(b_3^{}+ 2a_4^{} + b_5^{})u^2\;\!v^2   -
4(b_4^{} + 2a_5^{})uv^3   - 4b_5^{}v^4 
\,+ 
\hfill
$
\\[2.25ex]
\mbox{}\qquad\qquad\quad
$
+\,    b_1^{}u^5
 +  (b_2^{} - 2a_1^{})u^4\;\!v  -
 2a_2^{}u^3\;\!v^2  -2a_1^{}u^2\;\!v^3  -   (b_1^{} + 2a_2^{})uv^4   
- b_2^{}v^5 
\,+
 \hfill
$
\\[2.5ex]
\mbox{}\qquad\qquad\quad
$
 +\,  \dfrac{b_{_0}}{4}\; u^6   -
\dfrac{a_{_0}}{2}\; u^5\;\!v  +
\dfrac{b_{_0}}{4}\; u^4\;\!v^2  -
a_{_0}u^3\;\!v^3  -  \dfrac{b_{_0}}{4}\;u^2\;\!v^4  -
\dfrac{a_{_0}}{2}\; uv^5  -  \dfrac{b_{_0}}{4}\; v^6\equiv V_{_0} (u,v),
\hfill
$
\\[2.25ex]
where $(u^2 + v^2)^2\,d\tau  =  dt.$
\vspace{0.75ex}

If
\vspace{0.75ex}
$a_4^{}  =  b_5^{}  -  b_3^{}, \  a_5^{}  =   a_3^{} -  b_4^{}, \
|b_2^{} - a_1^{}| +|b_1^{} + a_2^{}|+
|2a_3^{}- b_4^{}| +  |b_3^{} + b_5^{}| \ne 0,
$ 
and
$
|a_{3}^{}| + | b_5^{}  -  b_3^{}| + |a_3^{} -  b_4^{}|+|b_3^{}| + |b_4^{}| + |b_{5}^{}|\ne 0,$
\vspace{0.75ex}
then the stereographically conjugate system to the system (4.9) is the differential system  
\\[2ex]
\mbox{}\qquad\quad
$
\dfrac{du}{d\theta}  =
{}- 4a_3^{}u^2  -  4(b_3^{} + b_5^{})uv  +
4(a_3^{} - b_4^{})v^2  -  a_1^{}u^3  -
(a_2^{} + 2b_1^{})u^2\;\!v
\, +
\hfill
$
\\[1.75ex]
\mbox{}\qquad\qquad\qquad
$
+ \,
(a_1^{} - 2b_2^{})uv^2  + 2a_2^{}v^3   -
\dfrac{a_{_0}}{4}\; u^4  -  \dfrac{b_{_0}}{2}\; u^3\;\!v
-  \dfrac{b_{_0}}{2}\; uv^3  +  \dfrac{a_{_0} }{4}\; v^4\equiv U_{1}^{} (u,v),
\hfill
$
\\[1ex]
\mbox{}\hfill (4.11)
\\[-0.15ex]
\mbox{}\qquad\quad
$
\dfrac{dv}{d\theta}  =
{} -  4b_3^{}u^2  -
4(b_4^{} - 2a_3^{})uv  + 4b_5^{}v^2  +  b_1^{}u^3  +  (b_2^{} - 2a_1^{})u^2\;\!v 
\, -
\hfill
$
\\[1.75ex]
\mbox{}\qquad\qquad\qquad
$
- \, 
(b_1^{} + 2a_2^{})uv^2  -   b_2^{}v^3  +
\dfrac{b_{_0}}{4}\; u^4  -  \dfrac{a_{_0}}{2}\; u^3\;\!v  -
\dfrac{a_{_0}}{2}\; uv^3  -  \dfrac{b_{_0}}{4}\; v^4\equiv V_{1}^{} (u,v),
\hfill
$
\\[2.25ex]
where
\vspace{0.5ex}
$(u^2 + v^2)\,d\theta =  dt.$ 

If
\vspace{0.75ex}
$a_4^{} = 2b_5^{}, \ a_5^{} = {}-a_3^{}, \ b_3^{} = {}- b_5^{},
\ b_4^{} = 2a_3^{}, \
 b_2^{} = a_1^{}, \ b_1^{} =  {}- a_1^{},\  
|a_{3}^{}| + | b_5^{}|\ne 0,$
then the stereographically conjugate system to the system (4.9) is the differential system  
\\[2.5ex]
\mbox{}\hfill
$
\dfrac{du}{dt}  = {} - 4a_3^{}  - a_1^{}u  +  a_2^{}v  -  \dfrac{a_{_0}}{4}\; u^2  -
\dfrac{b_{_0}}{2}\; uv  +  \dfrac{a_{_0}}{4}\; v^2\equiv U_{2}^{} (u,v),
\hfill
$
\\[0.5ex]
\mbox{}\hfill (4.12)
\\[0.5ex]
\mbox{}\hfill
$
\dfrac{dv}{dt}  =
4b_5^{}  - a_2^{}u  -  a_1^{}v  +  \dfrac{b_{_0}}{4}\;u^2  -
\dfrac{a_{_0}}{2}\; uv  -  \dfrac{b_{_0}}{4}\; v^2\equiv V_{2}^{} (u,v).
\quad\;\, 
\hfill
$
\\[4.25ex]
\centerline{\bf  5. Stereographic atlas of trajectories for differential system}
\\[1.25ex]
\indent
Using the plane 
\\[1ex]
\mbox{}\hfill
$
\{(x^{\star},y^{\star},z^{\star})\colon z^{\star}=1-\varepsilon_1^{}\},
\ \ \
0<\varepsilon_1^{}< 1,
\hfill
$ 
\\[1.5ex]
we divide the sphere (1.1) on two parts and take the part
\\[1.5ex]
\mbox{}\hfill
$
S_1^2=\{(x^{\star},y^{\star},z^{\star})\colon\,
x^{\star}{}^{\,^{\scriptstyle 2}}+
y^{\star}{}^{\,^{\scriptstyle 2}}+
z^{\star}{}^{\,^{\scriptstyle 2}}=1, \ \,
{}-1 \leq z^{\star} \leq 1-\varepsilon_1^{}\}
\hfill
$
\\[1.75ex]
without the northern pole $N(0,0,1).$
\vspace{0.35ex}
Number $\varepsilon_1^{}\in (0;1)$ such that 
stereographic images of all equilibrium states and the isolated closed trajectories of system (D), 
\vspace{0.35ex}
which are lying in the final part of the extended phase plane $\overline{Oxy},$  
\vspace{0.5ex}
are located on the part $S_1^2$ of sphere (1.1).

Let the circle $K (x, y)$ be the circle lying on the phase plane $Oxy$
\vspace{0.35ex}
with the centre in the origin of coordinates $O$ and being pre-image of the part $S_1^2$ of sphere  (1.1) 
\vspace{0.35ex}
by the stereographic projection with the projection centre in the northern pole $N(0,0,1)$ 
\vspace{0.5ex}
(Fig. 5.1).

Using the plane 
\\[1ex]
\mbox{}\hfill
$
\{(x^{\star},y^{\star},z^{\star})\colon\, 
z^{\star}=\varepsilon_2^{}-1\},
\ \ \
0<\varepsilon_2^{}< 1,
\hfill
$ 
\\[1.5ex]
we divide the sphere (1.1) on two parts and take the part
\\[1.5ex]
\mbox{}\hfill
$
S_2^2=\{(x^{\star},y^{\star},z^{\star})\colon\,
x^{\star}{}^{\,^{\scriptstyle 2}}+
y^{\star}{}^{\,^{\scriptstyle 2}}+
z^{\star}{}^{\,^{\scriptstyle 2}}=1, \ \,
\varepsilon_2^{}-1 \leq z^{\star} \leq 1\},
\hfill
$
\\[1.75ex]
without the southern pole $S(0,0,{}-1).$
\vspace{0.35ex}
Number $\varepsilon_2^{}\in (0;1)$ such that 
stereographic images of all equilibrium states and the isolated closed trajectories of system (3.3), 
\vspace{0.5ex}
which are lying in the final part of the extended phase plane $\overline{O^{\ast} uv},$  
are located on the part $S_2^2$ of sphere (1.1).
\\[4ex]
\mbox{}\hfill
{\unitlength=1mm
\begin{picture}(106,57)
\put(0,0){\includegraphics[width=106mm,height=57mm]{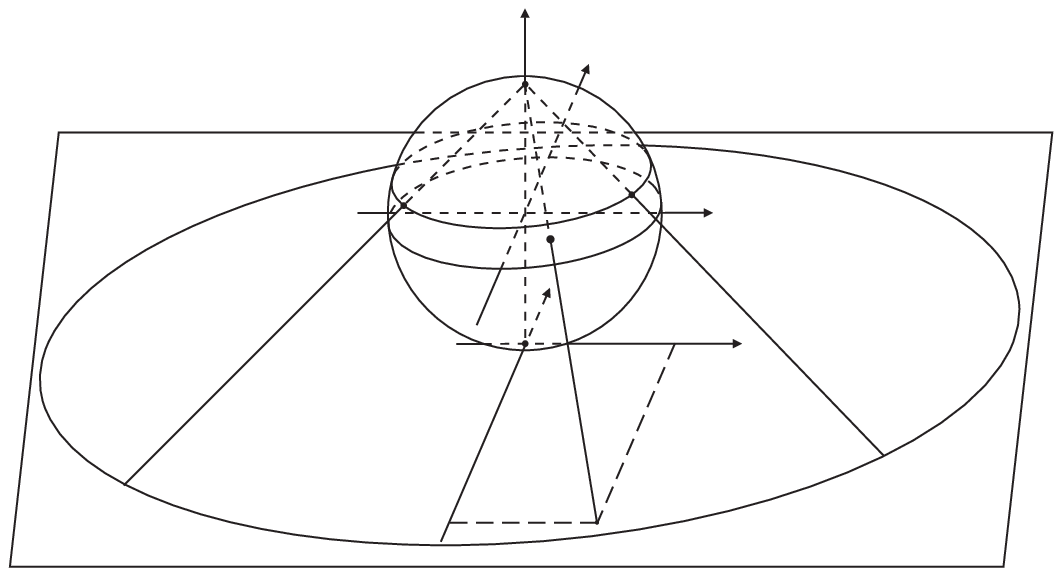}}

\put(53.9,20.7){\makebox(0,0)[cc]{\scriptsize $O$}}
\put(50.7,24.7){\makebox(0,0)[cc]{\scriptsize $S$}}
\put(68.5,21.5){\makebox(0,0)[cc]{\scriptsize $x$}}
\put(73,21){\makebox(0,0)[cc]{\scriptsize $x$}}
\put(43,4.9){\makebox(0,0)[cc]{\scriptsize $y$}}
\put(53.6,29){\makebox(0,0)[cc]{\scriptsize $y$}}
\put(63,6.6){\makebox(0,0)[cc]{\scriptsize $M$}}

\put(50.5,38){\makebox(0,0)[cc]{\scriptsize $O^{\star}$}}
\put(56.9,32.9){\makebox(0,0)[cc]{\scriptsize $P$}}
\put(70.6,34.1){\makebox(0,0)[cc]{\scriptsize $x^{\star}$}}
\put(61,50.3){\makebox(0,0)[cc]{\scriptsize $y^{\star}$}}

\put(50.6,51.5){\makebox(0,0)[cc]{\scriptsize $N$}}
\put(54.7,56.5){\makebox(0,0)[cc]{\scriptsize $z^{\star}$}}

\put(52.5,-7){\makebox(0,0)[cc]{\rm Fig. 5.1}}
\end{picture}}
\hfill\mbox{}
\\[7.75ex]
\indent
Let the circle $K (u, v)$ be the circle lying on the phase plane $O^{\ast} uv$
\vspace{0.35ex}
with the centre in the origin of coordinates $O^{\ast}$ and being pre-image of the part $S_2^2$ of sphere  (1.1) 
\vspace{0.35ex}
by the stereographic projection of the plane $O^{\ast} uv$ 
\vspace{0.35ex}
with the projection centre in the 
southern pole $S(0,0,{}-1)$ (Fig. 5.2).
\\[6ex]
\mbox{}\hfill
{\unitlength=1mm
\begin{picture}(80,58)
\put(0,0){\includegraphics[width=79.79mm,height=58.58mm]{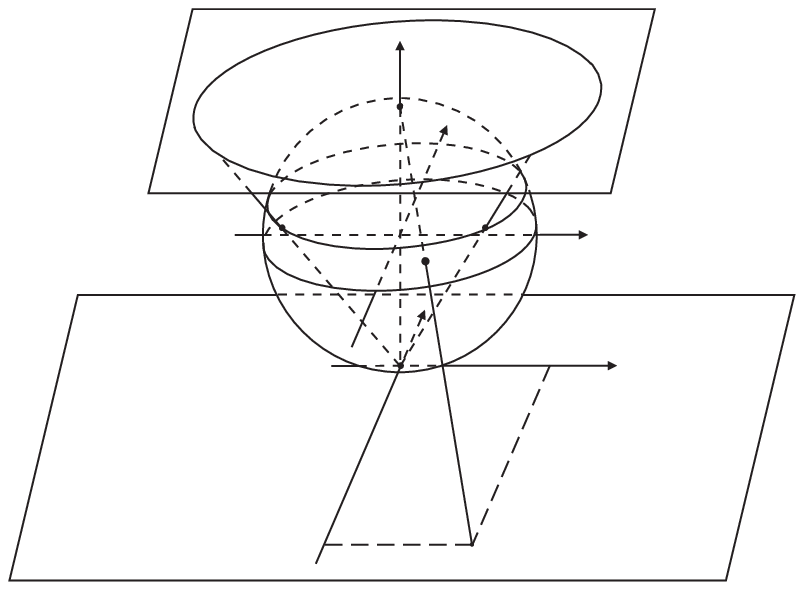}}

\put(41,19.7){\makebox(0,0)[cc]{\scriptsize $O$}}
\put(38.4,23.7){\makebox(0,0)[cc]{\scriptsize $S$}}
\put(56,20.5){\makebox(0,0)[cc]{\scriptsize $x$}}
\put(60.7,20){\makebox(0,0)[cc]{\scriptsize $x$}}
\put(30.8,3.9){\makebox(0,0)[cc]{\scriptsize $y$}}
\put(40.7,27.5){\makebox(0,0)[cc]{\scriptsize $y$}}
\put(49,2.6){\makebox(0,0)[cc]{\scriptsize $M$}}

\put(37.4,37.1){\makebox(0,0)[cc]{\scriptsize $O^{\star}$}}
\put(44.9,32.5){\makebox(0,0)[cc]{\scriptsize $P$}}
\put(58.6,33.8){\makebox(0,0)[cc]{\scriptsize $x^{\star}$}}
\put(42.5,46.7){\makebox(0,0)[cc]{\scriptsize $y^{\star}$}}

\put(37.6,47){\makebox(0,0)[cc]{\scriptsize $N$}}
\put(42,54.5){\makebox(0,0)[cc]{\scriptsize $z^{\star}$}}

\put(40,-7){\makebox(0,0)[cc]{\rm Fig. 5.2}}
\end{picture}}
\hfill\mbox{}
\\[7.75ex]
\indent
The ordered pair $(K(x,y), K(u,v))$ of the circles $K(x,y)$ and $K(u,v)$ 
\vspace{0.35ex}
with trajectories of systems (D) and (3.3) plotted on them 
\vspace{0.25ex}
is called {\it stereographic atlas of trajectories} for system (D).
Then (by Theorem 2.1) 
\vspace{0.25ex}
the ordered pair $(K(u,v), K(x,y))$ is stereographic atlas of trajectories for system (3.3).
\vspace{0.25ex}

Correspondence between the circles $K (x, y) $ and $K (u, v) $ is shown on Fig. 5.3.
\vspace{0.25ex}

Using the numbers $1, \ldots, 40,$ 
\vspace{0.15ex}
we reflect correspondences between halfneighbourhoods of points, 
\vspace{0.15ex}
which are lying on the coordinate axes and concentric circles with the centre in the origin of coordinates.

\newpage

\mbox{}
\\[-1.75ex]
\mbox{}\hfill
{\unitlength=1mm
\begin{picture}(65,65)
\put(0,0){\includegraphics[width=65mm,height=65mm]{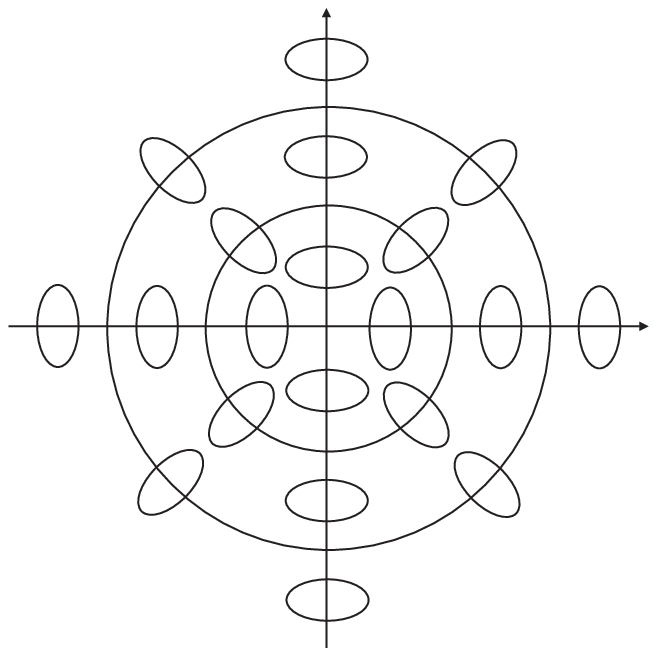}}

\put(38.9,34.5){\makebox(0,0)[cc]{\scriptsize $1$}}
\put(38.9,30.8){\makebox(0,0)[cc]{\scriptsize $2$}}
\put(49.9,34.5){\makebox(0,0)[cc]{\scriptsize $3$}}
\put(49.9,30.8){\makebox(0,0)[cc]{\scriptsize $4$}}
\put(59.9,34.5){\makebox(0,0)[cc]{\scriptsize $5$}}
\put(59.9,30.8){\makebox(0,0)[cc]{\scriptsize $6$}}
\put(39.9,39.9){\makebox(0,0)[cc]{\scriptsize $7$}}
\put(42.9,42.8){\makebox(0,0)[cc]{\scriptsize $8$}}
\put(46.9,46.9){\makebox(0,0)[cc]{\scriptsize $9$}}
\put(49.5,49.5){\makebox(0,0)[cc]{\scriptsize $10$}}
\put(34.2,38.7){\makebox(0,0)[cc]{\scriptsize $11$}}
\put(30.5,38.7){\makebox(0,0)[cc]{\scriptsize $12$}}
\put(34.2,49.9){\makebox(0,0)[cc]{\scriptsize $13$}}
\put(30.5,49.9){\makebox(0,0)[cc]{\scriptsize $14$}}
\put(34.2,59.7){\makebox(0,0)[cc]{\scriptsize $15$}}
\put(30.5,59.7){\makebox(0,0)[cc]{\scriptsize $16$}}
\put(25.1,39.9){\makebox(0,0)[cc]{\scriptsize $17$}}
\put(22.7,42.8){\makebox(0,0)[cc]{\scriptsize $18$}}
\put(18.1,47.1){\makebox(0,0)[cc]{\scriptsize $19$}}
\put(15.5,49.9){\makebox(0,0)[cc]{\scriptsize $20$}}
\put(26.3,34.5){\makebox(0,0)[cc]{\scriptsize $21$}}
\put(26.3,31){\makebox(0,0)[cc]{\scriptsize $22$}}
\put(15.2,34.5){\makebox(0,0)[cc]{\scriptsize $23$}}
\put(15.2,31){\makebox(0,0)[cc]{\scriptsize $24$}}
\put(5.2,34.5){\makebox(0,0)[cc]{\scriptsize $25$}}
\put(5.2 ,31){\makebox(0,0)[cc]{\scriptsize $26$}}
\put(25.1,24.9){\makebox(0,0)[cc]{\scriptsize $27$}}
\put(22.5,22.4){\makebox(0,0)[cc]{\scriptsize $28$}}
\put(17.8,18.2){\makebox(0,0)[cc]{\scriptsize $29$}}
\put(15.4,15.5){\makebox(0,0)[cc]{\scriptsize $30$}}
\put(34.2,26.3){\makebox(0,0)[cc]{\scriptsize $32$}}
\put(30.5,26.3){\makebox(0,0)[cc]{\scriptsize $31$}}
\put(34.2,15.1){\makebox(0,0)[cc]{\scriptsize $34$}}
\put(30.5,15.1){\makebox(0,0)[cc]{\scriptsize $33$}}
\put(34.2,5){\makebox(0,0)[cc]{\scriptsize $36$}}
\put(30.5,5){\makebox(0,0)[cc]{\scriptsize $35$}}
\put(40.1,25.1){\makebox(0,0)[cc]{\scriptsize $37$}}
\put(42.7,22.5){\makebox(0,0)[cc]{\scriptsize $38$}}
\put(47.2,18.1){\makebox(0,0)[cc]{\scriptsize $39$}}
\put(49.6,15.4){\makebox(0,0)[cc]{\scriptsize $40$}}
\put(64,31){\makebox(0,0)[cc]{\scriptsize $x$}}
\put(30.5,64){\makebox(0,0)[cc]{\scriptsize $y$}}

\end{picture}
}
\quad\qquad
{\unitlength=1mm
\begin{picture}(65,65)
\put(0,0){\includegraphics[width=65mm,height=65mm]{r05-03.eps}}

\put(38.7,34.5){\makebox(0,0)[cc]{\scriptsize $5$}}
\put(38.7,30.8){\makebox(0,0)[cc]{\scriptsize $6$}}
\put(49.9,34.5){\makebox(0,0)[cc]{\scriptsize $3$}}
\put(49.9,30.8){\makebox(0,0)[cc]{\scriptsize $4$}}
\put(59.9,34.5){\makebox(0,0)[cc]{\scriptsize $1$}}
\put(59.9,30.8){\makebox(0,0)[cc]{\scriptsize $2$}}
\put(40,40){\makebox(0,0)[cc]{\scriptsize $10$}}
\put(42.9,42.8){\makebox(0,0)[cc]{\scriptsize $9$}}
\put(46.9,46.9){\makebox(0,0)[cc]{\scriptsize $8$}}
\put(49.5,49.5){\makebox(0,0)[cc]{\scriptsize $7$}}
\put(34.2,38.7){\makebox(0,0)[cc]{\scriptsize $15$}}
\put(30.5,38.7){\makebox(0,0)[cc]{\scriptsize $16$}}
\put(34.2,49.9){\makebox(0,0)[cc]{\scriptsize $13$}}
\put(30.5,49.9){\makebox(0,0)[cc]{\scriptsize $14$}}
\put(34.2,59.7){\makebox(0,0)[cc]{\scriptsize $11$}}
\put(30.5,59.7){\makebox(0,0)[cc]{\scriptsize $12$}}
\put(25.2,40){\makebox(0,0)[cc]{\scriptsize $20$}}
\put(22.7,42.8){\makebox(0,0)[cc]{\scriptsize $19$}}
\put(18.1,47.1){\makebox(0,0)[cc]{\scriptsize $18$}}
\put(15.6,49.7){\makebox(0,0)[cc]{\scriptsize $17$}}
\put(26.3,34.5){\makebox(0,0)[cc]{\scriptsize $25$}}
\put(26.3,31){\makebox(0,0)[cc]{\scriptsize $26$}}
\put(15.2,34.5){\makebox(0,0)[cc]{\scriptsize $23$}}
\put(15.2,31){\makebox(0,0)[cc]{\scriptsize $24$}}
\put(5.2,34.5){\makebox(0,0)[cc]{\scriptsize $21$}}
\put(5.2 ,31){\makebox(0,0)[cc]{\scriptsize $22$}}
\put(25,25.1){\makebox(0,0)[cc]{\scriptsize $30$}}
\put(22.5,22.4){\makebox(0,0)[cc]{\scriptsize $29$}}
\put(17.8,18.2){\makebox(0,0)[cc]{\scriptsize $28$}}
\put(15.4,15.4){\makebox(0,0)[cc]{\scriptsize $27$}}
\put(34.2,26.3){\makebox(0,0)[cc]{\scriptsize $36$}}
\put(30.5,26.3){\makebox(0,0)[cc]{\scriptsize $35$}}
\put(34.2,15.1){\makebox(0,0)[cc]{\scriptsize $34$}}
\put(30.5,15.1){\makebox(0,0)[cc]{\scriptsize $33$}}
\put(34.2,5){\makebox(0,0)[cc]{\scriptsize $32$}}
\put(30.5,5){\makebox(0,0)[cc]{\scriptsize $31$}}
\put(40.1,25.1){\makebox(0,0)[cc]{\scriptsize $40$}}
\put(42.7,22.5){\makebox(0,0)[cc]{\scriptsize $39$}}
\put(47.2,18.1){\makebox(0,0)[cc]{\scriptsize $38$}}
\put(49.7,15.4){\makebox(0,0)[cc]{\scriptsize $37$}}
\put(64,31){\makebox(0,0)[cc]{\scriptsize $u$}}
\put(30.5,64){\makebox(0,0)[cc]{\scriptsize $v$}}

\end{picture}
}
\hfill\mbox{}
\\[0.5ex]
\mbox{}\hfill
Fig. 5.3
\hfill\mbox{}
\\[3ex]
\indent
{\bf Example 5.1.}
Trajectories of the differential system (4.4) are the parallel straight lines 
\\[1.5ex]
\mbox{}\hfill
$
b_{_0}x-a_{_0}y = C,
\quad 
C\in\R.
\hfill
$
\\[1.5ex]
\indent
Trajectories of the differential system (4.5) are the equilibrium state $O^{\ast} (0,0)$ and 
adjoining curves to this equilibrium state 
\\[1.75ex]
\mbox{}\hfill
$
\dfrac{b_{_0}u-a_{_0}v}{u^2+v^2}=C^{\ast},
\quad 
|u|+|v|\ne 0,
\quad 
C^{\ast}\in\R,
\quad 
C^{\ast}=4C.
\hfill
$
\\[2ex]
\indent
The behaviour of trajectories on the sphere (1.1) 
\vspace{0.35ex}
for the stereographically mutually conjugate systems (4.4) and (4.5) 
is represented on Fig. 5.4 for $a_{_0}=1,\ b_{_0}=0$.
\vspace{0.5ex}

The circles on Fig. 5.5 are stereografic atlases of trajectories for systems (4.4) and (4.5).
\\[4ex]
\parbox{50mm}{
\mbox{}\hfill
{\unitlength=1mm
\begin{picture}(42,42)
\put(0,0){\includegraphics[width=42mm,height=42mm]{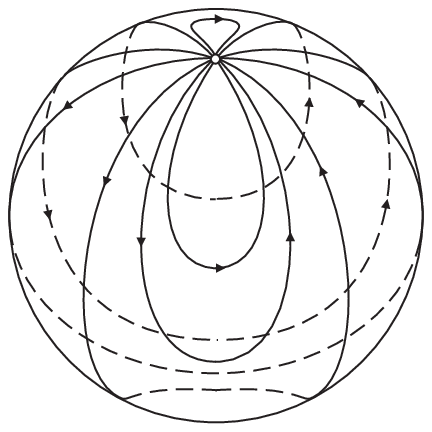}}
\end{picture}}
\hfill\mbox{}
\\[2ex]
\mbox{}\hfill
Fig. 5.4
\hfill\mbox{}
}
\parbox{104mm}{
\mbox{}\hfill
{\unitlength=1mm
\begin{picture}(42,42)
\put(0,0){\includegraphics[width=42mm,height=42mm]{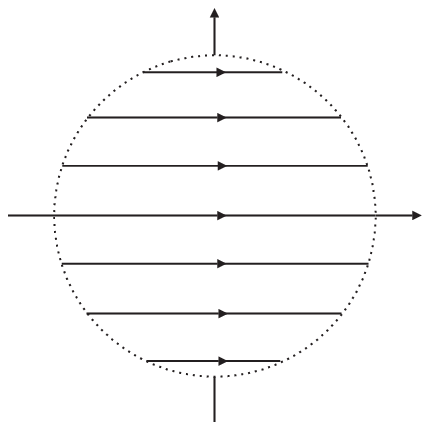}}
\put(41,19){\makebox(0,0)[cc]{\scriptsize $x$}}
\put(19,41){\makebox(0,0)[cc]{\scriptsize $y$}}
\end{picture}}
\qquad
{\unitlength=1mm
\begin{picture}(42,42)
\put(0,0){\includegraphics[width=42mm,height=42mm]{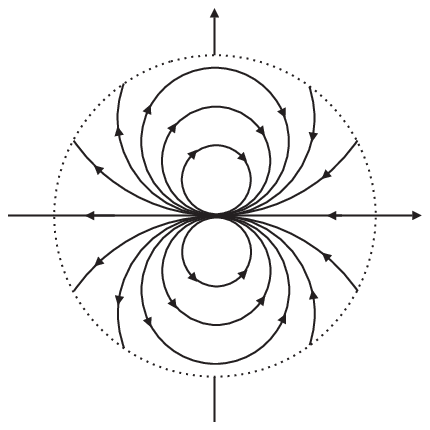}}
\put(41,19){\makebox(0,0)[cc]{\scriptsize $u$}}
\put(19,41){\makebox(0,0)[cc]{\scriptsize $v$}}
\end{picture}}
\hfill\mbox{}
\\[2ex]
\mbox{}\hfill
Fig. 5.5
\hfill\mbox{}
}
\\[4ex]
\indent
Notice that the  equilibrium state  $O^{\ast} (0,0)$ 
\vspace{0.25ex}
of the homogeneous quadratic system (4.5) is complicated and 
consists of two elliptic Bendixon's sectors, which are limited by  the trajectories-rays of the straight line $v=0.$
\vspace{0.75ex}

{\bf Example 5.2.}
Trajectories of the linear system
\\[2ex]
\mbox{}\hfill                        
$
\dfrac{dx}{dt}=x,
\qquad
\dfrac{dy}{dt}=y
$
\hfill (5.1)
\\[2.25ex]
are $O\!$-rays of the family of straight lines 
\\[1.5ex]
\mbox{}\hfill
$
C_1^{}y+C_2^{}x=0,
\quad 
C_1^{},C_2^{}\in \R,
\hfill
$
\\[1.5ex]
and the equilibrium state (unstable dicritical node) $O (0,0).$

\newpage

Trajectories of  the stereographically conjugate system
\\[2ex]
\mbox{}\hfill                        
$
\dfrac{du}{dt}={}-u,
\qquad
\dfrac{dv}{dt}={}-v
$
\hfill (5.2)
\\[2.25ex]
are $O\!$-rays of the family of straight lines 
\\[1.5ex]
\mbox{}\hfill
$
C_1^{}v+C_2^{}u=0,
\quad 
C_1^{},C_2^{}\in \R,
\hfill
$
\\[1.5ex]
and  the equilibrium state (stable dicritical node) $O^{\ast} (0,0).$
\\[4ex]
\parbox{50mm}{
\mbox{}\hfill
{\unitlength=1mm
\begin{picture}(42,42)
\put(0,0){\includegraphics[width=42mm,height=42mm]{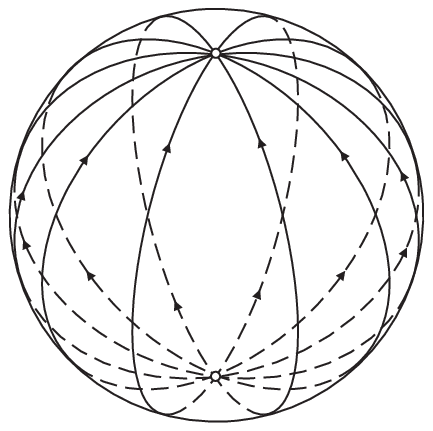}}
\end{picture}}
\hfill\mbox{}
\\[2ex]
\mbox{}\hfill
Fig. 5.6
\hfill\mbox{}
}
\parbox{104mm}{
\mbox{}\hfill
{\unitlength=1mm
\begin{picture}(42,42)
\put(0,0){\includegraphics[width=42mm,height=42mm]{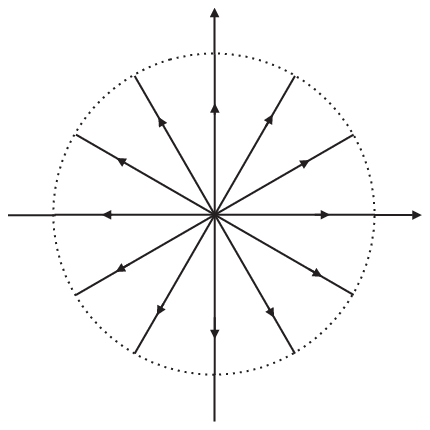}}
\put(41,19){\makebox(0,0)[cc]{\scriptsize $x$}}
\put(19,41){\makebox(0,0)[cc]{\scriptsize $y$}}
\end{picture}}
\qquad
{\unitlength=1mm
\begin{picture}(42,42)
\put(0,0){\includegraphics[width=42mm,height=42mm]{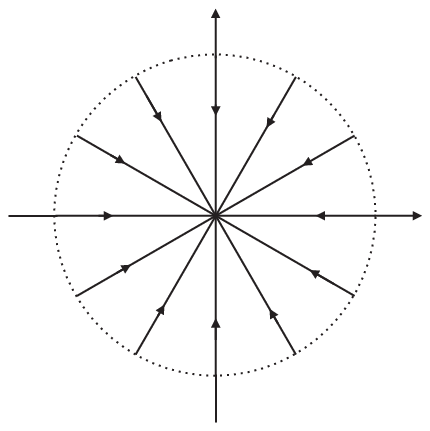}}
\put(41,19){\makebox(0,0)[cc]{\scriptsize $u$}}
\put(19,41){\makebox(0,0)[cc]{\scriptsize $v$}}
\end{picture}}
\hfill\mbox{}
\\[2ex]
\mbox{}\hfill
Fig. 5.7
\hfill\mbox{}
}
\\[3.5ex]
\indent
The behaviour of trajectories on the sphere (1.1) for the stereographically mutually conjugate systems (5.1) and (5.2)
is represented on Fig. 5.6. 
The stereographic images of trajectories for the differential systems (5.1) and (5.2) on the sphere (1.1) are 
semicircles of meridians of sphere (1.1), which are adjoining to the northern and to the southern poles, 
and also the northern and the southern poles of sphere (1.1).
\vspace{0.25ex}

The circles on Fig. 5.7 are  
\vspace{1ex}
stereographic atlases of trajectories for systems\! (5.1)\! and\! (5.2).

{\bf Example 5.3.}
Trajectories of the linear system
\\[2ex]
\mbox{}\hfill                        
$
\dfrac{dx}{dt}=y,
\qquad
\dfrac{dy}{dt}={}-x
$
\hfill (5.3)
\\[2.25ex]
are the concentric circles 
\\[1.5ex]
\mbox{}\hfill
$
x^2+y^2=C,
\quad 
C\in (0; {}+\infty),
\hfill
$
\\[1.5ex]
and  the equilibrium state (centre) $O(0,0).$
\\[4ex]
\parbox{50mm}{
\mbox{}\hfill
{\unitlength=1mm
\begin{picture}(42,42)
\put(0,0){\includegraphics[width=42mm,height=42mm]{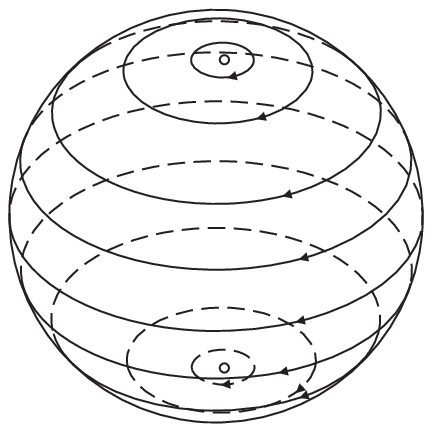}}
\end{picture}}
\hfill\mbox{}
\\[2ex]
\mbox{}\hfill
Fig. 5.8
\hfill\mbox{}
}
\parbox{104mm}{
\mbox{}\hfill
{\unitlength=1mm
\begin{picture}(42,42)
\put(0,0){\includegraphics[width=42mm,height=42mm]{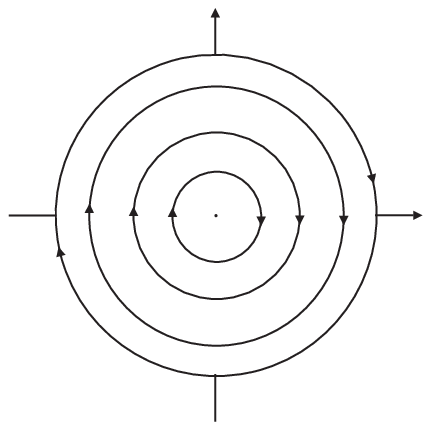}}
\put(41,19){\makebox(0,0)[cc]{\scriptsize $x$}}
\put(19,41){\makebox(0,0)[cc]{\scriptsize $y$}}
\end{picture}}
\qquad
{\unitlength=1mm
\begin{picture}(42,42)
\put(0,0){\includegraphics[width=42mm,height=42mm]{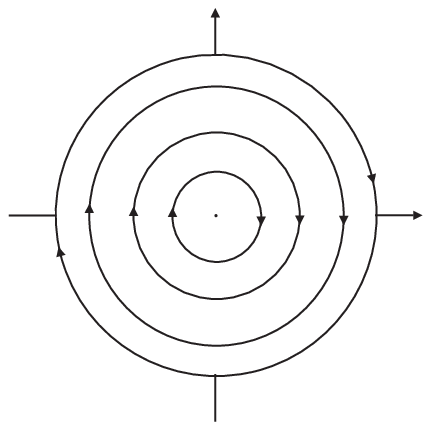}}
\put(41,19){\makebox(0,0)[cc]{\scriptsize $u$}}
\put(19,41){\makebox(0,0)[cc]{\scriptsize $v$}}
\end{picture}}
\hfill\mbox{}
\\[2ex]
\mbox{}\hfill
Fig. 5.9
\hfill\mbox{}
}
\\[3ex]
\indent
The direction of movement along trajectories for the differential system (5.3)
are defined by  the tangential vectors 
\\[1.25ex]
\mbox{}\hfill
$
\vec{a}(x,y)=(y,{}-x)
$ 
\ for all $(x,y)\in \R^2\backslash\{(0,0)\}.
\hfill
$
\\[-1.5ex]

\newpage

Trajectories of  the stereographically conjugate system
\\[2ex]
\mbox{}\hfill                        
$
\dfrac{du}{dt}=v,
\qquad
\dfrac{dv}{dt}={}-u
$
\hfill (5.4)
\\[1.75ex]
are the concentric circles 
\\[1.25ex]
\mbox{}\hfill
$
u^2+v^2=C^\ast,
\quad 
C^\ast \in (0;{}+\infty), 
\ \ C^\ast=4C,
\hfill
$
\\[1.5ex]
and the equilibrium state (centre) $O^{\ast} (0,0).$
\vspace{0.25ex}

The stereographic images 
\vspace{0.15ex}
of trajectories for systems (5.3) and (5.4) on the sphere (1.1) are 
parallels, the northern and the southern poles of sphere (1.1) (Fig. 5.8).
\vspace{0.25ex}

The circles on Fig. 5.9 are 
\vspace{1ex}
stereographic atlases of trajectories for systems\! (5.3)\! and\! (5.4).

{\bf Example 5.4.}
Trajectories of the linear system
\\[2ex]
\mbox{}\hfill                        
$
\dfrac{dx}{dt}=x-y,
\qquad
\dfrac{dy}{dt}=x+y
$
\hfill (5.5)
\\[2ex]
are the logarithmic spirals  
\\[1.75ex]
\mbox{}\hfill
$
(x^2+y^2)\exp \Bigl( {}-\arctan \dfrac{y}{x}\Bigr)=C,
\quad 
C\in (0;{}+\infty),
\hfill
$
\\[1.75ex]
and  the equilibrium state (unstable focus) $O (0,0).$
\vspace{0.25ex}

Trajectories of  the stereographically conjugate system
\\[2ex]
\mbox{}\hfill                        
$
\dfrac{du}{dt}={}-u-v,
\qquad
\dfrac{dv}{dt}=u-v
$
\hfill (5.6)
\\[2ex]
are the logarithmic spirals  
\\[1.75ex]
\mbox{}\hfill
$
(u^2+v^2)\exp \Bigl( {}-\arctan \dfrac{v}{u}\Bigr)=C^\ast,
\quad 
C^\ast\in (0;{}+\infty),
\quad 
C^\ast=4C,
\hfill
$
\\[1.75ex]
and the equilibrium state (stable focus) $O^{\ast} (0,0).$
\vspace{0.25ex}

Trajectories of the stereographically mutually conjugate systems (5.5) and (5.6)
on the sphere (1.1) are represented on Fig. 5.10.
The circles on Fig. 5.11 are stereographic atlases of trajectories for the differential systems (5.5) and (5.6).
\\[4ex]
\parbox{50mm}{
\mbox{}\hfill
{\unitlength=1mm
\begin{picture}(42,42)
\put(0,0){\includegraphics[width=42mm,height=42mm]{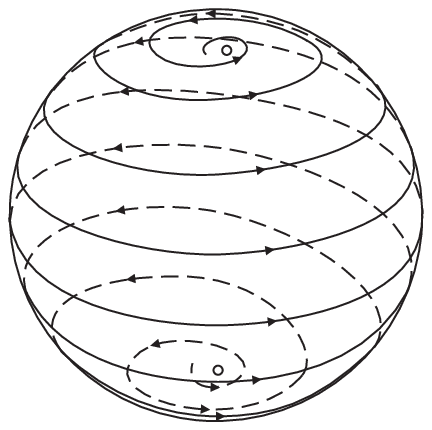}}
\end{picture}}
\hfill\mbox{}
\\[2ex]
\mbox{}\hfill
Fig. 5.10
\hfill\mbox{}
}
\parbox{104mm}{
\mbox{}\hfill
{\unitlength=1mm
\begin{picture}(42,42)
\put(0,0){\includegraphics[width=42mm,height=42mm]{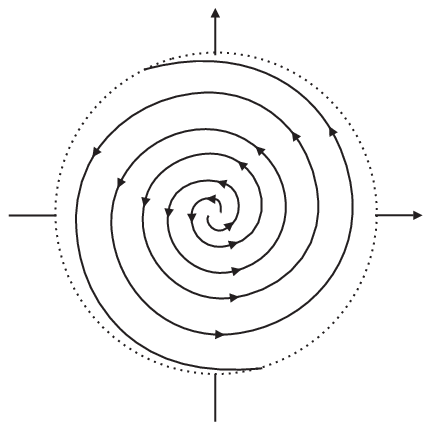}}
\put(41,19){\makebox(0,0)[cc]{\scriptsize $x$}}
\put(19,41){\makebox(0,0)[cc]{\scriptsize $y$}}
\end{picture}}
\qquad
{\unitlength=1mm
\begin{picture}(42,42)
\put(0,0){\includegraphics[width=42mm,height=42mm]{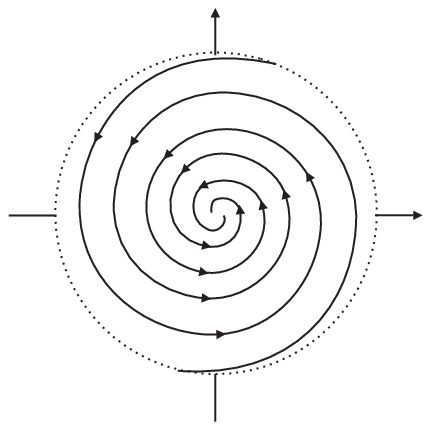}}
\put(41,19){\makebox(0,0)[cc]{\scriptsize $u$}}
\put(19,41){\makebox(0,0)[cc]{\scriptsize $v$}}
\end{picture}}
\hfill\mbox{}
\\[2ex]
\mbox{}\hfill
Fig. 5.11
\hfill\mbox{}
}
\\[4ex]
\indent
{\bf Example 5.5.}
Trajectories of the linear system
\\[2ex]
\mbox{}\hfill                        
$
\dfrac{dx}{dt}=x,
\qquad
\dfrac{dy}{dt}={}-y
$
\hfill (5.7)
\\[2.25ex]
are $O\!$-curves of the family  
\vspace{0.35ex}
$
xy=C,\ 
C\in \R,
$
and the equilibrium state (saddle such that its separatrices are the coordinate $O\!$-rays) $O(0,0).$ 

\newpage

Trajectories of  the stereographically conjugate system
\\[2ex]
\mbox{}\hfill                        
$
\dfrac{du}{d\tau}={}-u^3+3uv^2,
\qquad
\dfrac{dv}{d\tau}={}-3u^2\;\!v+v^3, 
\qquad 
(u^2+v^2)\, d\tau=dt,
$
\hfill (5.8)
\\[2ex]
are $O^{\ast}\!$-curves of the family  
\\[1.75ex]
\mbox{}\hfill
$
\dfrac{uv}{(u^2+v^2)^2 }=C^{\ast},
\quad 
C^{\ast}\in \R,
\quad 
C^{\ast}=16C,
\hfill
$
\\[1.75ex]
and the complicated  equilibrium state  $O^{\ast} (0,0).$
This equilibrium state is consisting of four elliptic Bendixon's sectors, which are 
limited by the coordinate $O^{\ast}\!$-rays and organised by lemniscates of Bernoulli.

The behaviour of trajectories on the sphere (1.1) for the stereographically mutually conjugate systems (5.7) and (5.8)
is represented on Fig. 5.12.
The circles on Fig. 5.13 are stereographic atlases of trajectories for the differential systems (5.7) and (5.8).
\\[3ex]
\mbox{}\hfill
{\unitlength=1mm
\begin{picture}(42,42)
\put(0,0){\includegraphics[width=42mm,height=42mm]{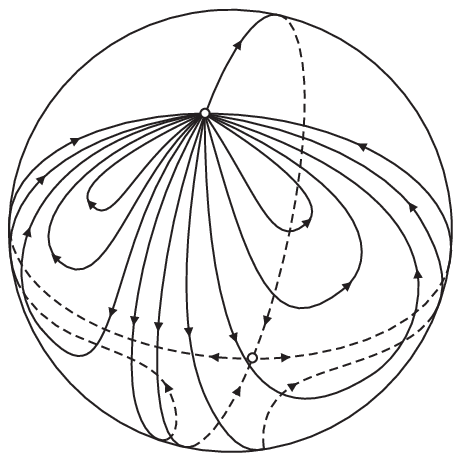}}
\end{picture}}
\qquad
\qquad
\qquad
{\unitlength=1mm
\begin{picture}(42,42)
\put(0,0){\includegraphics[width=42mm,height=42mm]{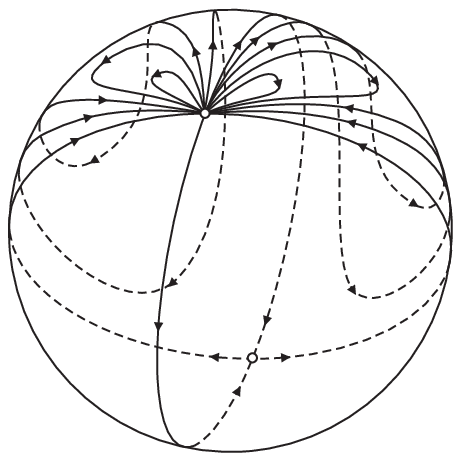}}
\end{picture}}
\hfill\mbox{}
\\[2ex]
\mbox{}\hfill
Fig. 5.12
\hfill\mbox{}
\\[7ex]
\mbox{}\hfill
{\unitlength=1mm
\begin{picture}(42,42)
\put(0,0){\includegraphics[width=42mm,height=42mm]{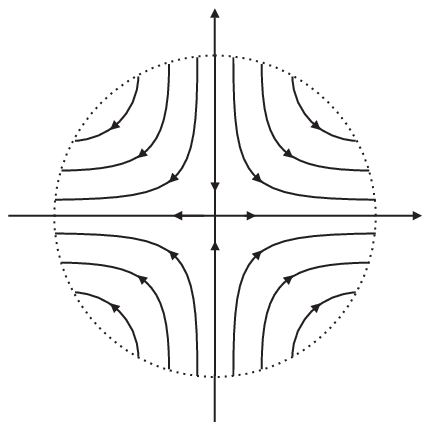}}
\put(41,19){\makebox(0,0)[cc]{\scriptsize $x$}}
\put(19,41){\makebox(0,0)[cc]{\scriptsize $y$}}
\end{picture}}
\qquad
\qquad
{\unitlength=1mm
\begin{picture}(42,42)
\put(0,0){\includegraphics[width=42mm,height=42mm]{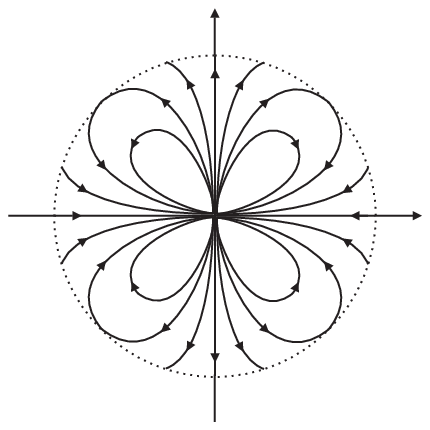}}
\put(41,19){\makebox(0,0)[cc]{\scriptsize $u$}}
\put(19,41){\makebox(0,0)[cc]{\scriptsize $v$}}
\end{picture}}
\hfill\mbox{}
\\[1.5ex]
\mbox{}\hfill
Fig. 5.13
\hfill\mbox{}
\\[3ex]
\indent
{\bf Example 5.6.}
Trajectories of the linear differential system
\\[2ex]
\mbox{}\hfill                        
$
\dfrac{dx}{dt}=x,
\qquad
\dfrac{dy}{dt}=2y
$
\hfill (5.9)
\\[2.5ex]
are $O\!$-curves of the family   
\\[1.5ex]
\mbox{}\hfill
$
C_1^{}y+C_2^{}x^2=0,
\quad  
C_1^{},C_2^{}\in \R,
\hfill
$
\\[1.5ex]
and the equilibrium state (simple unstable node) $O(0,0).$

Trajectories of  the stereographically conjugate system
\\[2ex]
\mbox{}\hfill                        
$
\dfrac{du}{d\tau}={}-u^3-3uv^2,
\qquad
\dfrac{dv}{d\tau}={}-2v^3, 
\qquad 
(u^2+v^2)\, d\tau=dt,
$
\hfill (5.10)
\\[2.5ex]
are $O^{\ast}\!$-curves of the family
\\[1.5ex]
\mbox{}\hfill
$
C_1^{\ast}v(u^2+v^2)+C_2^{\ast}\, u^2=0,
\quad  
C_1^{\ast},C_2^{\ast}\in \R,
\quad
C_1^{\ast}=C_1^{},
\quad 
C_2^{\ast}=4C_2^{},
\hfill
$
\\[1.5ex]
and the complicated equilibrium state (stable node) $O^{\ast} (0,0).$
\vspace{0.25ex}

The behaviour of trajectories on the sphere (1.1) for the stereographically mutually conjugate systems (5.9) and (5.10)
is represented on Fig. 5.14.
The circles on Fig. 5.15 are stereographic atlases of trajectories for the differential systems (5.9) and (5.10).
\\[4ex]
\parbox{50mm}{
\mbox{}\hfill
{\unitlength=1mm
\begin{picture}(42,42)
\put(0,0){\includegraphics[width=42mm,height=42mm]{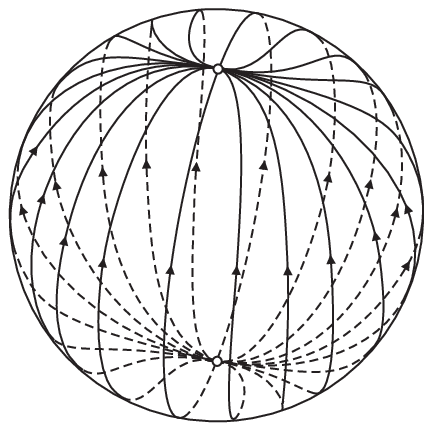}}
\end{picture}}
\hfill\mbox{}
\\[2ex]
\mbox{}\hfill
Fig. 5.14
\hfill\mbox{}
}
\parbox{104mm}{
\mbox{}\hfill
{\unitlength=1mm
\begin{picture}(42,42)
\put(0,0){\includegraphics[width=42mm,height=42mm]{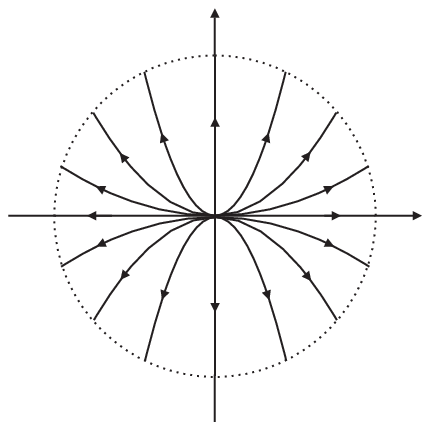}}
\put(41,19){\makebox(0,0)[cc]{\scriptsize $x$}}
\put(19,41){\makebox(0,0)[cc]{\scriptsize $y$}}
\end{picture}}
\qquad
{\unitlength=1mm
\begin{picture}(42,42)
\put(0,0){\includegraphics[width=42mm,height=42mm]{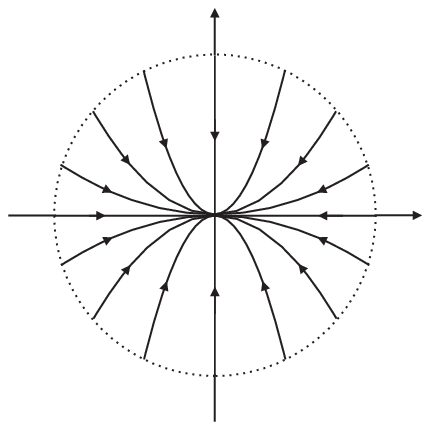}}
\put(41,19){\makebox(0,0)[cc]{\scriptsize $u$}}
\put(19,41){\makebox(0,0)[cc]{\scriptsize $v$}}
\end{picture}}
\hfill\mbox{}
\\[2ex]
\mbox{}\hfill
Fig. 5.15
\hfill\mbox{}
}
\\[3ex]
\indent
{\bf Example 5.7.}
Trajectories of the linear system
\\[2ex]
\mbox{}\hfill                        
$
\dfrac{dx}{dt}=x+y,
\qquad
\dfrac{dy}{dt}=y
$
\hfill (5.11)
\\[2ex]
are $O\!$-curves of the family  
\\[1.5ex]
\mbox{}\hfill
$
y \exp \Bigl({}-\dfrac{x}{y}\Bigr)=C,
\quad 
C\in \R,
\hfill
$
\\[1.5ex]
and the equilibrium state (unstable degenerate node) $O (0,0).$
\vspace{0.25ex}

Trajectories of  the stereographically conjugate system
\\[2ex]
\mbox{}\hfill                        
$
\dfrac{du}{d\tau}={}-u^3-u^2\;\!v-uv^2+v^3,
\quad
\dfrac{dv}{d\tau}={}-u^2\;\!v -2uv^2-v^3, \ \ (u^2+v^2)\, d\tau=dt,
$
\hfill (5.12)
\\[2ex]
are $O^{\ast}\!$-curves of the family  
\\[1.5ex]
\mbox{}\hfill
$
\dfrac{v}{u^2+v^2 }\, \exp\Bigl({}- \dfrac{u}{v}\Bigr)=C^{\ast},
\quad
C^{\ast}\in\R, 
\quad 
C^{\ast}=4C,
\hfill
$
\\[1.5ex]
and the complicated  equilibrium state (stable node) $O^{\ast} (0,0).$
\vspace{0.25ex}

Trajectories of the stereographically mutually conjugate systems (5.11) and (5.12)
on the sphere (1.1) are represented on Fig. 5.16.
The circles on Fig. 5.17 are stereographic atlases of trajectories for the differential systems (5.11) and (5.12).
\\[3.5ex]
\parbox{50mm}{
\mbox{}\hfill
{\unitlength=1mm
\begin{picture}(42,42)
\put(0,0){\includegraphics[width=42mm,height=42mm]{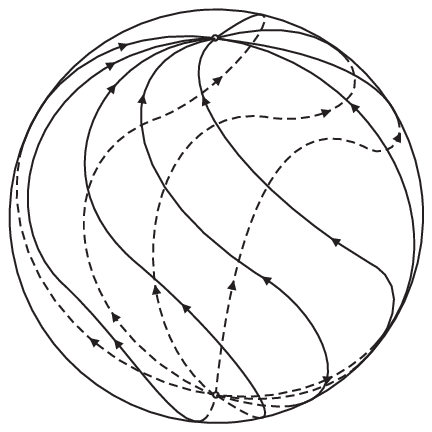}}
\end{picture}}
\hfill\mbox{}
\\[2ex]
\mbox{}\hfill
Fig. 5.16
\hfill\mbox{}
}
\parbox{104mm}{
\mbox{}\hfill
{\unitlength=1mm
\begin{picture}(42,42)
\put(0,0){\includegraphics[width=42mm,height=42mm]{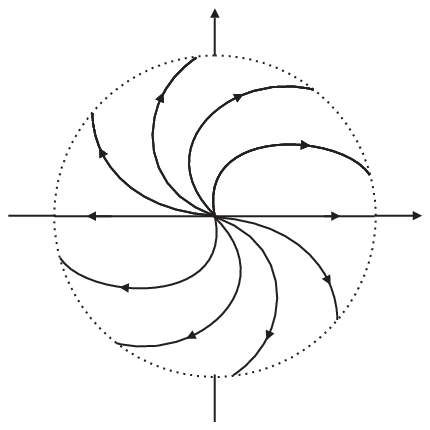}}
\put(41,19){\makebox(0,0)[cc]{\scriptsize $x$}}
\put(19,41){\makebox(0,0)[cc]{\scriptsize $y$}}
\end{picture}}
\qquad
{\unitlength=1mm
\begin{picture}(42,42)
\put(0,0){\includegraphics[width=42mm,height=42mm]{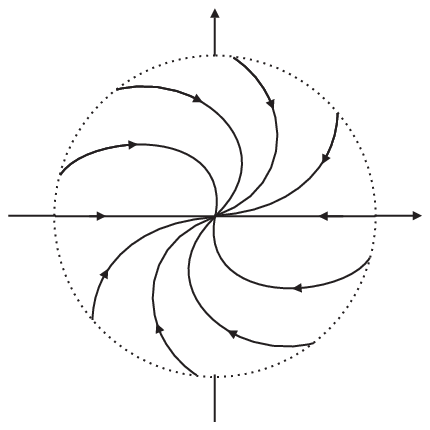}}
\put(41,19){\makebox(0,0)[cc]{\scriptsize $u$}}
\put(19,41){\makebox(0,0)[cc]{\scriptsize $v$}}
\end{picture}}
\hfill\mbox{}
\\[2ex]
\mbox{}\hfill
Fig. 5.17
\hfill\mbox{}
}
\\[-4ex]

\newpage

\mbox{}
\\[-2.5ex]
\centerline{
{\bf\large \S\;\!3. Trajectories of stereographically conjugate differential systems}}
\\[1.5ex]
\centerline{
\bf  6. Regular points and equilibrium states }
\\[0.25ex]
\centerline{
\bf of stereographically conjugate differential systems}
\\[1.5ex]
\indent
Suppose the point $\!M(x,y)\!$ 
\vspace{0.5ex}
of the phase plane $\!Oxy\!$ is distinct from the origin of coordi- nates $O(0,0).$
\vspace{0.5ex}
Then the image of the point $M(x,y)$ of the phase plane $Oxy$ for system (D)
under the diffeomorphic map (2.6) is
the point $M^* \Bigl(\dfrac{4x}{x^{2}+y^{2}}\,,\ \dfrac{4y}{x^{2}+y^{2}}\Bigr)$
\vspace{0.75ex}
of the phase plane $O^*uv$ for system (3.3) and 
\vspace{0.75ex}
the point $M^*$ is distinct from the origin of coordinates $O^*(0,0).$

Let the image of the origin of coordinates $O(0,0)$ of the extended plane $\overline{Oxy}$ be
\vspace{0.35ex}
the infinitely remote point $M^{\ast}_{\infty}$ of the stereographically conjugate extended plane $\overline{O^{\ast}uv},$
\vspace{0.35ex}
and let the image of  the infinitely remote point  $M^{}_{\infty}$ of the extended plane $\overline{Oxy}$ be
\vspace{0.5ex}
the the origin of coordinates $O^{\ast}(0,0)$ of the stereographically conjugate extended plane $\overline{O^{\ast}uv}.$
\vspace{0.5ex}

Then such extension of the map $\varphi_{21}^{}$ is the bijective map $\overline{\varphi_{21}^{}}$ 
\vspace{0.75ex}
of the extended phase plane $\overline{Oxy}$ of system (D) to the extended phase plane $\overline{O^{\ast}uv}$ of system (3.3).
\vspace{0.5ex}
The points $M$ and $M^{\ast},\ O$ and $M^{\ast}_{\infty}\,,\ M^{}_{\infty}$ and $O^{\ast}$ are 
\vspace{0.5ex}
stereographically mutually conjugate.
Also a curve  $l$ on the extended plane $\overline{Oxy}$ and its image $l^\ast$ on the extended plane $\overline{O^{\ast}uv}$ 
\vspace{0.35ex}
under the map $\overline{\varphi_{21}^{}}$ are stereographically mutually conjugate curves.
\vspace{0.5ex}

The map (2.6) is the superposition (2.5) of stereographic maps. Then,
\vspace{0.15ex}
taking into account the basic property of stereographic projection (Subsection 1), we get
\vspace{0.35ex}

{\bf Property 6.1.}
{\it
Angle between curves is equal to angle between stereographically conjugate curves to them.}
\vspace{0.35ex}

A point of the phase plane $Oxy$ is 
\vspace{0.35ex}
a regular point of system (D) if this point is not an equilibrium state of system (D).
\vspace{0.5ex}
The infinitely remote point $M_{\infty}^{}$ of the extended phase plane $\overline{Oxy}$ is called  
\vspace{0.35ex}
{\it regular infinitely remote point} of system (D) if  
the origin of coordinates $O^{\ast}(0,0)$ of the phase plane $O^{\ast}uv$ is not an  equilibrium state of system (3.3). 
\vspace{0.35ex}
If the point $O^{\ast}(0,0)$ is an equilibrium state of system (3.3), 
\vspace{0.5ex}
then the infinitely remote point $M^{}_{\infty}$ of the extended phase plane $\overline{Oxy}$ is 
infinitely remote equilibrium state of system (D) of the same form.
\vspace{0.5ex}

Thus every point 
\vspace{0.35ex}
of the extended phase plane $\overline{Oxy}$ is either a regular point or an equilibrium state of system (D). 
\vspace{0.15ex}
Every trajectory of system (D) on  the extended phase plane $\overline{Oxy}$ is an equilibrium state or consists of regular points. 
\vspace{0.35ex}

The map (2.6) of the plane $Oxy$ without the origin of coordinates $O(0,0)$ to 
\vspace{0.35ex}
the stereographically conjugate plane $O^{\ast}uv$ 
\vspace{0.35ex}
without the origin of coordinates $O^{\ast}(0,0)$ is diffeomorphic.
Then, by Property 6.1, 
\vspace{0.5ex}
for trajectories of the stereographically mutually conjugate systems (D) and (3.3)
on the extended phase planes 
\vspace{0.5ex}
$\overline{Oxy}$ and $\overline{O^{\ast}uv},$ we have 

{\bf Property 6.2.}
\vspace{0.5ex}
{\it
Stereographically mutually conjugate points of  the extended phase planes 
$\overline{Oxy}$ and $\overline{O^{\ast}uv}$ are simultaneously either regular poits  or 
\vspace{0.25ex}
equilibrium states for the same form of the differential systems {\rm (D)} and {\rm (3.3)}.}
\vspace{0.5ex}

Since the map (2.6) is the superposition (2.5), we have  
\vspace{0.5ex}

{\bf Property\! 6.3\! (6.4).}\!
\vspace{0.25ex}
{\it
The image of closed curve on the plane $\!Oxy,\!\!$ which is passing {\rm(}not pas\-sing{\rm)} 
through the origin of coordinates $\!O(0,0),\!\!$ is  not closed {\rm(}closed{\rm)} curve 
\vspace{0.25ex}
on the stereog\-ra\-p\-hi\-cal\-ly conjugate plane $\!O^{\ast}uv,\!\!$ which is not passing through  
\vspace{0.35ex}
the origin of coordinates $\!O^{\ast}(0,0).$}

For example, using calculations and the Bendixon transformation (3.1), we obtain
\vspace{0.35ex}

{\bf Property 6.5.}
{\it
The image of\;\!{\rm:}

a{\rm)} a circle 
\\[0.15ex]
\mbox{}\hfill
$
(x+a)^2+(y+b)^2=a^2+b^2,
\quad 
|a|+|b|\ne 0,
\hfill
$
\\[1ex]
\mbox{}\qquad\ \
which is passing through the origin of coordinates $O(0,0)$ of the plane $Oxy;$


b{\rm)} a circle 
\\[0.5ex]
\mbox{}\hfill
$
(x+a)^2+(y+b)^2=r^2,
\quad 
r>0,
\quad 
r^2\ne a^2+b^2,
\hfill
$
\\[1.5ex]
\mbox{}\qquad\ \
which is not passing through the origin of coordinates $O(0,0)$ of the plane $Oxy;$
\vspace{0.5ex}

c{\rm)} a circle 
\\[0.5ex]
\mbox{}\hfill
$
x^2+y^2=r^2,
\quad 
r>0,
\hfill
$
\\[1.5ex]
\mbox{}\qquad\ \
with the center in the origin of coordinates $O(0,0)$ of the plane $Oxy;$
\vspace{0.75ex}

d{\rm)} 
a point $A(a,b),\ |a|+|b|\ne 0,$ which lies inside {\rm(}outside{\rm)} of a circle 
\\[1.5ex]
\mbox{}\hfill
$
x^2+y^2=r^2,
\quad 
r>0;
\hfill
$
\\[1ex]
\indent
e{\rm)}
a straight line 
\\[0.5ex]
\mbox{}\hfill
$
Ax+By=0,
\quad 
|A|+|B|\ne 0,
\hfill
$
\\[1.5ex]
\mbox{}\qquad\ \
which  is passing through  the origin of coordinates $O(0,0)$ of the plane $Oxy;$
\vspace{0.75ex}

f{\rm)} 
the coordinate straight line $x=0\ (y=0)$ of the plane $Oxy;$
\vspace{0.5ex}

g{\rm)} 
a straight line  
\\[1ex]
\mbox{}\hfill
$
Ax+By+C=0,
\quad 
|A|+|B|\ne 0,
\quad 
C\ne 0,
\hfill
$
\\[1ex]
\mbox{}\qquad\ \
which is not passing through  the origin of coordinates $O(0,0)$ of the plane $Oxy;$
\\[0.75ex]
on the stereographically conjugate plane $O^\ast uv$ is\;\!{\rm:}
\vspace{0.5ex}

a{\rm)}
the straight line 
\\[0.5ex]
\mbox{}\hfill
$
au+bv+2=0;
\hfill
$
\\[0.75ex]
\indent
b{\rm)}
the circle  
\\[1ex]
\mbox{}\hfill
$
\Bigl( u+\dfrac{4a}{a^2+b^2-r^2}\Bigr)^{\!2}+\Bigl( v+\dfrac{4b}{a^2+b^2-r^2}\Bigr)^{\!2}=\dfrac{4r^2}{(a^2+b^2-r^2)^2}\,;
\hfill
$
\\[1.5ex]
\indent
c{\rm)}
the circle  
\\[1ex]
\mbox{}\hfill
$
u^2+v^2=\dfrac{4}{r^2}\,;
\hfill
$
\\[1.5ex]
\indent
d{\rm)}
the stereographically conjugate point 
\\[1.5ex]
\mbox{}\hfill
$
A^\ast \biggl( \dfrac{4a}{a^2+b^2}\,,\, \dfrac{4b}{a^2+b^2}\biggr),
\hfill
$
\\[1.5ex]
\mbox{}\qquad\ \
which lies outside {\rm(}inside{\rm)} of the circle 
\\[1.25ex]
\mbox{}\hfill
$
u^2+v^2=\dfrac{4}{r^2}\,;
\hfill
$
\\[1.25ex]
\indent
e{\rm)}
the straight line 
\\[0.5ex]
\mbox{}\hfill
$
Au+Bv=0,
\hfill
$
\\[1ex]
\mbox{}\qquad\ \
which is passing through  the origin of coordinates $O^\ast(0,0);$ 
\vspace{0.5ex}

f{\rm)}
the coordinate straight line $u=0\ (v=0);$
\vspace{0.5ex}

g{\rm)}
the circle  
\\[1ex]
\mbox{}\hfill
$
\Bigl( u+\dfrac{2A}{C}\Bigr)^{\!2}+\Bigl(v+\dfrac{2B}{C}\Bigr)^{\!2}=\dfrac{4(A^2+B^2)}{C^2}\,,
\hfill
$
\\[2ex]
\mbox{}\qquad\ \
which is passing through the origin of coordinates $O^\ast(0,0).$}
\vspace{1.25ex}

The stereographically mutually conjugate differential systems (5.1) and (5.2) 
\vspace{0.15ex}
were considered
in Example 5.2. 
\vspace{0.25ex}
Trajectories of these systems  (Fig. 5.6 and 5.7) lie on the straight lines, 
which are passing through  the origins of coordinates $O(0,0)$ and $O^\ast(0,0)$ 
\vspace{0.25ex}
of the phase planes $Oxy$ and $O^\ast uv$ accordingly (Property 6.5, case e). 
\vspace{0.75ex}

\newpage

{\bf Example 6.1.}
\vspace{0.25ex}
If we replace $x$ by $x-1$ and $y$ by $y-1$ in the linear differential system (5.1), 
then we obtain the the linear differential system 
\\[1.5ex]
\mbox{}\hfill                        
$
\dfrac{dx}{dt}=x-1,
\qquad
\dfrac{dy}{dt}=y-1.
$
\hfill (6.1)
\\[1.75ex]
\indent
Trajectories of the linear differential system (6.1) are the unstable dicritical node $A(1,1)$
and $A\!$-rays of  the family of lines
\\[1.5ex]
\mbox{}\hfill
$
\dfrac{y-1}{x-1}=C,
\quad 
{}-\infty\leq C\leq {}+\infty.
\hfill
$ 
\\[1.5ex]
\indent
Trajectories of the stereographically conjugate differential system 
\\[2ex]
\mbox{}\hfill                        
$
\dfrac{du}{dt}={}-u+\dfrac{1}{4}\; u^2+\dfrac{1}{2}\; uv-\dfrac{1}{4}\; v^2,
\qquad
\dfrac{dv}{dt}={}-v-\dfrac{1}{4}\; u^2+\dfrac{1}{2}\; uv+\dfrac{1}{4}\; v^2
$
\hfill (6.2)
\\[2.5ex]
are the dicritical nodes $O^\ast (0,0)$ (stable) and $A^\ast (2,2)$ (unstable),
\vspace{0.35ex}
the segment $O^\ast\, A^\ast$ without endpoints, 
$O^\ast\!$-ray, $A^\ast\!$-ray of the straight line $v=u,$ and 
the arcs of the circles 
\\[1.75ex]
\mbox{}\hfill
$
\Bigl( u-\dfrac{2C}{C-1}\Bigr)^{\!2}+\Bigl( v+\dfrac{2}{C-1}\Bigr)^{\!2}=
\dfrac{4(C^2+1)}{(C-1)^2}\,,
\quad 
C\in \R\backslash\{1\},
\hfill
$
\\[1.75ex]
which are adjoining to the equilibrium states  $O^\ast$ and $A^\ast$
\vspace{0.25ex}
(the centres of circles lie on the straight line $u+v-2=0).$
\vspace{0.25ex}

The circles on Fig. 6.1 are stereographic atlases of trajectories for the differential systems (6.1) and (6.2).
\\[3ex]
\mbox{}\hfill
{\unitlength=1mm
\begin{picture}(45,45)
\put(0,0){\includegraphics[width=45mm,height=45mm]{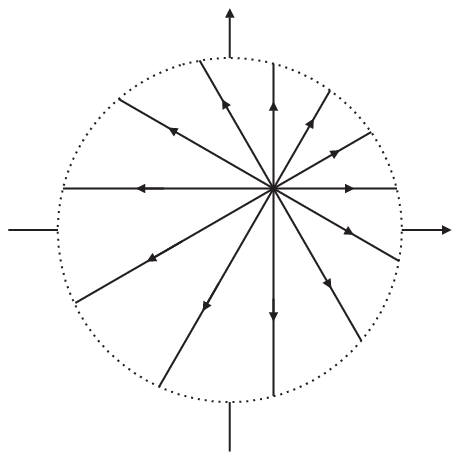}}
\put(44,20.5){\makebox(0,0)[cc]{\scriptsize $x$}}
\put(20.5,44){\makebox(0,0)[cc]{\scriptsize $y$}}
\end{picture}}
\qquad\qquad
{\unitlength=1mm
\begin{picture}(45,45)
\put(0,0){\includegraphics[width=45mm,height=45mm]{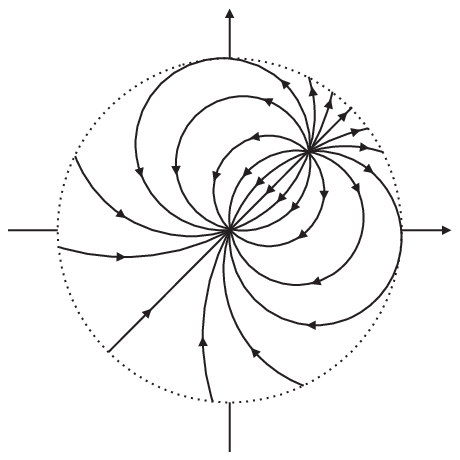}}
\put(44,20.5){\makebox(0,0)[cc]{\scriptsize $u$}}
\put(20.5,44){\makebox(0,0)[cc]{\scriptsize $v$}}
\end{picture}}
\hfill\mbox{}
\\[1.5ex]
\mbox{}\hfill
Fig. 6.1
\hfill\mbox{}
\\[3.75ex]
\indent
The stereographically mutually conjugate differential systems (5.3) and (5.4) 
\vspace{0.25ex}
were considered
in Example 5.3. 
\vspace{0.25ex}
Trajectories of these systems  (Fig. 5.8 and 5.9) are  the concentric circles, which are collapsing to the centres 
$O$ and $O^\ast$ accordingly (Property 6.5, case c).
\vspace{1ex}

{\bf Example 6.2.}
\vspace{0.25ex}
If we replace $x$ by $x-1$ and $y$ by $y-1$ in the linear differential system (5.3), 
then we get the the linear differential system 
\\[1.5ex]
\mbox{}\hfill                        
$
\dfrac{dx}{dt}=y-1,
\qquad
\dfrac{dy}{dt}={}-x+1.
$
\hfill (6.3)
\\[1.5ex]
\indent
Trajectories of  system (6.3) are the center  $A(1,1)$ and the concentric circles  
\\[1.5ex]
\mbox{}\hfill
$
(x-1)^2\,+\,(y-1)^2=C,
\quad 
C\in (0;{}+\infty).
\hfill
$ 
\\[1.5ex]
\indent
Trajectories of the stereographically conjugate system 
\\[2ex]
\mbox{}\hfill                        
$
\dfrac{du}{dt}=v+\dfrac{1}{4}\; u^2-\dfrac{1}{2}\; uv-\dfrac{1}{4}\; v^2,
\qquad
\dfrac{dv}{dt}={}-u+\dfrac{1}{4}\; u^2+\dfrac{1}{2}\; uv-\dfrac{1}{4}\; v^2
$
\hfill (6.4)
\\[2.25ex]
are the straight line $u+v-2=0,$ the circles  
\\[1.75ex]
\mbox{}\hfill
$
\Bigl( u-\dfrac{4}{2-C}\Bigr)^{\!2}+\Bigl( v-\dfrac{4}{2-C}\Bigr)^{\!2}=\dfrac{16\;\!C}{(2-C)^2}\,,
\quad 
C\in (0;2)\sqcup (2;{}+\infty),
\hfill
$
\\[1.5ex]
the centres $O^\ast (0,0)$ and $A^\ast (2,2).$
\vspace{0.25ex}

The circles on Fig. 6.2 are stereographic atlases of trajectories for the differential systems (6.3) and (6.4).
\\[3ex]
\mbox{}\hfill
{\unitlength=1mm
\begin{picture}(45,45)
\put(0,0){\includegraphics[width=45mm,height=45mm]{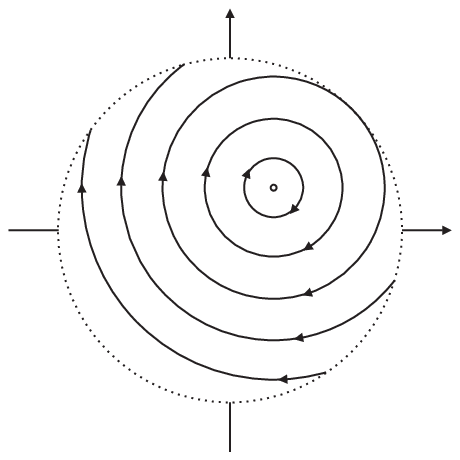}}
\put(44,20.5){\makebox(0,0)[cc]{\scriptsize $x$}}
\put(20.5,44){\makebox(0,0)[cc]{\scriptsize $y$}}
\end{picture}}
\qquad\qquad
{\unitlength=1mm
\begin{picture}(45,45)
\put(0,0){\includegraphics[width=45mm,height=45mm]{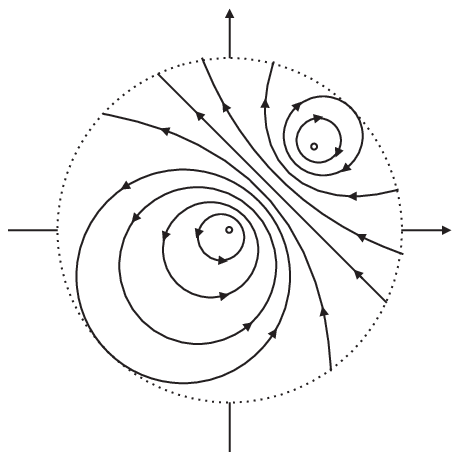}}
\put(44,20.5){\makebox(0,0)[cc]{\scriptsize $u$}}
\put(20.5,44){\makebox(0,0)[cc]{\scriptsize $v$}}
\end{picture}}
\hfill\mbox{}
\\[0ex]
\mbox{}\hfill
Fig. 6.2
\hfill\mbox{}
\\[2.5ex]
\indent
{\bf  Example 6.3.}
Trajectories of Jacobi's system
\\[2ex]
\mbox{}\hfill        
$
\dfrac{dx}{dt}=1+x-y+x(x+y-1),
\qquad 
\dfrac{dy}{dt}=y(x+y-1)
$
\hfill(6.5)
\\[2ex]
are the curves of family 
\\[1.5ex]
\mbox{}\hfill
$
\dfrac{x^2+(y-1)^2}{y^2}\,\exp \Bigl({}-2\arctan \dfrac{y-1}{x}\Bigr)=C,
\quad
0\leq C\leq {}+\infty.
\hfill
$
\\[1.5ex]
Among these curves are the straight line $y=0$ and the unstable focus $A(0,1).$
\\[3ex]
\mbox{}\hfill
{\unitlength=1mm
\begin{picture}(45,45)
\put(0,0){\includegraphics[width=45mm,height=45mm]{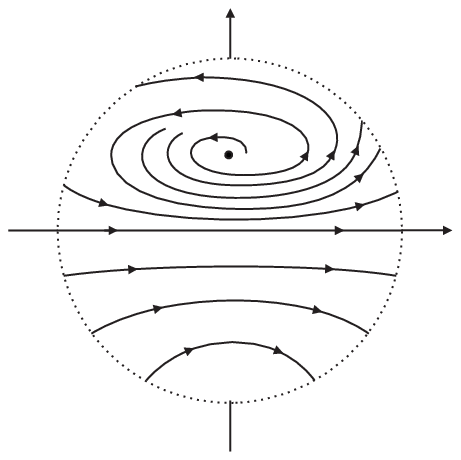}}
\put(44,20.5){\makebox(0,0)[cc]{\scriptsize $x$}}
\put(20.5,44){\makebox(0,0)[cc]{\scriptsize $y$}}
\end{picture}}
\qquad\qquad
{\unitlength=1mm
\begin{picture}(45,45)
\put(0,0){\includegraphics[width=45mm,height=45mm]{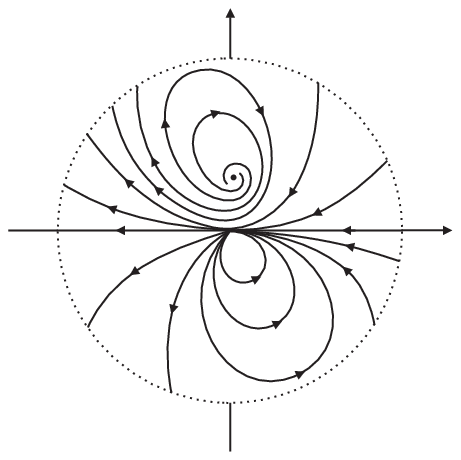}}
\put(44,20.5){\makebox(0,0)[cc]{\scriptsize $u$}}
\put(20.5,44){\makebox(0,0)[cc]{\scriptsize $v$}}
\end{picture}}
\hfill\mbox{}
\\[0ex]
\mbox{}\hfill
Fig. 6.3
\hfill\mbox{}
\\[2.5ex]
\indent
Trajectories of the stereographically conjugate system 
\\[1.75ex]
\mbox{}\hfill        
$
\dfrac{du}{d\theta}={}-4u^2-4uv+u^2v+2uv^2-v^3-\dfrac{1}{4}\, u^4+\dfrac{1}{4}\, v^4,
\hfill
$
\\
\mbox{}\hfill (6.6)
\\
\mbox{}\hfill
$
\dfrac{dv}{d\theta}=v\Bigl({}-4u-4v-u^2+2uv+v^2-\dfrac{1}{2}\, u^3-\dfrac{1}{2}\,uv^2\Bigr),
\hfill
$
\\[2.5ex]
where $(u^2+v^2)\;\!d\theta=dt,$
are the curves of family  
\\[2ex]
\mbox{}\hfill
$
\dfrac{16u^2+(u^2+v^2-4v)^2}{v^2}\,\exp \Bigl({}-2\arctan \dfrac{4v-u^2-v^2}{4u}\Bigr)=C^\ast,
\ \ 
0\leq C^\ast\leq {}+\infty,\ C^\ast=16C.
\hfill
$
\\[2.25ex]
Among these curves are $O^\ast\!$-rays of the straight line $v=0,$
the unstable focus $A^\ast(0,2),$
and the complicated equilibrium state $O^\ast(0,0),$
which is consisting from hyperbolic, elliptic and accompanying it two parabolic Bendixon's sectors.
 
The circles on Fig. 6.3 are stereographic atlases of trajectories for the differential systems (6.5) and (6.6).
\vspace{0.75ex}

{\bf  Example 6.4.}
\vspace{0.25ex}
If we replace $x$ by $x$ and $y$ by $y+1$ in the differential system (6.5), 
then we have Jacobi's differential system 
\\[2ex]
\mbox{}\hfill        
$
\dfrac{dx}{dt}=x-y+x(x+y),
\qquad 
\dfrac{dy}{dt}=(y+1)(x+y).
$
\hfill(6.7)
\\[2.5ex]
\indent
The curves of family   
\\[1.5ex]
\mbox{}\hfill
$
\dfrac{x^2+y^2}{(y-1)^2}\,\exp \Bigl({}-2\arctan \dfrac{y}{x}\Bigr)=C,
\quad
0\leq C\leq {}+\infty,
\hfill
$
\\[1.5ex]
are trajectories of Jacobi's differential system (6.7). 
\vspace{0.25ex}
Among these curves are the  straight line $y={}-1$ and the unstable focus $O(0,0).$
\vspace{0.35ex}

Trajectories of the stereographically conjugate system 
\\[2ex]
\mbox{}\hfill        
$
\dfrac{du}{d\theta}={}-(4u^2+4uv+u^3+u^2v+uv^2+v^3),
\quad 
\dfrac{dv}{d\theta}={}-4uv-4v^2+u^3-u^2v+uv^2-v^3,
$
\hfill(6.8)
\\[2.5ex]
where $(u^2+v^2)\;\!d\theta=dt,$ are the curves of family 
\\[2ex]
\mbox{}\hfill
$
\dfrac{u^2+v^2}{(u^2+v^2-4v)^2}\,\exp \Bigl({}-2\arctan \dfrac{v}{u}\Bigr)=C^\ast,
\quad
0\leq C^\ast\leq {}+\infty,
\quad 
16C^\ast=C.
\hfill
$
\\[2ex]
Among these curves are the complicated equilibrium state $O^\ast(0,0),$ 
\vspace{0.25ex}
which is consisting from hyperbolic, elliptic and accompanying it two parabolic Bendixon's sectors .
\vspace{0.35ex}

The circles on Fig. 6.4 are stereographic atlases of trajectories for the differential systems (6.7) and (6.8).
\\[3ex]
\mbox{}\hfill
{\unitlength=1mm
\begin{picture}(45,45)
\put(0,0){\includegraphics[width=45mm,height=45mm]{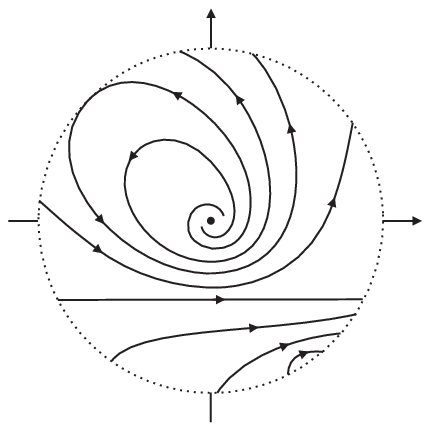}}
\put(44,20.5){\makebox(0,0)[cc]{\scriptsize $x$}}
\put(20.5,44){\makebox(0,0)[cc]{\scriptsize $y$}}
\end{picture}}
\qquad\qquad
{\unitlength=1mm
\begin{picture}(45,45)
\put(0,0){\includegraphics[width=45mm,height=45mm]{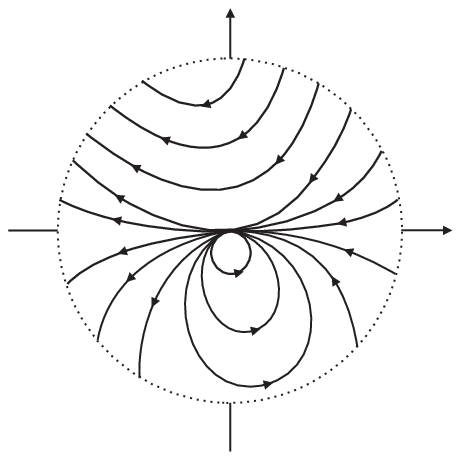}}
\put(44,20.5){\makebox(0,0)[cc]{\scriptsize $u$}}
\put(20.5,44){\makebox(0,0)[cc]{\scriptsize $v$}}
\end{picture}}
\hfill\mbox{}
\\[2ex]
\mbox{}\hfill
Fig. 6.4
\hfill\mbox{}
\\[4.75ex]
\centerline{\bf  7. Stereographic cycles}
\\[1.25ex]
\indent
Using Property 6.4 and the diffeomorphic map (2.6), we obtain
\vspace{0.5ex}

{\bf Property 7.1 (7.2).}
\vspace{0.25ex}
{\it
The image of a cycle {\rm(}limit cycle{\rm)} of the differential system {\rm(D)}, which is
\vspace{0.25ex}
not passing through the origin of coordinates $O(0,0)$ of the phase plane $Oxy,$  
on the phase plane $O^{\ast}uv$ is a cycle {\rm(}limit cycle{\rm)} 
\vspace{0.25ex}
of the differential system {\rm(3.3)}, which is 
not passing through the origin of coordinates $O^{\ast}(0,0).$}
\vspace{0.5ex}

Using Property 6.3 and the diffeomorphic map (2.6), we get
\vspace{0.5ex}

{\bf Property 7.3 (7.4).}
\vspace{0.25ex}
{\it
The image of a cycle {\rm(}limit cycle{\rm)} of the differential system {\rm(D)}, which is
\vspace{0.25ex}
passing through the origin of coordinates $O(0,0)$ of the phase plane $Oxy,$ 
on  the phase plane $O^{\ast}uv$ is an open trajectory of the differential system {\rm(3.3)}, 
\vspace{0.25ex}
which is not passing through the origin of coordinates $O^{\ast}(0,0).$
}

\newpage

{\bf Definition 7.1.}
A {\it stereographic cycle} of system (D) is a trajectory of system (D) 
such that the stereographic image of this trajectory on the sphere (1.1) is a closed curve and 
each point of this curve is image of a regular point (final or infinitely remove) of system (D). 
\vspace{0.75ex}

{\bf Definition 7.2.}
A stereographic cycle of system (D) is called a
{\it limit stereographic cycle} of system (D) if
the stereographic image of this stereographic cycle on the sphere (1.1) has an neighbourhood 
without an stereographic image of other stereographic cycle of system (D).
\vspace{0.75ex}

{\bf Property 7.5 (7.6).}
\vspace{0.15ex}
{\it
A cycle {\rm(}limit cycle{\rm)} of system {\rm(D)} is a stereographic cycle {\rm(}limit stereographic cycle{\rm)}
of system {\rm(D)}.
}
\vspace{0.75ex}

{\bf Definition 7.3.}
\vspace{0.35ex}
A stereographic cycle {\rm(}limit stereographic cycle{\rm)} 
is called {\it  open} if this cycle is passing through the infinitely remote point $M^{}_{\infty}$
of the extended phase plane $\overline{Oxy}.$ 
\vspace{0.75ex}

{\bf Property 7.7 (7.8).}
\vspace{0.25ex}
{\it
The image of a stereographic cycle {\rm(}limit stereographic cycle{\rm)} of the differential system {\rm(D)},
\vspace{0.25ex}
which is passing through the origin of coordinates $O(0,0)$ of the phase plane $Oxy,$ 
\vspace{0.25ex}
on the phase plane $O^{\ast}uv$
is an open stereographic cycle {\rm(}open limit stereographic cycle{\rm)}
of the differential system} (3.3).
\vspace{0.75ex}

{\bf Property 7.9 (7.10).}
\vspace{0.25ex}
{\it
The image of a cycle {\rm(}limit cycle{\rm)} of system {\rm(D)}, 
which is passing through the origin of coordinates $O(0,0)$ of the phase plane $Oxy,$
\vspace{0.25ex}
on the phase plane $O^{\ast}uv$ 
is an open stereographic cycle {\rm(}open limit stereographic cycle{\rm)}
of system} (3.3).
\vspace{0.5ex}

For instance, 
\vspace{0.15ex}
the straight line-trajectory $u+v-2=0$ is an open stereographic cycle of 
the differential system (6.4).
\vspace{1ex}

{\bf Example 7.1.}
Trajectories of the differential system [8, p. 88]
\\[2ex]
\mbox{}\hfill        
$
\dfrac{dx}{dt}=x(x^2+y^2-1)-y(x^2+y^2+1),
\qquad
\dfrac{dy}{dt}=x(x^2+y^2+1)+y(x^2+y^2-1)
$
\hfill (7.1)
\\[2ex]
are the curves [9; 1]
\\[2ex]
\mbox{}\hfill
$
\dfrac{x^2+y^2}{(x^2+y^2-1)^2}\,\exp \Bigl(2\arctan \dfrac{y}{x}\Bigr)=C,
\quad
0\leq C\leq {}+\infty.
\hfill
$
\\[2ex]
\indent
Among these curves are the equilibrium state $O(0,0)$ (stable focus) 
\vspace{0.35ex}
and the unstable limit cycle $x^2+y^2=1.$
\\[5ex]
\parbox{50mm}{
\mbox{}\hfill
{\unitlength=1mm
\begin{picture}(45,45)
\put(0,0){\includegraphics[width=45mm,height=45mm]{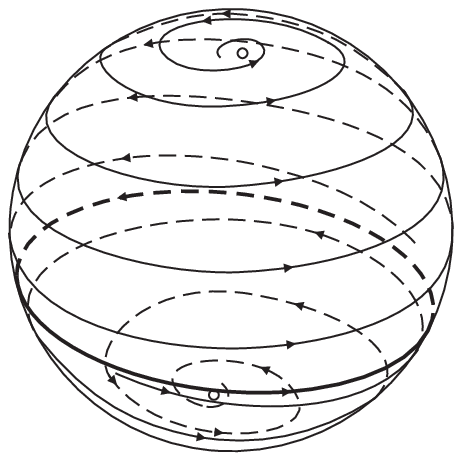}}
\end{picture}}
\hfill\mbox{}
\\[2ex]
\mbox{}\hfill
Fig. 7.1
\hfill\mbox{}
}
\parbox{104mm}{
\mbox{}\hfill
{\unitlength=1mm
\begin{picture}(45,45)
\put(0,0){\includegraphics[width=45mm,height=45mm]{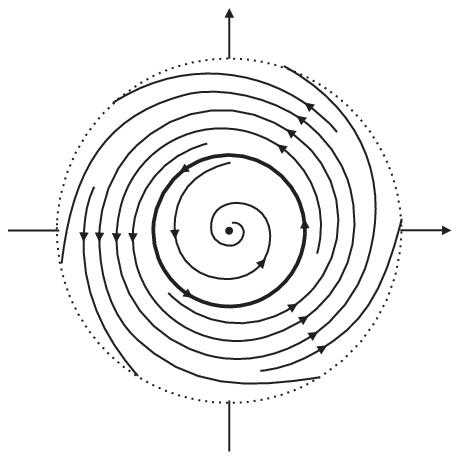}}
\put(44,20.5){\makebox(0,0)[cc]{\scriptsize $x$}}
\put(20.5,44){\makebox(0,0)[cc]{\scriptsize $y$}}
\end{picture}}
\qquad
{\unitlength=1mm
\begin{picture}(45,45)
\put(0,0){\includegraphics[width=45mm,height=45mm]{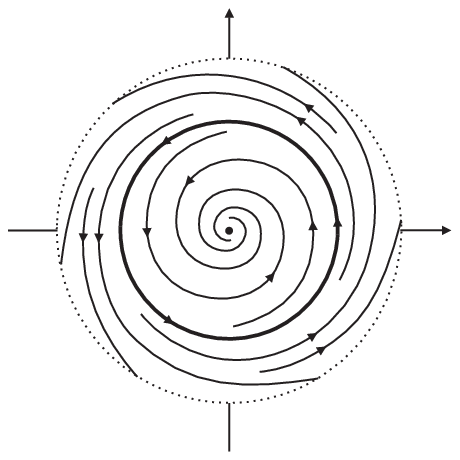}}
\put(44,20.5){\makebox(0,0)[cc]{\scriptsize $u$}}
\put(20.5,44){\makebox(0,0)[cc]{\scriptsize $v$}}
\end{picture}}
\hfill\mbox{}
\\[2ex]
\mbox{}\hfill
Fig. 7.2
\hfill\mbox{}
}
\\[5ex]
\indent
Trajectories of the stereographically conjugate system 
\\[2ex]
\mbox{}\hfill        
$
\dfrac{du}{d\theta}=u(u^2+v^2-16)-v(u^2+v^2+16),
\quad\
\dfrac{dv}{d\theta}=u(u^2+v^2+16)+v(u^2+v^2-16),
$
\hfill (7.2)
\\[3ex]
where $(u^2+v^2)\,d\theta=dt,$ are the curves
\\[2ex]
\mbox{}\hfill
$
\dfrac{u^2+v^2}{(u^2+v^2-16)^2}\,\exp \Bigl(2\arctan \dfrac{v}{u}\Bigr)=C^\ast,
\quad
0\leq C^\ast\leq {}+\infty,
\quad 
16C^\ast=C.
\hfill
$
\\[2ex]
\indent
Among these curves are the equilibrium state $O^\ast(0,0)$ (stable focus) 
\vspace{0.25ex}
and the unstable limit cycle $u^2+v^2=16.$
\vspace{0.25ex}

Trajectories on the sphere (1.1) 
\vspace{0.15ex}
for the stereographically mutually conjugate differential systems (7.1) and (7.2) 
are represented on Fig. 7.1. 
\vspace{0.15ex}
The circles on Fig. 7.2 are stereographic atlases of trajectories for the differential systems (7.1) and (7.2).
\vspace{1ex}

{\bf  Example 7.2.}
Trajectories of Darboux's sustem [10]
\\[2ex]
\mbox{}\hfill                        
$
\dfrac{dx}{dt}={}-y-x(x^2+y^2-1),
\qquad
\dfrac{dy}{dt}=x-y(x^2+y^2-1)
$
\hfill (7.3)
\\[2ex]
are the curves [9; 1]
\\[2ex]
\mbox{}\hfill
$
\dfrac{x^2+y^2}{1-x^2-y^2}\, \exp\Bigl({}-2\arctan\dfrac{y}{x}\Bigr)=C,
\quad 
{}-\infty\leq C\leq {}+\infty.
\hfill
$
\\[2ex]
\indent
Among these curves are the equilibrium state $O(0,0)$
\vspace{0.25ex}
(unstable focus) and the stable limit cycle $x^2+y^2=1.$
\vspace{0.35ex}

Trajectories of the stereographically conjugate system 
\\[2ex]
\mbox{}\hfill                        
$
\dfrac{du}{d\theta}=
16u-u^3-u^2v-uv^2-v^3,   
\qquad
\dfrac{dv}{d\theta}=16v+u^3-u^2v+uv^2-v^3,
$
\hfill (7.4)
\\[2ex]
where $(u^2+v^2)\;\!d\theta=dt,$ are the curves
\\[2ex]
\mbox{}\hfill
$
\dfrac{1}{u^2+v^2-16}\, \exp\Bigl({}-2\arctan\dfrac{v}{u}\Bigr)=C^\ast,
\quad
{}-\infty  \leq C^\ast \leq {}+\infty,
\quad
16C^\ast=C.
\hfill
$
\\[2ex]
\indent
Among these curves are the equilibrium state $O^{\ast}(0,0)$ 
\vspace{0.25ex}
(unstable dicritical node) and the stable limit circle $u^2+v^2=16.$
\vspace{0.35ex}

Trajectories on the sphere (1.1) 
\vspace{0.15ex}
for the stereographically mutually conjugate differential systems (7.3) and (7.4) 
are represented on Fig. 7.3. 
\vspace{0.15ex}
The circles on Fig. 7.4 are stereographic atlases of trajectories for the differential systems (7.3) and (7.4).
\\[4.5ex]
\parbox{50mm}{
\mbox{}\hfill
{\unitlength=1mm
\begin{picture}(45,45)
\put(0,0){\includegraphics[width=45mm,height=45mm]{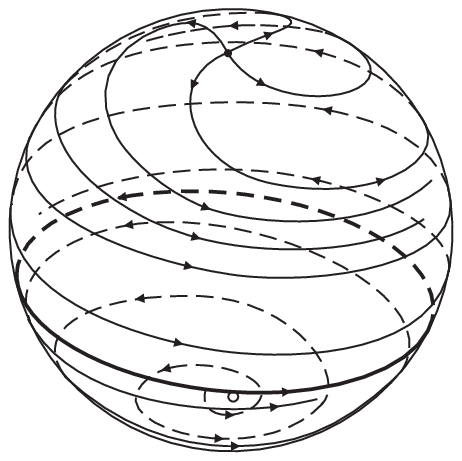}}
\end{picture}}
\hfill\mbox{}
\\[2ex]
\mbox{}\hfill
Fig. 7.3
\hfill\mbox{}
}
\parbox{104mm}{
\mbox{}\hfill
{\unitlength=1mm
\begin{picture}(45,45)
\put(0,0){\includegraphics[width=45mm,height=45mm]{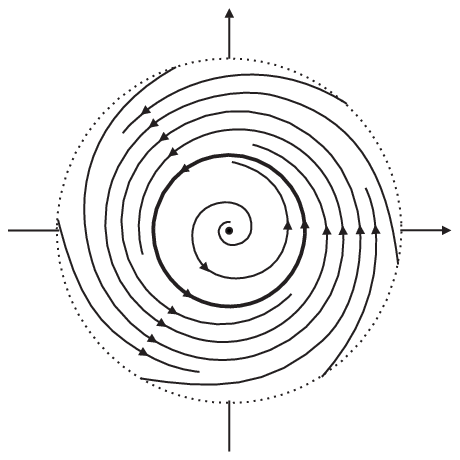}}
\put(44,20.5){\makebox(0,0)[cc]{\scriptsize $x$}}
\put(20.5,44){\makebox(0,0)[cc]{\scriptsize $y$}}
\end{picture}}
\qquad
{\unitlength=1mm
\begin{picture}(45,45)
\put(0,0){\includegraphics[width=45mm,height=45mm]{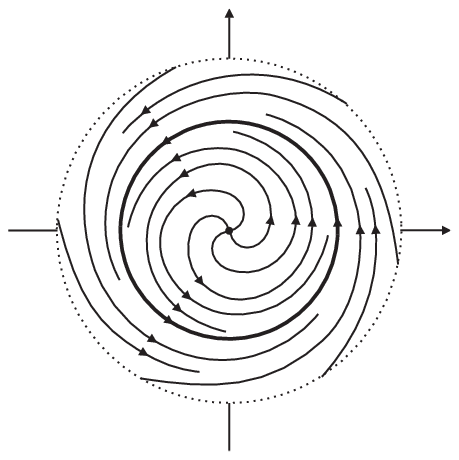}}
\put(44,20.5){\makebox(0,0)[cc]{\scriptsize $u$}}
\put(20.5,44){\makebox(0,0)[cc]{\scriptsize $v$}}
\end{picture}}
\hfill\mbox{}
\\[2ex]
\mbox{}\hfill
Fig. 7.4
\hfill\mbox{}
}
\\[6.5ex]
\indent
{\bf  Example 7.3.}
\vspace{0.25ex}
If we replace $x$ by $x-1$ and $y$ by $y$ in the differential system (7.3), 
then we obtain Darboux's differential system 
\\[2.25ex]
\mbox{}\hfill                        
$
\dfrac{dx}{dt}={}-2x-y+3x^2+y^2-x(x^2+y^2),
\qquad
\dfrac{dy}{dt}={}-1+x+2xy-y(x^2+y^2).
$
\hfill (7.5)
\\[-2ex]

\newpage

The curves of family 
\\[2ex]
\mbox{}\hfill
$
\dfrac{(x-1)^2+y^2}{1-(x-1)^2-y^2}\, \exp\Bigl({}-2\arctan\dfrac{y}{x-1}\Bigr)=C,
\quad 
{}-\infty\leq C\leq {}+\infty,
\hfill
$
\\[2ex]
are trajectories of the differential system (7.5). 
\vspace{0.35ex}
Further, the equilibrium state $A(1,0)$ is an unstable focus of system (7.5) and 
\vspace{0.35ex}
the circle $(x-1)^2+y^2=1,$ which is passing through the origin of coordinates 
$O(0,0),$ is a stable limit cycle of system (7.5).
\vspace{0.35ex}

Trajectories of the stereographically conjugate system
\\[2ex]
\mbox{}\hfill                        
$
\dfrac{du}{d\theta}=
16u-12u^2+4v^2+2u^3-u^2v-2uv^2-v^3+\dfrac{1}{2}\,u^3v+\dfrac{1}{2}\,uv^3,   
\hfill
$
\\[0.5ex]
\mbox{}\hfill (7.6)
\\[0.5ex]
\mbox{}\hfill
$
\dfrac{dv}{d\theta}=
16v-16uv+u^3+4u^2v+uv^2-\dfrac{1}{4}\,u^4+\dfrac{1}{4}\,v^4,   
\hfill
$
\\[2.25ex]
where $(u^2+v^2)\;\!d\theta=dt,$ are the curves of family
\\[2ex]
\mbox{}\hfill
$
\dfrac{(u-4)^2+v^2}{u-2}\, \exp\Bigl(2\;\!\arctan\dfrac{4v}{(u-2)^2+v^2-4}\Bigr)=C^\ast,
\quad 
{}-\infty\leq C^\ast\leq {}+\infty,
\quad
 C^\ast=8C.
\hfill
$
\\[2ex]
\indent
Among these curves are the straight line $u=2,$ 
\vspace{0.25ex}
which is an open limit stereographic cycle,  
the unstable dicritical node $O^{\ast}(0,0),$ and  the unstable focus $A^{\ast}(4,0).$
\vspace{0.35ex}

The circles on Fig. 7.5 are stereographic atlases of trajectories for the differential systems (7.5) and (7.6).
\\[4.5ex]
\mbox{}\hfill
{\unitlength=1mm
\begin{picture}(45,45)
\put(0,0){\includegraphics[width=45mm,height=45mm]{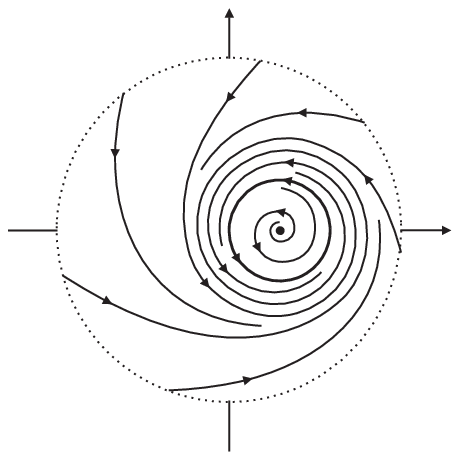}}
\put(44,20.5){\makebox(0,0)[cc]{\scriptsize $x$}}
\put(20.5,44){\makebox(0,0)[cc]{\scriptsize $y$}}
\end{picture}}
\qquad\qquad
{\unitlength=1mm
\begin{picture}(45,45)
\put(0,0){\includegraphics[width=45mm,height=45mm]{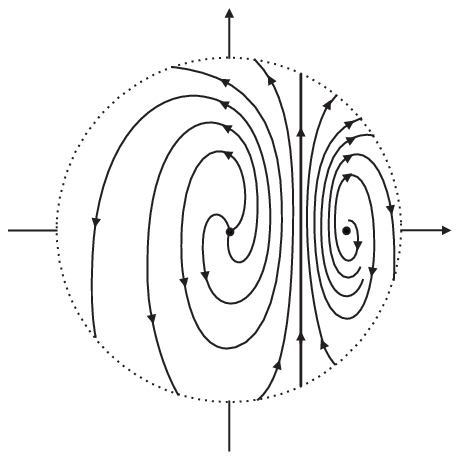}}
\put(44,20.5){\makebox(0,0)[cc]{\scriptsize $u$}}
\put(20.5,44){\makebox(0,0)[cc]{\scriptsize $v$}}
\end{picture}}
\hfill\mbox{}
\\[1.5ex]
\mbox{}\hfill
Fig. 7.5
\hfill\mbox{}
\\[4ex]
\indent
{\bf Example 7.4.}
\vspace{0.25ex}
If we replace $u$ by $u+2$ and $v$ by $v$ in the differential system (7.6), 
then we get the differential system 
\\[2ex]
\mbox{}\hfill                        
$
\dfrac{du}{d\theta}=
{}-8u+2uv+2u^3+2u^2v-2uv^2+\dfrac{1}{2}\,u^3v+\dfrac{1}{2}\,uv^3,   
\hfill
$
\\[0.75ex]
\mbox{}\hfill (7.7)
\\[0.75ex]
\mbox{}\hfill
$
\dfrac{dv}{d\theta}=
4+4u+2v^2-u^3+4u^2v+uv^2-\dfrac{1}{4}\,u^4+\dfrac{1}{4}\,v^4.   
\hfill
$
\\[2.5ex]
\indent
The curves of family 
\\[2ex]
\mbox{}\hfill
$
\dfrac{(u-2)^2+v^2}{u}\, \exp\Bigl(2\arctan\dfrac{4v}{u^2+v^2-4}\Bigr)=C^{\ast},
\quad 
{}-\infty\leq C^{\ast}\leq {}+\infty,
\hfill
$
\\[2.25ex]
are trajectories of system (7.7).\!
\vspace{0.5ex}
Moreover, the straight line $u\!=\!0$ is an open limit stereo- graphic cycle,  
$A^{}_1({}-2,0)$ is an unstable dicritical node, and 
$A^{}_2(2,0)$ is an unstable focus. 

\newpage

Trajectories of the stereographically conjugate system 
\\[2ex]
\mbox{}\hfill                        
$
\dfrac{dx}{dt}=
{}-32x-8xy+8x^3-8x^2y-8xy^2-2x^3y-2xy^3,   
\hfill
$
\\[0.5ex]
\mbox{}\hfill (7.8)
\\[0.5ex]
\mbox{}\hfill
$
\dfrac{dy}{dt}=
{}-16-16x-8y^2+4x^3+16x^2y-4xy^2+x^4-y^4,   
\hfill
$
\\[2.25ex]
where $(x^2+y^2)\;\!dt=d\theta,$ are the curves of family
\\[2ex]
\mbox{}\hfill
$
\dfrac{(x-2)^2+y^2}{x}\, \exp\Bigl(2\arctan\dfrac{4y}{4-x^2-y^2}\Bigr)=C^\ast,
\quad
{}-\infty\leq C^\ast\leq {}+\infty.
\hfill
$
\\[2ex]
\indent
Among these curves are the straight line $x=0,$ 
\vspace{0.35ex}
which is open limit stereographic cycle,  
the unstable dicritical node $A^{\ast}_1({}-2,0),$ and the unstable focus $A^{\ast}_2(2,0).$
\vspace{0.35ex}

The circles on Fig. 7.6 are stereographic atlases of trajectories for systems (7.7) and (7.8).
\\[3.25ex]
\mbox{}\hfill
{\unitlength=1mm
\begin{picture}(45,45)
\put(0,0){\includegraphics[width=45mm,height=45mm]{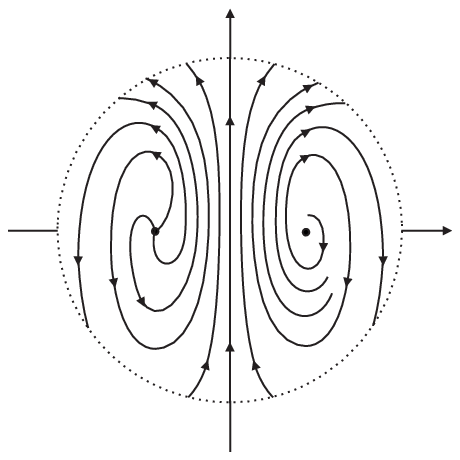}}
\put(44,20.5){\makebox(0,0)[cc]{\scriptsize $x$}}
\put(20.5,44){\makebox(0,0)[cc]{\scriptsize $y$}}
\end{picture}}
\qquad\qquad
{\unitlength=1mm
\begin{picture}(45,45)
\put(0,0){\includegraphics[width=45mm,height=45mm]{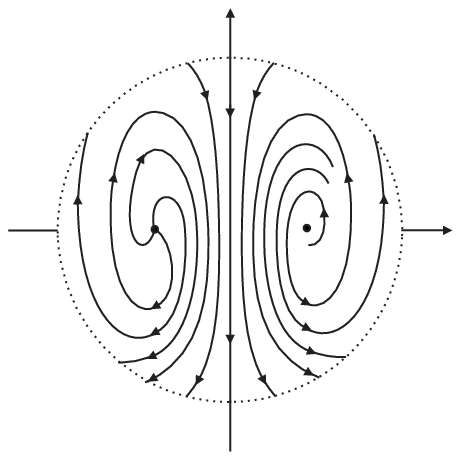}}
\put(44,20.5){\makebox(0,0)[cc]{\scriptsize $u$}}
\put(20.5,44){\makebox(0,0)[cc]{\scriptsize $v$}}
\end{picture}}
\hfill\mbox{}
\\[1ex]
\mbox{}\hfill
Fig. 7.6
\hfill\mbox{}
\\[6ex]
\centerline{\bf  8. Symmetry of phase directional field}
\\[0.35ex]
\centerline{\bf  for stereographically conjugate differential systems}
\\[1ex]
\indent
Using the analytical conditions of symmetry of phase directional field for differential system [11; 1], 
we obtain criteria of symmetry 
\vspace{0.5ex}
for stereographically conjugate differential systems.

{\bf Property 8.1.} 
{\it 
The following statements are equivalent\;\!{\rm:}
\vspace{0.35ex}

{\rm 1.} 
The phase directional field of system {\rm (D)}
\vspace{0.25ex}
is  symmetric with respect to the origin of coordinate of the phase plane $Oxy;$ 
\vspace{0.75ex}

{\rm 2.} 
\vspace{0.75ex}
The identity 
$
X(x,y)\;\!Y({}-x,{}-y)  - X({}-x,{}-y)\;\!Y(x,y)  = 0
$
for all 
$(x,y)\!\in\! \R^2$
is true{\rm;}

{\rm 3.} 
The phase directional field of system  {\rm (3.3)} 
\vspace{0.35ex}
is  symmetric with respect to the origin of coordinate of the phase plane $O^{\ast}uv;$ 
\vspace{0.75ex}

{\rm 4.} 
\vspace{0.75ex}
The identity
$
U(u,v)\;\!V({}-u,{}-v)  - U({}-u,{}-v)\;\!V(u,v)  = 0
$
for all 
$(u,v)\in \R^2$
is true.}

For example, 
the phase directional fields of stereographically conjugate differential systems  
(4.4) and (4.5) with $a_{_0}=b_{_0},$ 
\vspace{0.25ex}
(5.1) and (5.2), 
(5.3) and (5.4), 
(5.5) and (5.6), 
(5.7) and (5.8), 
(5.9) and (5.10),
(5.11) and (5.12),
(7.1) and (7.2), 
(7.3) and (7.4)
\vspace{0.75ex}
have such type of symmetry.

{\bf Property 8.2.}
{\it 
The following statements are equivalent\;\!{\rm:}
\vspace{0.35ex}

{\rm 1.} 
\vspace{0.25ex}
The phase directional field of the differential system {\rm (D)} is  symmetric with respect to the coordinate axis $Ox;$ 
\vspace{0.5ex}

{\rm 2.} 
\vspace{0.75ex}
The identity
$
X(x,y)\;\!Y(x,{}-y)  + X(x,{}-y)\;\!Y(x,y)  = 0
$
for all 
$(x,y)\in \R^2$
is true{\rm;}

{\rm 3.} 
The phase directional field of the differential system {\rm (3.3)} is symmetric
with respect to the coordinate axis $O^{\ast}\;\!\!u;$ 

{\rm 4.} 
The identity
\vspace{0.5ex}
$
U(u,v)\;\!V(u,{}-v)  + U(u,{}-v)\;\!V(u,v)  = 0
$
for all 
$(u,v)\in \R^2$ 
is true.}

For example, the phase directional fields of stereographically conjugate differential systems 
(4.4) and (4.5) with $a_{_0}b_{_0}=0,$ 
(5.1) and (5.2), 
(5.3) and (5.4), 
(5.7) and (5.8), 
(5.9) and (5.10)
\vspace{0.75ex}
have such type of symmetry.

{\bf Property 8.3.}
{\it 
The following statements are equivalent\;\!{\rm:}
\vspace{0.35ex}

{\rm 1.} 
The phase directional field of the differential system {\rm (D)} 
is  symmetric with respect to the coordinate axis $Oy;$ 
\vspace{0.5ex}

{\rm 2.} 
\vspace{0.75ex}
The identity
$
X(x,y)\;\!Y({}-x,y)  + X({}-x,y)\;\!Y(x,y)  = 0
$
for all 
$(x,y)\in \R^2$
is true{\rm;}

{\rm 3.} 
The phase directional field of the differential system {\rm (3.3)} is symmetric
with respect to the coordinate axis $O^{\ast}v;$ 
\vspace{0.5ex}

{\rm 4.} 
The identity
$
U(u,v)\;\!V({}-u,v)  + U({}-u,v)\;\!V(u,v)  = 0
$
\vspace{0.75ex}
for all $(u,v)\in \R^2$
is true.}

For example, 
\vspace{0.15ex}
the phase directional fields of stereographically conjugate differential systems
(4.4) and (4.5) with $a_{_0}b_{_0}=0,$ 
\vspace{0.25ex}
(5.1) and (5.2), 
(5.3) and (5.4), 
(5.7) and (5.8), 
(5.9) and (5.10)
\vspace{0.75ex}
have such type of symmetry.

{\bf Property 8.4.}
{\it 
The following statements are equivalent\;\!{\rm:}
\vspace{0.35ex}

{\rm 1.} 
The phase directional field of the differential system {\rm (D)} is symmetric
with respect to the straight line $y = x;$
\vspace{0.5ex}

{\rm 2.} 
\vspace{0.75ex}
The identity
$
X(x,y)\;\!X(y,x)  -  Y(x,y)\;\!Y(y,x)  =  0
$
for all 
$(x,y)\in \R^2$
is true{\rm;}

{\rm 3.} 
The phase directional field of the differential system {\rm (3.3)} is symmetric
with respect to the straight line $v = u;$
\vspace{0.5ex}

{\rm 4.} 
The identity
$
U(u,v)\;\!U(v,u)  -  V(u,v)\;\!V(v,u)  =  0
$
for all 
\vspace{0.5ex}
$(u,v)\in \R^2$
is true.}

For example, 
\vspace{0.25ex}
the phase directional fields of stereographically conjugate differential systems
(4.4) and (4.5) with $|a_{_0}|=|b_{_0}|,$ 
\vspace{0.35ex}
(5.1) and (5.2), 
(5.3) and (5.4), 
(5.7) and (5.8), 
(6.1) and (6.2), 
(6.3) and (6.4)
\vspace{0.75ex}
have such type of symmetry.

{\bf Property 8.5.}
{\it 
The following statements are equivalent\;\!{\rm:}
\vspace{0.35ex}

{\rm 1.} 
The phase directional field of the differential system {\rm (D)} is symmetric
with respect to the straight line $y = {}-x;$
\vspace{0.5ex}

{\rm 2.} 
\vspace{0.75ex}
The identity
$
X(\!{}-x,\!{}-y)\;\!X(y,x)  -  Y(\!{}-x,\!{}-y)\;\!Y(y,x)  =  0
$ 
for all 
$(x,y)\in \R^2$
is true{\rm;}

{\rm 3.} 
The phase directional field of the differential system {\rm (3.3)} is symmetric
with respect to the  straight line $v = {}-u;$
\vspace{0.5ex}

{\rm 4.} 
The identity
$
U({}-u,{}-v)\;\!U(v,u)  -  V({}-u,{}-v)\;\!V(v,u)  =  0
$ 
for all 
$(u,v)\in \R^2
$
\vspace{0.75ex}
is true.}

For example, 
\vspace{0.15ex}
the phase directional fields of stereographically conjugate differential systems
(4.4) and (4.5) with $|a_{_0}|=|b_{_0}|,$ 
(5.1) and (5.2), 
(5.3) and (5.4), 
(5.7) and (5.8)
have such type of symmetry. 
\\[3ex]
\centerline{\bf  9. Infinitely remote equilibrium state}
\\[1.5ex]
\indent
By Property 6.2, 
\vspace{0.25ex}
the infinitely remote point $M_\infty^{}$ of the extended phase plane $\overline{Oxy}$ is an equilibrium state of system (D) 
if and only if 
\vspace{0.25ex}
the point $O^\ast(0,0)$ is an equilibrium state of system (3.3). 
Note also that the equilibrium states $M_\infty^{}$ and $O^\ast$ have the same type.
\vspace{0.5ex}

Behaviour of trajectories for the differential system (D) 
\vspace{0.35ex}
in an neighbourhood of the infinitely remote point $M_\infty^{}$ 
of the extended phase plane $\overline{Oxy}$  
\vspace{0.25ex}
is defined by 
behaviour of trajectories for the differential system (D) in an neighbourhood 
\vspace{0.35ex}
of the infinitely remote straight line of the projective phase plane $\P\R(x,y)$ [12; 1].
\vspace{0.5ex}

Let $L$ be an infinitely remote equilibrium state of system (D) 
\vspace{0.25ex}
on the projective phase plane $\P\R(x,y).$ 
\vspace{0.25ex}
Then, consider two forms of Bendixon's sectors for the equilibrium state $L\colon$ interior form and exterior form. 
Points of infinitely remote straight line of the projective phase plane $\P\R(x,y),$ 
\vspace{0.35ex}
which are lying in a punctured neighborhood of the equilibrium state  $L,$ 
belong to exterior Bendixon's sectors  
and don't belong to interior Bendixon sectors. 
\vspace{0.25ex}

Thus we have the following 
\vspace{0.35ex}

{\bf Property 9.1.}
\vspace{0.15ex}
{\it Interior Bendixon's sector of the equilibrium state $L$ corresponds to Bendixon's sector of  the equilibrium state $M_\infty^{}$
with the same type and the same direction of movement along trajectories.}
\vspace{0.35ex}

Exterior Bendixon's sectors 
\vspace{0.15ex}
of infinitely remote equilibrium states of  the projective phase plane $\P\R(x, y)$ 
\vspace{0.35ex}
are collapsing  and forming Bendixon's sectors of the infinitely remote equilibrium state $M_\infty^{}$ of
the extended phase plane $\overline{Oxy}\,.$ 
\vspace{0.35ex}
In addition, note that [12; 1]:
if the system (D) is projectively nonsingular, then the infinitely remote straight line 
\vspace{0.25ex}
of the projective phase plane $\P\R(x, y)$ consists of trajectories of this system. 
\vspace{0.15ex}

Thus we have the following assertions. 
\vspace{0.15ex}

Suppose $L_1^{}$ and $L_2^{}$ are 
\vspace{0.35ex}
adjacent infinitely remove equilibrium states on the circle of  the projective circle  $\P\K(x,y)$ [12; 1] 
for the differential system (D). 
\vspace{0.5ex}

{\bf Property 9.2.}
\vspace{0.25ex}
{\it
If  the projectively nonsingular system {\rm (D)} has 
adjacent exterior Bendixon's sectors for  the equilibrium states $L_1^{}$ and $L_2^{}$
such that these sectors are{\rm:}

a{\rm)} parabolic{\rm;} 
\vspace{0.15ex}

b{\rm)} hy\-perbolic{\rm;} 
\vspace{0.15ex}

c{\rm)} one is parabolic, another is hyperbolic{\rm;} 
\vspace{0.15ex}

d{\rm)} one is hyperbolic, another is elliptic{\rm,} 
\\[0.5ex]
then these sectors are collapsing to Bendixon's sector of the type{\rm:} 
\vspace{0.15ex}

a{\rm)} elliptic; 
\vspace{0.15ex}

b{\rm)} hyperbolic; 
\vspace{0.15ex}

c{\rm)} parabolic; 
\vspace{0.15ex}

d{\rm)} elliptic
\\[0.5ex]
for the equilibrium state $M_\infty^{}.$}
\vspace{1.25ex}

{\bf  Example 9.1.}
The differential system [5, pp. 84 --- 85; 1]
\\[2ex]
\mbox{}\hfill        
$
\dfrac{dx}{dt}=1-x^2-y^2,
\qquad 
\dfrac{dy}{dt}=xy-1
$
\hfill(9.1)
\\[2.35ex]
on the projective phase plane $\P\R(x,y)$ has one equilibrium state (node), 
\vspace{0.15ex}
which is lying on  <<extremities>> of the axis $Ox.$
\\[3.5ex]
\mbox{}\hfill
{\unitlength=1mm
\begin{picture}(45,45)
\put(0,0){\includegraphics[width=45mm,height=45mm]{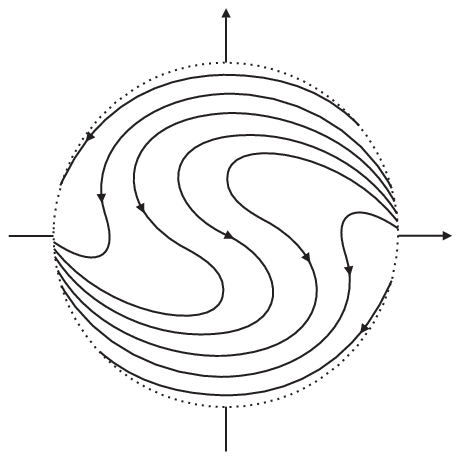}}
\put(44,20.5){\makebox(0,0)[cc]{\scriptsize $x$}}
\put(20.5,44){\makebox(0,0)[cc]{\scriptsize $y$}}
\end{picture}}
\qquad\qquad
{\unitlength=1mm
\begin{picture}(45,45)
\put(0,0){\includegraphics[width=45mm,height=45mm]{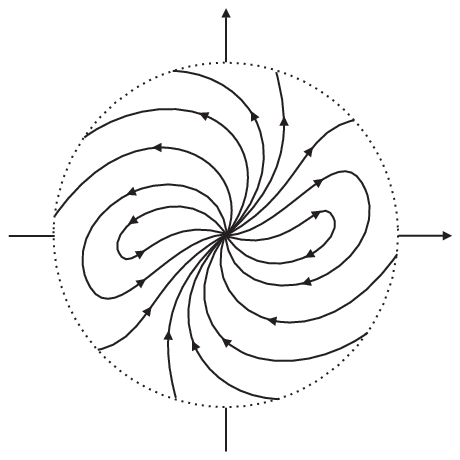}}
\put(44,20.5){\makebox(0,0)[cc]{\scriptsize $u$}}
\put(20.5,44){\makebox(0,0)[cc]{\scriptsize $v$}}
\end{picture}}
\hfill\mbox{}
\\[0ex]
\mbox{}\hfill
Fig. 9.1
\hfill\mbox{}
\\[4ex]
\indent
The equilibrium state $O^\ast(0,0)$ of the stereographically conjugate system
\\[2ex]
\mbox{}\hfill        
$
\dfrac{du}{d\theta}=4u^4-8u^2\;\!v^2-4v^4-\dfrac{1}{4}\,u^6+\dfrac{1}{2}\,u^5\;\!v-
\dfrac{1}{4}\,u^4\;\!v^2+u^3\;\!v^3+\dfrac{1}{4}\, u^2\;\!v^4+\dfrac{1}{2}\,uv^5+\dfrac{1}{4}\, v^6,
\hfill
$
\\[0.75ex]
\mbox{}\hfill (9.2)
\\[0.5ex]
\mbox{}\hfill
$
\dfrac{dv}{d\theta}=
12u^3\;\!v+4uv^3-\dfrac{1}{4}\,u^6-\dfrac{1}{2}\,u^5\;\!v-\dfrac{1}{4}\,u^4\;\!v^2-
u^3\;\!v^3+\dfrac{1}{4}\, u^2\;\!v^4-\dfrac{1}{2}\,uv^5+\dfrac{1}{4}\, v^6,
\hfill
$
\\[2.35ex]
where $(u^2+v^2)^2\;\!d\theta=dt,$
\vspace{0.35ex}
consists of two elliptic sectors, which are divided with the help of two parabolic sectors (Property 9.2, case  a).
\vspace{0.35ex}

Using the qualitative research [1] 
\vspace{0.15ex}
of behaviour of trajectories for system (9.1), 
we get the circles on Fig. 9.1 compose stereografic atlases of trajectories for systems (9.1) and (9.2). 
\vspace{1.25ex}

{\bf  Example 9.2.}
The differential system [13, pp. 61 --- 65]
\\[2ex]
\mbox{}\hfill        
$
2\,\dfrac{dx}{dt}=2y+i\;\!(x-i\;\!y)^q-i\;\!(x+i\;\!y)^q,
\qquad 
2\,\dfrac{dy}{dt}={}-2x+(x-i\;\!y)^q+ (x+i\;\!y)^q
$
\hfill (9.3)
\\[2.25ex]
at $i=\sqrt{{}-1}\,,\ q=4$ and $q=5$  
is projectively nonsingular [12]. 
\vspace{0.5ex}

Moreover, 
\vspace{0.35ex}
all infinitely remove equilibrium states on the projective phase plane $\P\R(x,y)$ are nodes 
(see  Fig. 2.12 in [13, p. 65]  or Fig. 8.9 and Fig 8.10 in [12]).
\vspace{0.75ex}

If $q=4,$ then the differential system (9.3) has the form
\\[2ex]
\mbox{}\hfill        
$
\dfrac{dx}{dt}=y+4x^3y-4xy^3,
\qquad 
\dfrac{dy}{dt}={}-x+x^4-6x^2y^2+y^4
$
\hfill (9.4)
\\[2.5ex]
and the stereographically conjugate system  to the system (9.4) is the differential system
\\[2.25ex]
\mbox{}\qquad\qquad        
$
\dfrac{du}{d\theta}=
v({}-384u^5+1280u^3v^2-384uv^4+
u^8+4u^6v^2+6u^4v^4+4u^2v^6+v^8),
\hfill
$
\\[0.5ex]
\mbox{}\hfill (9.5)
\\[0.5ex]
\mbox{}\qquad\qquad
$
\dfrac{dv}{d\theta}=
64u^6-960u^4v^2+960u^2v^4-64v^6
-u^9-4u^7v^2-6u^5v^4-4u^3v^6-uv^8,
\hfill
$
\\[2.5ex]
where $(u^2+v^2)^4\;\!d\theta=dt.$
\\[4ex]
\mbox{}\hfill
{\unitlength=1mm
\begin{picture}(45,45)
\put(0,0){\includegraphics[width=45mm,height=45mm]{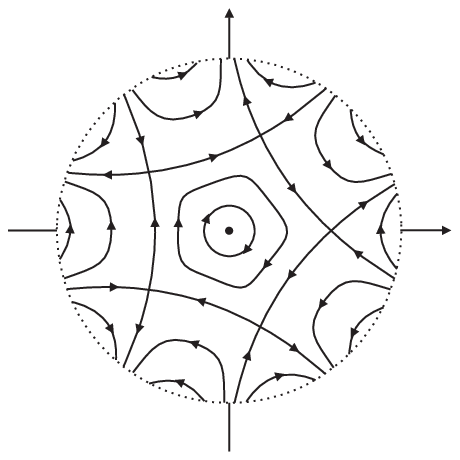}}
\put(44,20.5){\makebox(0,0)[cc]{\scriptsize $x$}}
\put(20.5,44){\makebox(0,0)[cc]{\scriptsize $y$}}
\end{picture}}
\qquad\qquad\qquad
{\unitlength=1mm
\begin{picture}(45,45)
\put(0,0){\includegraphics[width=45mm,height=45mm]{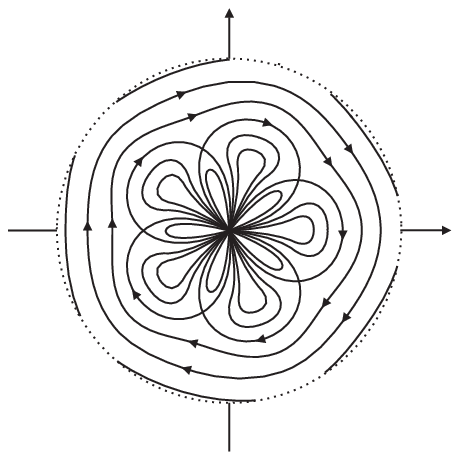}}
\put(44,20.5){\makebox(0,0)[cc]{\scriptsize $u$}}
\put(20.5,44){\makebox(0,0)[cc]{\scriptsize $v$}}
\end{picture}}
\hfill\mbox{}
\\[2.25ex]
\mbox{}\hfill
Fig. 9.2
\hfill\mbox{}
\\[7ex]
\mbox{}\hfill
{\unitlength=1mm
\begin{picture}(45,45)
\put(0,0){\includegraphics[width=45mm,height=45mm]{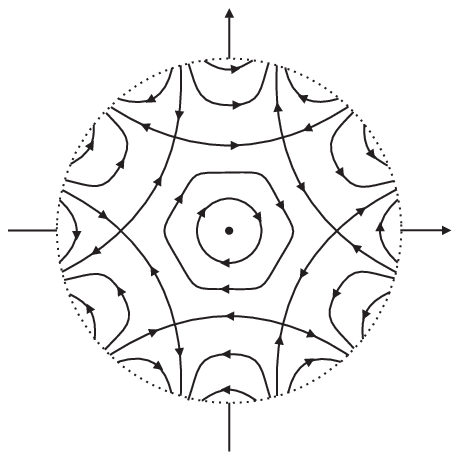}}
\put(44,20.5){\makebox(0,0)[cc]{\scriptsize $x$}}
\put(20.5,44){\makebox(0,0)[cc]{\scriptsize $y$}}
\end{picture}}
\qquad\qquad\qquad
{\unitlength=1mm
\begin{picture}(45,45)
\put(0,0){\includegraphics[width=45mm,height=45mm]{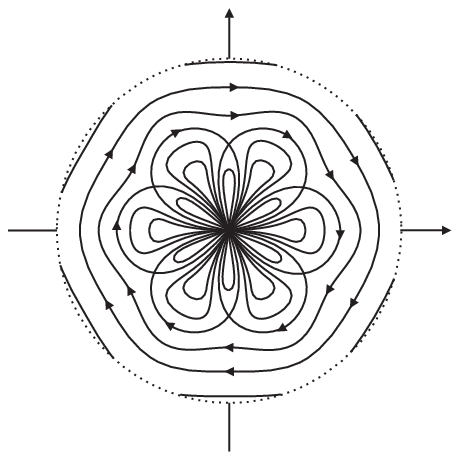}}
\put(44,20.5){\makebox(0,0)[cc]{\scriptsize $u$}}
\put(20.5,44){\makebox(0,0)[cc]{\scriptsize $v$}}
\end{picture}}
\hfill\mbox{}
\\[2.25ex]
\mbox{}\hfill
Fig. 9.3
\hfill\mbox{}
\\[-2ex]

\newpage

If $q=5,$ then the differential system (9.3) has the form
\\[2ex]
\mbox{}\hfill        
$
\dfrac{dx}{dt}=y+5x^4\;\!y-10x^2\;\!y^3+y^5,
\qquad 
\dfrac{dy}{dt}={}-x+x^5-10x^3\;\!y^2+5xy^4
$
\hfill (9.6)
\\[2.5ex]
and stereographically conjugate system to the system (9.6) is the differential system
\\[2.25ex]
\mbox{}        
$
\dfrac{du}{d\theta}=
v({}-1792u^6+8960u^4v^2-5376u^2v^4+256v^6+
u^{10}+5u^8v^2+10u^6v^4+10u^4v^6+5u^2v^8+v^{10}),
\ \ 
\hfill
$
\\
\mbox{}\hfill (9.7)
\\
\mbox{}
$
\dfrac{dv}{d\theta}=
u(256u^6-5376u^4v^2+8960u^2v^4-1792v^6
-u^{10}-5u^8v^2-10u^6v^4-10u^4v^6-5u^2v^8-v^{10}),
\ \ 
\hfill
$
\\[2.5ex]
where $(u^2+v^2)^5\;\!d\theta=dt.$
\vspace{0.5ex}

By Property 9.2, the equilibrium state $O^\ast(0,0)$ both for the differential system (9.5) 
\vspace{0.35ex}
and for the differential system (9.7) consists of elliptic Bendixon's sectors.
\vspace{0.35ex}

The circles, 
\vspace{0.15ex}
which are forming stereografic atlases of trajectories for the differential systems (9.4) and (9.5), are constructed on Fig. 9.2,
\vspace{0.15ex}
and the circles, which are forming stereografic atlases of trajectories for the differential systems (9.6) and (9.7), 
\vspace{1.25ex}
are constructed on Fig. 9.3.

{\bf Example 9.3.}
The differential system [5, pp. 85 --- 87; 14, pp. 209 --- 212]
\\[2ex]
\mbox{}\hfill        
$
\dfrac{dx}{dt}={}-1+x^2+y^2,
\qquad 
\dfrac{dy}{dt}={}-5+5xy
$
\hfill(9.8)
\\[2.25ex]
on the projective phase plane $\P\R(x,y)$ has three equilibrium states:
\vspace{0.25ex}
a saddle and two stable nodes.
The saddle is lying on the <<extremities>> of the straight line $y=0.$
\vspace{0.25ex}
The first stable node lies on the <<extremities>> of the straight line $y={}-2x$
\vspace{0.25ex}
and the second stable node lies on the <<extremities>> of the straight line $y=2x.$
\vspace{0.5ex}

The stereographically conjugate differential system
\\[2.75ex]
\mbox{}\hfill        
$
\dfrac{du}{d\theta}={}-4u^4-40u^2v^2+4v^4+\dfrac{1}{4}\,u^6+\dfrac{5}{2}\,u^5v+\dfrac{1}{4}\,u^4v^2+
5u^3v^3-\dfrac{1}{4}\, u^2v^4+\dfrac{5}{2}\,uv^5-\dfrac{1}{4}\, v^6,
\hfill
$
\\[0.5ex]
\mbox{}\hfill (9.9)
\\[0.5ex]
\mbox{}\hfill
$
\dfrac{dv}{d\theta}=
12u^3v-28uv^3-\dfrac{5}{4}\,u^6+\dfrac{1}{2}\,u^5v-\dfrac{5}{4}\,u^4v^2+
u^3v^3+\dfrac{5}{4}\, u^2v^4+\dfrac{1}{2}\,uv^5+\dfrac{5}{4}\, v^6,
\hfill
$
\\[2.75ex]
where $(u^2+v^2)^2\;\!d\theta=dt,$ 
\vspace{0.35ex}
has the equilibrium state $O^\ast(0,0),$ 
which is consisting of two elliptic Bendixon's sectors and two parabolic Bendixon's sectors 
\vspace{0.5ex}
(Property 9.2, cases a, c).

Using the qualitative research [1] of behaviour of trajectories 
\vspace{0.25ex}
for the differential system (9.8), we can build the circles on Fig. 9.4, 
\vspace{0.15ex}
which are forming stereografic atlases of trajectories 
for the differential systems (9.8) and (9.9).
\\[3.75ex]
\mbox{}\hfill
{\unitlength=1mm
\begin{picture}(45,45)
\put(0,0){\includegraphics[width=45mm,height=45mm]{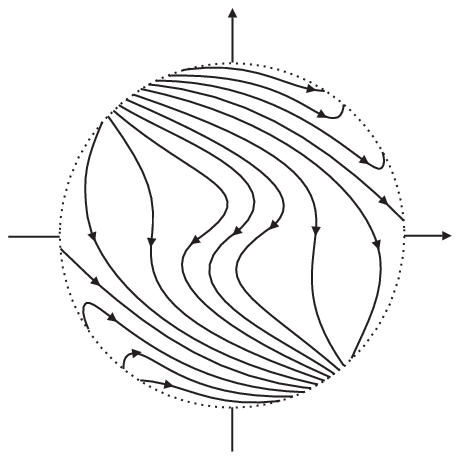}}
\put(44,20.5){\makebox(0,0)[cc]{\scriptsize $x$}}
\put(20.5,44){\makebox(0,0)[cc]{\scriptsize $y$}}
\end{picture}}
\qquad\qquad\qquad
{\unitlength=1mm
\begin{picture}(45,45)
\put(0,0){\includegraphics[width=45mm,height=45mm]{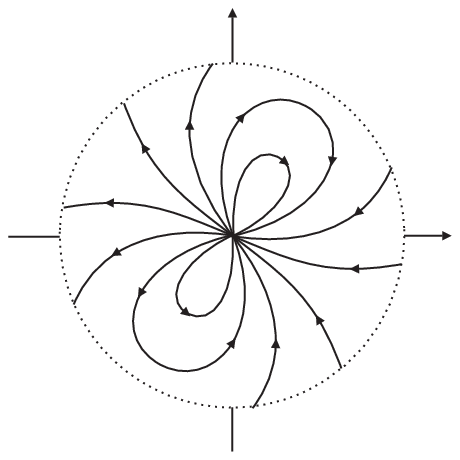}}
\put(44,20.5){\makebox(0,0)[cc]{\scriptsize $u$}}
\put(20.5,44){\makebox(0,0)[cc]{\scriptsize $v$}}
\end{picture}}
\hfill\mbox{}
\\[2ex]
\mbox{}\hfill
Fig. 9.4
\hfill\mbox{}
\\[-2ex]

\newpage

{\bf Property 9.3.}
\vspace{0.25ex}
{\it
If the boundary circle of the projective circle $\P\K(x,y)$ is orthogonal intersected by each trajectory 
of the projectively singular differential system {\rm (D)}, 
\vspace{0.25ex}
then the infinitely remote equilibrium state $M_\infty^{} $ is a dicritical node.}
\vspace{0.5ex}

For example,  such differential system is Darboux's system (7.3). 
\vspace{0.15ex}
The projective atlas of trajectories for this system
is constructed on Fig. 16.3 in [1], and the stereographic atlas of trajectories for this system is constructed on Fig. 7.4.
\vspace{0.5ex}

{\bf  Example 9.4.} 
\vspace{0.25ex}
In the papers [11; 1], we did the qualitative research of behaviour of tra\-jec\-to\-ries on the projective phase plane $\P\R(x,y)$
for the differential systems
\\[1.5ex]
\mbox{}\hfill        
$
\dfrac{dx}{dt}={}-y+x^3,
\qquad 
\dfrac{dy}{dt}=x(1+xy),
$
\hfill(9.10)
\\[2.75ex]
\mbox{}\hfill                                       
$
\begin{array}{l}
\dfrac{dx}{dt}=x(x^2+y^2-1)(x^2+y^2-9)-y(x^2+y^2-2x-8),
\\[3ex]
\dfrac{dy}{dt}=y(x^2+y^2-1)(x^2+y^2-9)+x(x^2+y^2-2x-8),
\end{array}
$
\hfill (9.11)
\\[2.75ex]
\mbox{}\hfill                                       
$
\begin{array}{l}
\dfrac{dx}{dt}=x(2x^2+2y^2+1)\Bigl((x^2+y^2)^2+x^2-y^2+\dfrac{1}{10}\Bigr)-y(2x^2+2y^2-1),
\\[3ex]
\dfrac{dy}{dt}=y(2x^2+2y^2-1)\Bigl((x^2+y^2)^2+x^2-y^2+\dfrac{1}{10}\Bigr)+x(2x^2+2y^2+1),
\end{array}
$
\hfill (9.12)
\\[2ex]
and also the projective atlases of trajectories for these systems were built. 

Trajectories of each of these systems intersect orthogonally the
boundary circle of the projective circle $\P\K(x, y).$  

The stereographically conjugate systems to the differential systems (9.10), (9.11), (9.12) 
accordingly are the differential systems
\\[1.75ex]
\mbox{}\hfill        
$
\dfrac{du}{d\theta}={}-(16u^3+u^4v+2u^2v^3+v^5),
\qquad 
\dfrac{dv}{d\theta}=u({}-16uv+u^4+2u^2v^2+v^4),
$
\hfill(9.13)
\\[2.25ex]
where $(u^2+v^2)^2\;\!d\theta=dt,$
\\[2.25ex]
\mbox{}\hfill                                       
$
\begin{array}{l}
\dfrac{du}{d\theta}=
{}-256u +160u^3-16u^2v+160uv^2-16v^3
+8u^3v+8uv^3\ -
\\[2.5ex]
\mbox{}\qquad \ \ 
-\ 9u^5 +8u^4v-18u^3v^2+16u^2v^3-9uv^4+8v^5,
\\[2.75ex]
\dfrac{dv}{d\theta}=
{}-256v +16u^3+160u^2v +16uv^2+160v^3
-8u^4-8u^2v^2\ -
\\[2.5ex]
\mbox{}\qquad \ \
- \ 8u^5-9u^4v-16u^3v^2-18u^2v^3-8uv^4-9v^5,
\end{array}
$
\hfill (9.14)
\\[2.25ex]
where $(u^2+v^2)^2\;\!d\theta=dt,$
\\[2.25ex]
\mbox{}                                       
$
\dfrac{du}{d\theta}=
{}-8192u-768u^3+1280uv^2
-\dfrac{96}{5}\,u^5-32u^4v+\dfrac{288}{5}\,u^3v^2-64u^2v^3 -\
\hfill
$
\\[2.5ex]
\mbox{}\qquad
$
-\ \dfrac{256}{5}\,uv^4-32v^5 - \dfrac{1}{10}\,u^7-3u^6v+
\dfrac{1}{10}\,u^5v^2-5u^4v^3+\dfrac{1}{2}\,u^3v^4-u^2v^5+\dfrac{3}{10}\,uv^6+v^7,
\hfill
$
\\[0.5ex]
\mbox{}\hfill (9.15)
\\[0.5ex]
\mbox{}
$
\dfrac{dv}{d\theta}=
{}-8192v-1280u^2v+768v^3
+32u^5-\dfrac{256}{5}\,u^4v+64u^3v^2+\dfrac{288}{5}\,u^2v^3\ +
\hfill
$
\\[2.5ex]
\mbox{}\qquad
$
+\ 32uv^4-\dfrac{96}{5}\,v^5 +
u^7-\dfrac{3}{10}\,u^6v-u^5v^2-\dfrac{1}{2}\,u^4v^3-5u^3v^4-\dfrac{1}{10}\,u^2v^5
-3uv^6+\dfrac{1}{10}\,v^7,
\hfill
$
\\[2.35ex]
where $(u^2+v^2)^3\;\!d\theta=dt.$
\vspace{0.35ex}

The circles, which are forming stereografic atlases of trajectories for the differential systems (9.10) and (9.13), are constructed on Fig. 9.5.
The circles, which are forming  stereografic atlases of trajectories for the differential systems (9.11) and (9.14), are constructed on 
Fig.~9.6. And the circles, which are forming  stereografic atlases 
of trajectories for the dif\-fe\-ren\-ti\-al systems (9.12) and (9.15), are constructed on Fig. 9.7.      
\\[4.75ex]
\mbox{}\hfill
{\unitlength=1mm
\begin{picture}(45,45)
\put(0,0){\includegraphics[width=45mm,height=45mm]{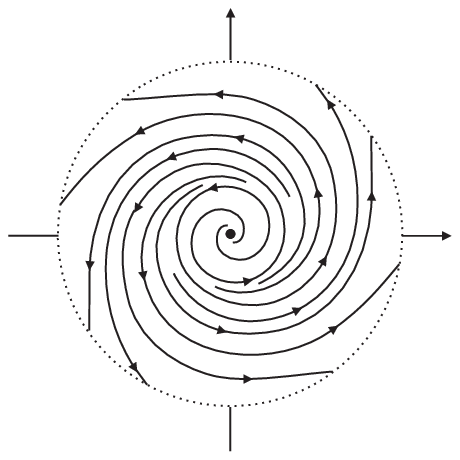}}
\put(44,20.5){\makebox(0,0)[cc]{\scriptsize $x$}}
\put(20.5,44){\makebox(0,0)[cc]{\scriptsize $y$}}
\end{picture}}
\qquad\qquad\qquad
{\unitlength=1mm
\begin{picture}(45,45)
\put(0,0){\includegraphics[width=45mm,height=45mm]{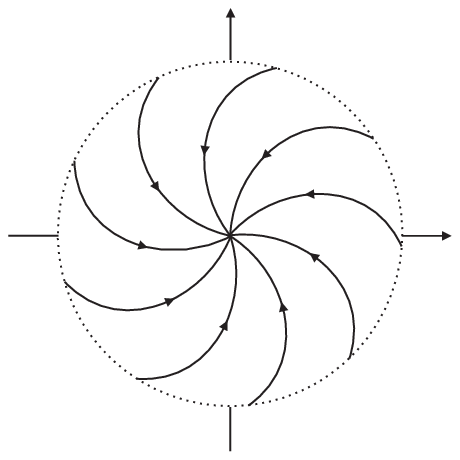}}
\put(44,20.5){\makebox(0,0)[cc]{\scriptsize $u$}}
\put(20.5,44){\makebox(0,0)[cc]{\scriptsize $v$}}
\end{picture}}
\hfill\mbox{}
\\[1ex]
\mbox{}\hfill
Fig. 9.5
\hfill\mbox{}
\\[9ex]
\mbox{}\hfill
{\unitlength=1mm
\begin{picture}(45,45)
\put(0,0){\includegraphics[width=45mm,height=45mm]{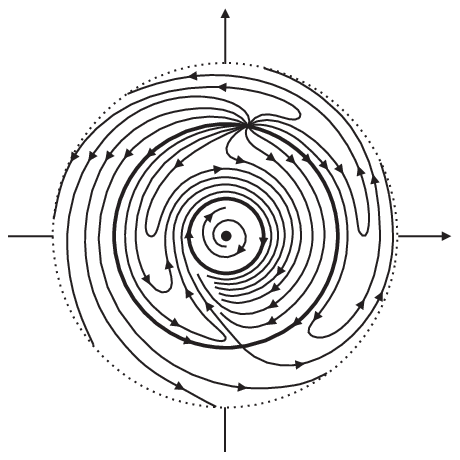}}
\put(44,20.5){\makebox(0,0)[cc]{\scriptsize $x$}}
\put(20.5,44){\makebox(0,0)[cc]{\scriptsize $y$}}
\end{picture}}
\qquad\qquad\qquad
{\unitlength=1mm
\begin{picture}(45,45)
\put(0,0){\includegraphics[width=45mm,height=45mm]{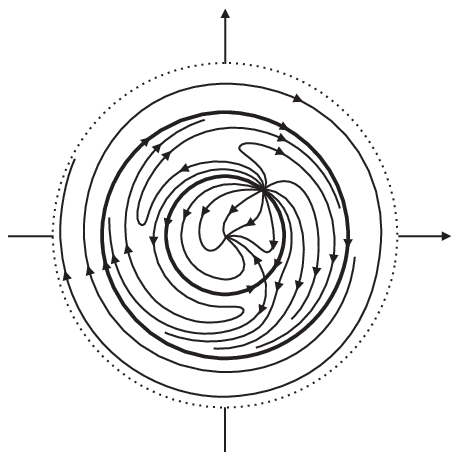}}
\put(44,20.5){\makebox(0,0)[cc]{\scriptsize $u$}}
\put(20.5,44){\makebox(0,0)[cc]{\scriptsize $v$}}
\end{picture}}
\hfill\mbox{}
\\[1ex]
\mbox{}\hfill
Fig. 9.6
\hfill\mbox{}
\\[9ex]
\mbox{}\hfill
{\unitlength=1mm
\begin{picture}(45,45)
\put(0,0){\includegraphics[width=45mm,height=45mm]{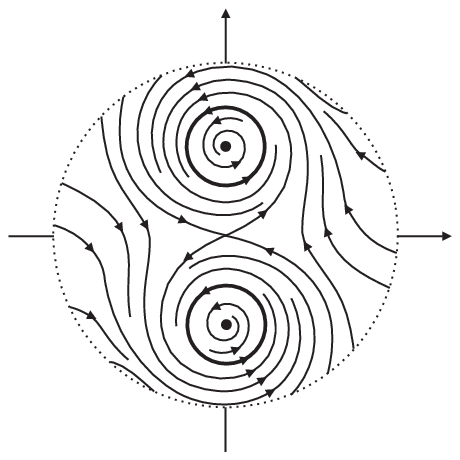}}
\put(44,20.5){\makebox(0,0)[cc]{\scriptsize $x$}}
\put(20.5,44){\makebox(0,0)[cc]{\scriptsize $y$}}
\end{picture}}
\qquad\qquad\qquad
{\unitlength=1mm
\begin{picture}(45,45)
\put(0,0){\includegraphics[width=45mm,height=45mm]{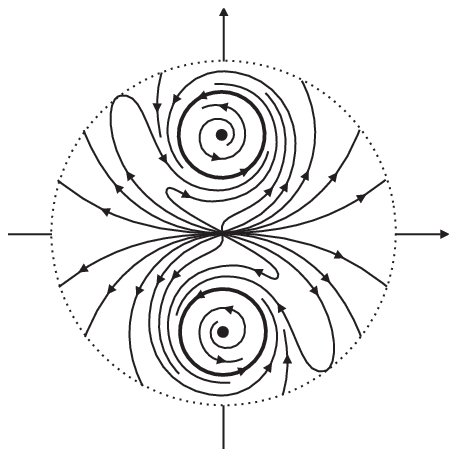}}
\put(44,20.5){\makebox(0,0)[cc]{\scriptsize $u$}}
\put(20.5,44){\makebox(0,0)[cc]{\scriptsize $v$}}
\end{picture}}
\hfill\mbox{}
\\[1ex]
\mbox{}\hfill
Fig. 9.7
\hfill\mbox{}
\\[-2ex]

\newpage

\mbox{}
\\[-2.25ex]

{\Large\bf  References}
\vspace{1.75ex}

1. 
V.N. Gorbuzov, 
Projective atlas of trajectories of differential systems,   
{\it Mathematics. Dy\-na\-mi\-cal Sys\-tems} (arXiv: 1401.1000v1 [math.DS].  Cornell Univ., 
Ithaca, New York), 
\linebreak
2014,  61 p.
\vspace{0.75ex}

2. 
V.N. Gorbuzov, 
Stereographic atlas of trajectories for differential systems of the second order  (Russian),
{\it Vestnik of the Yanka Kupala Grodno State University}, 2014, Ser. 2,  
\linebreak
No. 1(170), 12-20. 
\vspace{0.75ex}

3. 
V.N. Gorbuzov, 
Trajectories of stereographic conjugate differential systems  (Russian),
{\it Vestnik of the Yanka Kupala Grodno State University}, 2014, Ser. 2,  No. 3(180), 27-36. 
\vspace{0.75ex}

4.
M.A. Lavrentiev  and B.V. Shabat, 
{\it Methods of the theory of functions of a complex variable} (Russian),
Nauka, Moscow, 1987.
\vspace{0.75ex}

5.
{\it The mathematical encyclopaedia} (Russian), Vol. 5, 
\vspace{0.75ex}
Soviet encyclopaedia, Moscow, 1984.

6.
А.S. Mishchenko and А.Т. Fomenko,
{\it A course of differential geometry and topology}  (Russian), 
Moscow State University, Moscow, 1980.
\vspace{0.75ex}

7.
A.A. Andronov, E.A. Leontovich, I.I. Gordon, and A.G. Maier, 
{\it Qualitative theory of second order dynamical systems} (Russian),
Nauka, Moscow, 1966.
\vspace{0.75ex}

8.
H. Poincar\'{e},
{\it On curves defined by differential equations} (Russian),
GITTL, Moscow-Leningrad, 1947.
\vspace{0.75ex}

9.
V.N. Gorbuzov and P.B. Pavlyuchik, 
 Linear and open limit cycles of differential systems (Russian),
{\it Vestnik of the Yanka Kupala Grodno State University}, 2013, Ser. 2,  
\linebreak
No. 3(159), 23-32. 
\vspace{0.75ex}

10.
V.N. Gorbuzov and А.А. Samodurov, 
{\it The Darboux equation and its analogs} (Russian),
Grodno State University, Grodno, 1985. 
\vspace{0.75ex}

11.
V.N. Gorbuzov, 
Trajectories of projective reduced differential systems (Russian),
{\it Vestnik of the Yanka Kupala Grodno State University}, 2012, Ser. 2,  No. 1(126), 39-52. 
\vspace{0.75ex}

12.
V.N. Gorbuzov, 
The projective atlas of trajectories of differential systems of the second order (Russian), 
{\it Vestnik of the Yanka Kupala Grodno State University}, 2011, Ser. 2,  
\linebreak
No. 2(111), 15-26. 
\vspace{0.75ex}

13.
K.S. Sibirsky, 
{\it  
Algebraic invariants of differential equations and matrixes} (Russian), 
Shtiinca, Chisinau, 1976.
\vspace{0.75ex}

14.
 S. Lefschetz,
{\it Differential equations{\rm:} geometric theory}  (Russian), IL, Moscow, 1961.
\vspace{0.75ex}

}

\newpage

{\normalsize

\pagenumbering{arabic}

\sloppy

\lhead
    [\scriptsize В.Н. Горбузов]
    {\scriptsize В.Н. Горбузов}
\rhead
    [\it \scriptsize Стереографически сопряженные дифференциальные системы]
    {\it \scriptsize Стереографически сопряженные дифференциальные системы}

\thispagestyle{empty}

\mbox{}
\\[-0.15ex]
\centerline{
{\large
\bf
СТЕРЕОГРАФИЧЕСКИ\;\! СОПРЯЖЕННЫЕ 
}}
\\[0.5ex]
\centerline{
{\large
\bf
ДИФФЕРЕНЦИАЛЬНЫЕ\;\! СИСТЕМЫ
}
}
\\[2.25ex]
\centerline{
\bf 
В.Н. Горбузов
}
\\[2ex]
\centerline{
\it 
Факультет математики и информатики,
}
\\[0.5ex]
\centerline{
\it 
Гродненский государственный университет имени Янки Купалы,
}
\\[0.5ex]
\centerline{
\it 
Ожешко {\rm22}, Гродно, Беларусь, {\rm 230023}
}
\\[1.5ex]
\centerline{
E-mail: gorbuzov@grsu.by
}
\\[4.25ex]
\centerline{{\large\bf Резюме}}
\\[1ex]
\indent
Изложены топологические основы поведения траекторий автономных дифференциальных систем второго 
порядка на сфере. Построен стереографический атлас траекторий. 
Установлены дифференциальные связи между траекториями стереографически сопряженных 
дифференциальных систем. Исследовано поведение траекторий в окрестности бесконечно удаленной 
точки фазовой плоскости.  
Приведены примеры глобального качественного исследования траекторий стереографически сопряженных
дифференциальных систем. 
\\[1.5ex]
\indent
{\it Ключевые слова}:
дифференциальная система, стереографическая проекция, атлас карт многообразия.
\\[1.25ex]
\indent
{\it 2000 Mathematics Subject Classification}: 34A26, 34C05.
\\[4.25ex]
\centerline{{\large\bf Содержание}}
\\[1.25ex]
{\bf  Введение}                   \dotfill\ 2
\\[1ex]
{\bf \S 1. 
Стереографический атлас сферы}
                                                 \dotfill \ 2
\\[0.75ex]
\mbox{}\hspace{1.35em}
1. Стереографическая проекция плоскости
                                                 \dotfill \ 2
\\[0.5ex]
\mbox{}\hspace{1.35em}
2. Стереографический атлас сферы
                                                 \dotfill \ 4
\\[1ex]
\noindent
{\bf \S 2. 
Стереографически сопряженная дифференциальная система
}
                                                 \dotfill \ 6
\\[0.75ex]
\mbox{}\hspace{1.35em}
3. Преобразование Бендиксона
                                                 \dotfill \ 6
\\[0.5ex]
\mbox{}\hspace{1.35em}
4. Вид стереографически сопряженной дифференциальной системы 
                                                 \dotfill \ 7
\\[0.5ex]
\mbox{}\hspace{1.35em}
5. Стереографический атлас траекторий дифференциальной системы
                                                 \dotfill \ 10
\\[1ex]
\noindent
{\bf \S 3. 
Траектории стереографически сопряженных 
\\
\mbox{}\hspace{1.25em}
дифференциальных систем
}
                                                 \dotfill \ 16
\\[0.75ex]
\mbox{}\hspace{1.35em}
6. Регулярные точки и состояния равновесия
стереографически сопряженных 
\\
\mbox{}\hspace{2.6em}
 дифференциальных систем
                                                 \dotfill \ 16
\\[0.5ex]
\mbox{}\hspace{1.35em}
7. Стереографические циклы
                                                 \dotfill \ 20
\\[0.5ex]
\mbox{}\hspace{1.35em}
8. Симметpичность фазового поля направлений
стереографически сопряженых
\\
\mbox{}\hspace{2.6em}
дифференциальных систем
                                                 \dotfill \ 24
\\[0.5ex]
\mbox{}\hspace{1.35em}
9. Бесконечно удаленное состояние равновесия
                                                 \dotfill \ 26
\\[1ex]
{\bf Список литературы}
                                              \dotfill \ 31

\newpage

\mbox{}
\\[-1.75ex]
\centerline{\large\bf  Введение}
\\[1.5ex]
\indent
Объектом исследования является обыкновенная автономная полиномиальная
дифференциальная система второго порядка
\\[2ex]
\mbox{}\hfill                                   
$
\displaystyle 
\dfrac{dx}{dt} =
\sum \limits_{k=0}^{n}\, X_k^{}(x,y)\equiv
 X(x,y), 
\qquad 
\dfrac{dy}{dt} =
\sum \limits_{k=0}^{n}\,Y_k^{}(x,y)\equiv
Y(x,y),
$
\hfill (D)
\\[2.25ex]
где 
\vspace{0.5ex}
$X_k^{}$ и $Y_k^{}$ --- однородные полиномы по переменным $x$ и $y$ степени $k,\ k=0,1,\ldots,n,$ такие, что 
\vspace{0.35ex}
$|X_n^{}(x,y)|+|Y_n^{}(x,y)|\not\equiv0$ на $\R^2,$ а полиномы 
$X$ и $Y$ --- взаимно простые, т.е. не имеют общих делителей, отличных от вещественных чисел. 

Данная работа является продолжением исследований, изложенных в [1], а ее основные результаты 
опубликованы в статьях [2] и [3].
\\[5ex]
\centerline{
{\bf\large \S\;\!1. Стереографический атлас сферы}}
\\[2ex]
\centerline{
{\bf  1. 
Стереографическая проекция плоскости
}
}
\\[1.5ex]
\indent
Введем трехмерную прямоугольную декартову систему координат
\vspace{0.15ex}
$O^{\star}x^{\star}y^{\star}z^{\star},$ 
совмещенную с правой прямоугольной декартовой системой координат $Oxy,$
\vspace{0.15ex}
соблюдая условия:
прямая $OO^{\star}$ ортогональна плоскости $Oxy,$
\vspace{0.15ex}
длина отрезка $OO^{\star}$ равна одной единице масштаба системы
координат $Oxy;$
\vspace{0.15ex}
ось $O^{\star}x^{\star}$ сонаправлена с осью $Ox,$ ось
$O^{\star}y^{\star}$ сонаправлена с осью $Oy,$ а
ось $O^{\star}z^{\star}$ направлена так, что система
координат $O^{\star}x^{\star}y^{\star}z^{\star}$ будет правой;
масштаб в системе координат $O^{\star}x^{\star}y^{\star}z^{\star}$
\vspace{0.15ex}
такой же, как и в системе координат $Oxy.$
Построим сферу с центром $O^{\star}$
единичного радиуса: 
\\[2ex]
\mbox{}\hfill                                   
$
S^2=\bigl\{
(x^{\star},y^{\star},z^{\star})
\colon
x^{\star}{}^{\,^{\scriptstyle 2}}+
y^{\star}{}^{\,^{\scriptstyle 2}}+
z^{\star}{}^{\,^{\scriptstyle 2}}=1\bigr\}.
$
\hfill (1.1)
\\[2ex]
\noindent
Точки 
\vspace{0.15ex}
$N(0,0,1)$ и $S(0,0,{}-1)$ ---  соответственно северный и южный полюсы этой сферы. 
При этом южный полюс $S(0,0,{}-1)$ 
\vspace{0.15ex}
совпадает с началом $O(0,0)$ системы координат $Oxy.$
Уравнение $z^{\star}={}-1 $
\vspace{0.15ex}
является уравнением  в системе координат $O^{\star}x^{\star}y^{\star}z^{\star}$
плоскости $Oxy.$
Плоскость $Oxy$ касается
сферы (1.1) в южном полюсе $S(0,0,{}-1).$
\vspace{0.25ex}

На плоскости $Oxy$  произвольным образом
\vspace{0.15ex}
выберем точку $M(x,y)$ и проведем луч с
на\-ча\-лом в точке $M$ через северный полюс $N.$ Луч $MN$ пересекает
сферу (1.1) в неко\-то\-рой точке $P.$ Тем самым, каждой точке 
пло\-с\-кос\-ти $Oxy$ сопоставляется одна точка сфе\-ры, а каждой точке
сферы, отличной от северного полюса, сопоставляется, одна точка
плоскости $Oxy.$
\vspace{0.15ex}
Такую проекцию назовем
{\it стереографической проекцией} плоскости на сфе\-ру [4, c. 83 --- 84].
\vspace{0.15ex}
Точку $N$ будем называть
{\it центром} стереографической проекции (рис. 1.1).
\vspace{0.15ex}

{\bf Лемма 1.1.}
{\it
Стереографическая проекция 
плоскости является биекцией между плоскостью и сферой с
выколотым северным полюсом --- центром этой проекции.}
\vspace{0.25ex}

Чтобы распространить соответствие на всю сферу (1.1), на плоскости $Oxy$ введем условную
{\it бесконечно удаленную точку} $M^{}_{\infty},$ 
\vspace{0.15ex}
считая ее прообразом  северного полюса $N(0,0,1)$ при стереографическом проецировании.
\vspace{0.15ex}
Плоскость $Oxy,$ пополненную бесконечно удаленной точкой $M^{}_{\infty},$ 
\vspace{0.15ex}
образом которой при стереографическом проецировании на сфере (1.1) является северный полюс $N(0,0,1),$ 
\vspace{0.35ex}
назовем {\it расширенной плоскостью} $Oxy$ и обозначим $\overline{Oxy},$ 
\vspace{0.25ex}
т.е. $\overline{Oxy} = Oxy\sqcup M^{}_{\infty}\,.$
Расширенная плоскость $\overline{Oxy}$ состоит из двух частей: конечной --- собственно плоскость $Oxy,$ и
бесконечной --- точка $M^{}_{\infty}.$ 
\vspace{0.25ex}
Тогда любой последовательности $\{M_k^{}\}$ точек $M_k^{}$ плоскости $Oxy,$ 
\vspace{0.35ex}
удаляющихся в бесконечность (стремящихся к точке $M^{}_{\infty}),$ 
\vspace{0.35ex}
с помощью стереографической проекции однозначно сопоставляется последовательность 
\vspace{0.15ex}
$\{P_k^{}\}$ точек $P_k^{}$ сферы (1.1),  стремящихся к северному полюсу $N(0,0,1).$
\vspace{0.15ex}

Пусть точка $M$ расположена  на плоскости $Oxy$ и имеет
\vspace{0.15ex}
координаты $M(x,y).$ Тогда в
пространственной системе координат $O^{\star}x^{\star}y^{\star}z^{\star}$ эта же
точка имеет коодинаты $M(x,y,{}-1).$ Прямая $MN$
в системе координат $O^{\star}x^{\star}y^{\star}z^{\star}$ задается
системой уравнений
\\[2ex]
\mbox{}\hfill
$
\dfrac{x^{\star}}{x}=\dfrac{y^{\star}}{y}=\dfrac{z^{\star}-1}{{}-2}\,.
\hfill
$
\\[2ex]
При стереографическом проецировании 
\vspace{0.15ex}
с центром в северном полюсе $N(0,0,1)$ точке
$M(x,y,{}-1)$ соответствует  точка  $P(x^{\star},y^{\star},z^{\star}),$
\vspace{0.25ex}
являющаяся точкой пресечения прямой $MN$ и сферы (1.1). Поэтому
координаты $x^{\star}, y^{\star}, z^{\star}$ точки $P$ суть решения алгебраической  системы уравнений 
\\[2ex]
\mbox{}\hfill                 
$
\dfrac{x^{\star}}{x} =\dfrac{y^{\star}}{y}=\dfrac{1-z^{\star}}{2}\,,
\qquad 
x^{\star}{}^{\,^{\scriptstyle 2}}+
y^{\star}{}^{\,^{\scriptstyle 2}}+
z^{\star}{}^{\,^{\scriptstyle 2}}=1.
$
\hfill (1.2)
\\[2.25ex]
\indent
Будем считать, что точка 
\vspace{0.35ex}
$P(x^{\star},y^{\star},z^{\star})$  лежит  
на сфере (1.1) и не является северным полюсом  $N(0,0,1).$
\vspace{0.35ex}
Тогда ее аппликата $z^{\star}\in [{}-1;1).$
Разрешив систему уравнений (1.2) относительно 
$x^{\star},\,y^{\star},\,z^{\star}$ при ${}-1\leq z^{\star}<1,$
получим биективное отражение 
\\[2ex]
\mbox{}\hfill                                   
$
\psi_{_N} \colon (x,y)\to \bigl(x^{\star}(x,y), y^{\star}(x,y), z^{\star}(x,y)\bigr),
\hfill                                  
$
\\[-0.75ex]
\mbox{}\hfill {\rm (1.3)}
\\[0.75ex]
\mbox{}\hfill 
$
x^{\star}(x,y)=\dfrac{4x}{x^2+y^2+4}\,, \ \ \
y^{\star}(x,y)=\dfrac{4y}{x^2+y^2+4}\,, \ \ \
z^{\star}(x,y)=\dfrac{x^2+y^2-4}{x^2+y^2+4}
\quad
\forall (x,y)\in \R^2
\hfill 
$
\\[2.75ex]
плоскости $Oxy$  на сферу (1.1), 
\vspace{0.15ex}
проколотую в северном полюсе  $N(0,0,1).$
Координатные функции отображения (1.3) непрерывно дифференцируемы.
Якобианы 
\\[2.25ex]
\mbox{}\hfill                                  
$
\dfrac{{\sf D} (x^{\star},y^{\star})}
{{\sf D} (x,y)}=
{}-16\, \dfrac{x^2+y^2-4}{(x^2+y^2+4)^3}\,,
\qquad 
\dfrac{{\sf D} (x^{\star},z^{\star})}{{\sf D} (x,y)}=
64\, \dfrac{y}{(x^2+y^2+4)^3}\,,
\hfill                                  
$
\\[2.75ex]
\mbox{}\hfill                                  
$
\dfrac{{\sf D} (y^{\star},z^{\star})}{{\sf D} (x,y)}=
{}-64\, \dfrac{x}{(x^2+y^2+4)^3}
\quad 
\forall (x,y) \in \R^2
\hfill
$
\\[2.5ex]
одновременно не обращаются в нуль в любой точке плоскости $Oxy.$
Следовательно, с учетом леммы 1.1 имеет место 
\vspace{0.5ex}

{\bf Теорема 1.1.}
\vspace{0.15ex}
{\it
Стереографическое отображение {\rm (1.3)}  
плоскости $Oxy$ 
на сферу {\rm (1.1),} проколотую 
в северном полюсе  $N(0,0,1),$
является диффеоморфизмом.
}
\vspace{0.5ex}

Основное свойство диффеоморфизма (1.3) состоит в том, что [5, стб. 222 --- 223] 
угол между кривыми на плоскости и угол между стереографическими образами этих кривых на сфере равны. 
\\[3.75ex]
\mbox{}\hfill
{\unitlength=1mm
\begin{picture}(52,56)
\put(-7,0){\includegraphics[width=52mm,height=56mm]{r01-01.eps}}

\put(20.9,19.7){\makebox(0,0)[cc]{\scriptsize $O$}}
\put(17.7,23.7){\makebox(0,0)[cc]{\scriptsize $S$}}
\put(35.5,20.5){\makebox(0,0)[cc]{\scriptsize $x$}}
\put(40,20){\makebox(0,0)[cc]{\scriptsize $x$}}
\put(10,3.9){\makebox(0,0)[cc]{\scriptsize $y$}}
\put(22,25){\makebox(0,0)[cc]{\scriptsize $y$}}
\put(28,2.6){\makebox(0,0)[cc]{\scriptsize $M$}}

\put(17.2,37){\makebox(0,0)[cc]{\scriptsize $O^{\star}$}}
\put(24,32.5){\makebox(0,0)[cc]{\scriptsize $P$}}
\put(37.6,33.3){\makebox(0,0)[cc]{\scriptsize $x^{\star}$}}
\put(26,44.3){\makebox(0,0)[cc]{\scriptsize $y^{\star}$}}

\put(17.6,47){\makebox(0,0)[cc]{\scriptsize $N$}}
\put(21.5,54.5){\makebox(0,0)[cc]{\scriptsize $z^{\star}$}}

\put(20,-5){\makebox(0,0)[cc]{\rm Рис. 1.1}}
\end{picture}}
\qquad\qquad
{\unitlength=1mm
\begin{picture}(50,55)
\put(0,0){\includegraphics[width=49.79mm,height=55.84mm]{r02-01.eps}}

\put(25.9,19.7){\makebox(0,0)[cc]{\scriptsize $O$}}
\put(22.7,23.7){\makebox(0,0)[cc]{\scriptsize $S$}}
\put(40,20.5){\makebox(0,0)[cc]{\scriptsize $x$}}
\put(44.7,20){\makebox(0,0)[cc]{\scriptsize $x$}}
\put(15,3.9){\makebox(0,0)[cc]{\scriptsize $y$}}
\put(27,25){\makebox(0,0)[cc]{\scriptsize $y$}}
\put(34,2.6){\makebox(0,0)[cc]{\scriptsize $M$}}

\put(21.7,33.5){\makebox(0,0)[cc]{\scriptsize $O^{\star}$}}
\put(28.9,32.5){\makebox(0,0)[cc]{\scriptsize $P$}}
\put(42.6,33.3){\makebox(0,0)[cc]{\scriptsize $x^{\star}$}}
\put(27,46){\makebox(0,0)[cc]{\scriptsize $y^{\star}$}}

\put(31.9,42){\makebox(0,0)[cc]{\scriptsize $M^{\ast}$}}
\put(21.3,46.3){\makebox(0,0)[cc]{\scriptsize $O^{\ast}$}}
\put(22.7,50.3){\makebox(0,0)[cc]{\scriptsize $N$}}
\put(38.3,46){\makebox(0,0)[cc]{\scriptsize $u$}}
\put(31.5,49.8){\makebox(0,0)[cc]{\scriptsize $u$}}
\put(27.9,51.5){\makebox(0,0)[cc]{\scriptsize $v$}}
\put(20,41.9){\makebox(0,0)[cc]{\scriptsize $v$}}
\put(26.5,55.5){\makebox(0,0)[cc]{\scriptsize $z^{\star}$}}

\put(25,-5){\makebox(0,0)[cc]{\rm Рис. 2.1}}
\end{picture}}
\hfill\mbox{}
\\[-7ex]

\newpage

\mbox{}
\\[-1.75ex]
\centerline{
{\bf  2. Стереографический атлас сферы}
}
\\[1.5ex]
\indent
\indent
На плоскости, касающейся сферы (1.1) в северном полюсе $N(0,0,1),$
\vspace{0.15ex}
введем правую прямоугольную декартову систему координат
\vspace{0.15ex}
$O^{\ast}uv$ так, что ее начало $O^{\ast}(0,0)$ совпадает с северным полюсом
\vspace{0.35ex}
$N(0,0,1)$ сферы (1.1), ось $O^{\ast}u$ сонаправлена с осью $O^{\star}x^{\star},$ ось
\vspace{0.35ex}
$O^{\ast}v$ сонаправлена с осью $O^{\star}y^{\star}$ (рис. 2.1).
Масштаб в  системах координат $Oxy, \ O^{\star}x^{\star}y^{\star}z^{\star},\, O^{\ast}uv$ одинаков.
\vspace{0.5ex}

Если в плоскости $O^{\ast}uv$ точка $M^{\ast}$  имеет
\vspace{0.25ex}
координаты $M^{\ast}(u,v),$ то в
пространственной системе координат $O^{\star}x^{\star}y^{\star}z^{\star}$ эта же
\vspace{0.35ex}
точка имеет коодинаты $M^{\ast}(u,v,1).$ Прямая $M^{\ast}S$
в системе координат $O^{\star}x^{\star}y^{\star}z^{\star}$ задается
системой уравнений
\\[2ex]
\mbox{}\hfill
$
\dfrac{x^{\star}}{u}=\dfrac{y^{\star}}{v}=\dfrac{z^{\star}+1}{2}\,.
\hfill
$
\\[2ex]
При стереографическом проецировании с центром в южном полюсе 
\vspace{0.15ex}
$S(0,0,{}-1)$ сферы (1.1) образом точки 
$M^{\ast}(u,v,1)$ будет   точка  $P(x^{\star},y^{\star},z^{\star}),$
\vspace{0.35ex}
являющаяся точкой пресечения прямой $M^{\ast}S$ и сферы (1.1). Поэтому
\vspace{0.15ex}
координаты $x^{\star},\ y^{\star},\ z^{\star}$ точки $P$ суть решения алгебраической  системы уравнений 
\\[2ex]
\mbox{}\hfill                 
$
\dfrac{x^{\star}}{u} =\dfrac{y^{\star}}{v}=\dfrac{z^{\star}+1}{2}\,,
\qquad 
x^{\star}{}^{\,^{\scriptstyle 2}}+
y^{\star}{}^{\,^{\scriptstyle 2}}+
z^{\star}{}^{\,^{\scriptstyle 2}}=1.
$
\hfill (2.1)
\\[2.25ex]
\indent
Будем считать, что точка 
\vspace{0.35ex}
$P(x^{\star},y^{\star},z^{\star})$  лежит  
на сфере (1.1) и не является 
южным полюсом  $S(0,0,{}-1).$
Тогда ее аппликата $z^{\star}\in ({}-1;1].$
\vspace{0.35ex}
Разрешив систему уравнений (2.1) относительно 
$x^{\star},\ y^{\star},\ z^{\star}$ при ${}-1< z^{\star}\leq 1,$
получим биективное отражение 
\\[2ex]
\mbox{}\hfill                                   
$
\psi_{_S} \colon (u,v)\to\ \bigl(x^{\star}(u,v), y^{\star}(u,v), z^{\star}(u,v)\bigr),
\hfill                                  
$
\\[-0.75ex]
\mbox{}\hfill {\rm (2.2)}
\\[0.75ex]
\mbox{}\hfill 
$
x^{\star}(u,v)=\dfrac{4u}{u^2+v^2+4}\,, \ \ \
y^{\star}(u,v)=\dfrac{4v}{u^2+v^2+4}\,, \ \ \
z^{\star}(u,v)={}-\dfrac{u^2+v^2-4}{u^2+v^2+4}
\quad
\forall (u,v)\in \R^2
\hfill 
$
\\[2.75ex]
плоскости $O^{\ast}uv $  на сферу (1.1), проколотую 
\vspace{0.15ex}
в южном полюсе  $S(0,0,{}-1).$
Координатные функции отображения (2.2) непрерывно дифференцируемы.
Якобианы 
\\[2.25ex]
\mbox{}\hfill                                  
$
\dfrac{{\sf D} (x^{\star},y^{\star})}
{{\sf D} (u,v)}=
{}-16\, \dfrac{u^2+v^2-4}{(u^2+v^2+4)^3}\,,
\qquad 
\dfrac{{\sf D} (x^{\star},z^{\star})}{{\sf D} (u,v)}=
{}-64\, \dfrac{v}{(u^2+v^2+4)^3}\,,
\hfill                                  
$
\\[2.75ex]
\mbox{}\hfill                                  
$
\dfrac{{\sf D} (y^{\star},z^{\star})}{{\sf D} (u,v)}=
64\, \dfrac{u}{(u^2+v^2+4)^3}
\quad 
\forall (u,v) \in \R^2
\hfill
$
\\[2.5ex]
одновременно не обращаются в нуль в любой точке плоскости $O^{\ast}uv.$
\vspace{0.15ex}
Следовательно, стереографическое отображение (2.2)   плоскости $O^{\ast}uv$ на сферу (1.1), проколотую 
\vspace{0.15ex}
в южном полюсе  $S(0,0,{}-1),$
является диффеоморфизмом.
\vspace{0.35ex}

Плоскость $O^{\ast}uv,$ 
\vspace{0.35ex}
пополненная бесконечно удаленной точкой $M^{\ast}_{\infty},$ образом которой 
при стереографическом проецировании плоскости $O^{\ast}uv$ на сферу (1.1) 
\vspace{0.25ex}
является южный полюс $S(0,0,{}-1),$ есть  расширенная плоскость $O^{\ast}uv,$ 
\vspace{0.5ex}
 которую обозначим $\overline{O^{\ast}uv},$ т.е. $\overline{O^{\ast}uv} = O^{\ast}uv\sqcup M^{\ast}_{\infty}\,.$
\vspace{0.5ex}

Покроем сферу (1.1) двумя сферами $U_1^{}$ и $U_2^{},$ 
\vspace{0.35ex}
проколотыми соответственно в северном  $N(0,0,1)$ и южном  $S(0,0,{}-1)$ полюсах:
\\[2ex]
\mbox{}\hfill
$
\displaystyle
U_1^{}=\Bigl\{
(x^{\star}, y^{\star}, z^{\star})\colon\,   
{x^{\star}}^{\,^{\scriptstyle 2}}+{y^{\star}}^{\,^{\scriptstyle 2}}+{z^{\star}}^{\,^{\scriptstyle 2}}=1, \ \,
{}-1\leq z^{\star}<1\Bigr\}
\hfill
$
\\[1ex]
и 
\\[1ex]
\mbox{}\hfill
$
\displaystyle
U_2^{}=\Bigl\{
(x^{\star}, y^{\star}, z^{\star})\colon\,   
{x^{\star}}^{\,^{\scriptstyle 2}}+{y^{\star}}^{\,^{\scriptstyle 2}}+{z^{\star}}^{\,^{\scriptstyle 2}}=1, \ \,
{}-1< z^{\star}\leq1\Bigr\}.
\hfill
$
\\[2.25ex]
\indent
Введем диффеоморфные отображения 
\\[2.25ex]
\mbox{}\hfill                                                              
$
\displaystyle
\varphi_1^{}\colon (x^{\star}, y^{\star}, z^{\star})\to \
\bigl(
x(x^{\star}, y^{\star}, z^{\star}),\, y(x^{\star}, y^{\star}, z^{\star})
\bigr),
\hfill
$
\\[0ex]
\mbox{}\hfill (2.3)
\\[0ex]
\mbox{}\hfill
$
\displaystyle
x(x^{\star}, y^{\star}, z^{\star})=\dfrac{2x^{\star}}{1-z^{\star}}\,, 
\qquad 
y(x^{\star}, y^{\star}, z^{\star})=\dfrac{2y^{\star}}{1-z^{\star}}
\quad
\forall (x^{\star}, y^{\star}, z^{\star})\in U_1^{},
\hfill
$
\\[3.5ex]
\mbox{}\hfill                                                              
$
\displaystyle
\varphi_2^{}\colon (x^{\star}, y^{\star}, z^{\star})\to\ 
\bigl(
u(x^{\star}, y^{\star}, z^{\star}),\, v(x^{\star}, y^{\star}, z^{\star})
\bigr),
\hfill
$
\\[0ex]
\mbox{}\hfill (2.4)
\\[0ex]
\mbox{}\hfill
$
\displaystyle
u(x^{\star}, y^{\star}, z^{\star})=\dfrac{2x^{\star}}{1+z^{\star}}\,, 
\qquad 
v(x^{\star}, y^{\star}, z^{\star})=\dfrac{2y^{\star}}{1+z^{\star}}
\quad
\forall (x^{\star}, y^{\star}, z^{\star})\in U_2^{}.
\hfill
$
\\[2.5ex]
\indent
Отображение (2.3) является обратным к отображению (1.3), 
\vspace{0.15ex}
а значит, отображение (2.3) есть диффеоморфное отображение сферы $U_1^{},$ 
\vspace{0.25ex}
проколотой в северном полюсе $N(0,0,1),$ на плоскость $Oxy$
\vspace{0.25ex}
--- стереографическое отображение сферы (1.1) на плоскость $Oxy$ из центра 
\vspace{0.25ex}
в северном полюсе $N(0,0,1)$
[5, стб. 222 --- 223; 6, c. 37].

Отображение (2.4) является обратным к отображению (2.2), 
\vspace{0.15ex}
а значит, отображение (2.4) есть диффеоморфное отображение сферы $U_2^{}\;\!\!,$ 
\vspace{0.25ex}
проколотой в южном полюсе $S(0,0,{}-1),$ на плоскость $O^{\ast}uv$
\vspace{0.25ex}
--- стереографическое отображение сферы (1.1) на плоскость $O^{*}uv$ из центра в южном полюсе $S(0,0,{}-1).$
\vspace{0.35ex}

Таким образом, 
\vspace{0.25ex}
построены две карты  $(U_1^{}, \varphi_1^{}) $ и $(U_2^{}, \varphi_2^{}) $ сферы (1.1). 
Множество карт 
$(U_1^{}, \varphi_1^{}) $ и $(U_2^{}, \varphi_2^{}) $
образуют {\it стереографический атлас} сферы  (1.1) [6, c. 103].
\vspace{0.5ex}

Установим связь между локальными системами координат $Oxy$ и $O^{*}uv$
стереографического атласа сферы  (1.1). 
\vspace{0.15ex}
Для этого используем стереографические отображения (1.3) и (2.2), 
\vspace{0.15ex}
которые являются взаимообратными соответственно с отображениями 
(2.3) и (2.4) стереографического атласа сферы  (1.1):
\\[1.5ex]
\mbox{}\hfill                                                      
$
\varphi_1^{}=\psi_{_N}^{{}-1},
\quad 
\psi_{_N}=\varphi_1^{{}-1}
$
\ \ и \ \ 
$
\varphi_2^{}=\psi_{_S}^{{}-1},
\quad 
\psi_{_S}=\varphi_2^{{}-1}.
\hfill
$
\\[2.5ex]
\indent
Выполним последовательно два диффеоморфных отображения $\psi_{_N}$ и $\psi_{_S}^{{}-1},$ которые составят 
диффеоморфное отображение
\\[2ex]
\mbox{}\hfill                                                     
$
\varphi_{21}^{}=\psi_{_S}^{{}-1}\circ\psi_{_N}= \varphi_2^{}\circ\psi_{_N}
$
\hfill (2.5)
\\[2.25ex]
плоскости $Oxy,$ 
\vspace{0.15ex}
проколотой в начале координат $O(0,0),$ на плоскость $O^{\ast}uv,$ проколотую в начале координат $O^{\ast}(0,0)$ (рис. 2.1). 
\vspace{0.25ex}

Учитывая аналитические задания (1.3) и (2.4) отображений  $\psi_{_N}$ и $\varphi_{2}^{},$
\vspace{0.15ex}
получаем аналитическое задание диффеоморфного отображения (2.5) в виде 
\\[2ex]
\mbox{}\hfill                           
$
\varphi_{21}^{} \colon (x,y) \to\
\Bigl(\;\!\dfrac{4x}{x^{2}+y^{2}}\,,\
\dfrac{4y}{x^{2}+y^{2}}\Bigr)
\quad
\forall (x,y)\in\R^2\backslash\{(0,0)\}.
$
\hfill (2.6)
\\[2ex]
\indent
Коодинатные функции 
\\[2ex]
\mbox{}\hfill                           
$
u \colon (x,y) \to\ \dfrac{4x}{x^{2}+y^{2}}\,,
\quad\
v \colon (x,y) \to\ \dfrac{4y}{x^{2}+y^{2}} 
\quad
\forall (x,y)\in\R^2\backslash\{(0,0)\}
$
\hfill (2.7)
\\[2ex]
диффеоморфного  отображения (2.6) 
\vspace{0.25ex}
есть функции перехода [6, c. 99] от локальных координат $(u,v)$ к 
\vspace{0.25ex}
локальным координатам $(x,y)$ стереографического атласа сферы  (1.1).

У диффеоморфного  отображения 
\\[2ex]
\mbox{}\hfill                           
$
\varphi_{12}^{} \colon (u,v) \to\
\Bigl(\;\!\dfrac{4u}{u^{2}+v^{2}}\,,\
\dfrac{4v}{u^{2}+v^{2}}\Bigr)
\quad
\forall (u,v)\in\R^2\backslash\{(0,0)\}
$
\hfill (2.8)
\\[2ex]
коодинатные функции 

\newpage

\mbox{}
\\[-2ex]
\mbox{}\hfill                           
$
x \colon (u,v) \to\ \dfrac{4u}{u^{2}+v^{2}}\,,
\quad \ 
y \colon (u,v) \to\ \dfrac{4v}{u^{2}+v^{2}} 
\quad
\forall (u,v)\in\R^2\backslash\{(0,0)\}
$
\hfill (2.9)
\\[2ex]
есть функции перехода от локальных координат $(x,y)$ к 
\vspace{0.15ex}
локальным координатам $(u,v)$ стереографического атласа сферы  (1.1).

Диффеоморфное  отображение (2.8), полученное как обратное отображение к отображению (2.6), 
есть аналитическое задание отображения  
\\[2ex]
\mbox{}\hfill
$
\varphi_{12}^{}=\psi_{_N}^{{}-1}\circ\psi_{_S}= \varphi_1^{}\circ\psi_{_S}
\hfill
$
\\[2ex]
плоскости $O^{\ast}uv,$  
\vspace{0.15ex}
проколотой в начале координат $O^{\ast}(0,0),$  на плоскость $Oxy,$ проколотую в начале координат $O(0,0).$
\vspace{0.5ex}

{\bf Теорема 2.1.}
\vspace{0.15ex}
{\it
Диффеоморфное отображение $(2.6)$ и тождественное отображение 
плоскости $Oxy$ на себя образуют группу второго порядка.
}
\vspace{0.35ex}

{\sl Действительно}, преобразование (2.6) взаимообратимо:
\\[2ex]
\mbox{}\hfill
$
(x,y)\ \stackrel{\stackrel{\scriptstyle \varphi_{21}^{}}{\mbox{}}}
{\longrightarrow}\ 
\Bigl(\;\!\dfrac{4x}{x^2+y^2}\,,\
\dfrac{4y}{x^2+y^2}\Bigr)\ 
\stackrel{\stackrel{\scriptstyle \varphi_{21}^{}}{\mbox{}}}
{\longrightarrow}\ 
\left(\dfrac{4\cdot\dfrac{4x}{x^2+y^2}}
{\Bigl(\dfrac{4x}{x^2+y^2}\Bigr)^2
+\Bigl(\dfrac{4y}{x^2+y^2}\Bigr)^2},
\right.
\hfill
$
\\[1ex]
\mbox{}\hfill
$
\left.
\dfrac{4\cdot\dfrac{4y}{x^2+y^2}}
{\Bigl(\dfrac{4x}{x^2+y^2}\Bigr)^2
+\Bigl(\dfrac{4y}{x^2+y^2}\Bigr)^2}\right)=(x,y)\quad
\forall (x,y)\in\R^2\backslash \{(0,0)\},
\hfill
$
\\[3ex]
а также, 
$
\varphi_{21}^{}\circ I=I\circ \varphi_{21}^{}=\varphi_{21}^{}$ и $I\circ I=I,
$
где $I\colon (x,y)\to (x,y)\; \; \forall (x,y)\in\R^2.\, \k$
\\[5.75ex]
\centerline{
{\bf\large 
\S\;\!2. Стереографически сопряженная дифференциальная система}
}
\\[2.25ex]
\centerline{
{\bf  3. Преобразование Бендиксона
}
}
\\[1.5ex]
\indent
{\it Преобразованием Бендиксона} фазовой плоскости $Oxy$ дифференциальной системы (D) назовем преобразование 
\\[1ex]
\mbox{}\hfill                           
$
x = \dfrac{4u}{u^{2}+v^{2}}\,,
\qquad 
y =\dfrac{4v}{u^{2}+v^{2}} \,,
$
\hfill (3.1)
\\[2.5ex]
построенное на основании функций перехода (2.9) от локальных координат $x,\,y$ 
\vspace{0.15ex}
к локальным координатам $u,\,v$ стереографического атласа сферы (1.1).
\vspace{0.35ex}

Заменой Бендиксона (3.1) дифференциальную систему (D)
\vspace{0.15ex}
приводим к дифференциальной системе [7, c. 239]
\\[2ex]
\mbox{}\hfill                           
$
\dfrac{du}{dt} ={}-\dfrac{\,u^{2}
-v^{2}}{4}\,  X\Bigl(\dfrac{4u}{u^{2}+v^{2}}\,,\
\dfrac{4v}{u^{2}+v^{2}}\Bigr)-
\dfrac{uv}{2}\,   Y\Bigl(\dfrac{4u}{u^{2}+v^{2}}\,,\
\dfrac{4v}{u^{2}+v^{2}}\Bigr)\equiv
\ \stackrel{\ast}{U}(u,v),
\hfill
$
\\[.5ex]
\mbox{}\hfill (3.2)
\\[.5ex]
\mbox{}\hfill
$
\dfrac{dv}{dt}=
{}-\dfrac{u v}{2}\,   X\Bigl(\dfrac{4u}{u^{2}+v^{2}}\,,\
\dfrac{4v}{u^{2}+v^{2}}\Bigr)+
\dfrac{u^{2}-v^{2}}{4}\,
Y\Bigl(\dfrac{4u}{u^{2}+v^{2}}\,,\
\dfrac{4v}{u^{2}+v^{2}}\Bigr)\equiv \ \stackrel{\ast}{V}(u,v).
\hfill
$
\\[2.75ex]
\indent
Поскольку  $X$ и $Y$ --- полиномы, то систему (3.2) можно записать в виде
\\[2.1ex]
\mbox{}\hfill                         
$
\dfrac{d u}{dt}=\dfrac{U (u,v )}{(u^2+v^2)^{m}}\,,
\qquad
\dfrac{dv}{dt}=\dfrac{V(u,v)}{(u^2+v^2)^{m}}\,,
\hfill
$
\\[2.5ex]
где $U$ и $V$ --- полиномы, 
\vspace{0.15ex}
не делящиеся одновременно на
$u^2+v^2,$ а число $m$ --- целое неотрицательное.

\newpage

Автономную полиномиальную дифференциальную
систему
\\[2.1ex]
\mbox{}\hfill                           
$
\dfrac{d u}{d\tau}=U(u,v),
\qquad
\dfrac{d v}{d\tau}=V(u,v),
$
\hfill (3.3)
\\[2.5ex]
где $(u^2+v^2)^{m} d\tau=dt,$ 
\vspace{0.35ex}
у которой правые части
$U$ и $V$ суть взаимно простые полиномы, назовем {\it стереографически сопряженной}
к дифференциальной системе (D). 
\vspace{0.35ex}

Согласно теореме 2.1 система (D) 
\vspace{0.15ex}
является стереографически сопряженной к системе (3.3), и с помощью преобразования Бендиксона 
\\[2ex]
\mbox{}\hfill                         
$
u = \dfrac{4x}{x^{2}+y^{2}}\,,
\qquad 
v =\dfrac{4y}{x^{2}+y^{2}}
\hfill
$
\\[2ex]
система (3.3) приводится к системе (D).
\vspace{0.15ex}
Дифференциальные системы  (D) и (3.3) являются стереографически взаимосопряженными.
\vspace{0.35ex}

Фазовую плоскость $Oxy$ 
\vspace{0.35ex}
(расширенную фазовую плоскость $\overline{Oxy}$) системы (D) и фазовую плоскость 
\vspace{0.15ex}
$O^{\ast}uv$ (расширенную фазовую плоскость $\overline{O^{\ast}uv}$) системы (3.3) назовем 
стереографически сопряженными.
\\[4.25ex]
\centerline{
{\bf  4.
Вид стереографически сопряженной дифференциальной системы}
}
\\[1.5ex]
\indent
Вид стереографически сопряженной 
\vspace{0.15ex}
к дифференциальной системе (D)
дифференциальной системы (3.3) зависит от того, делится ли на $x^2+y^2$ или нет полином
\\[2ex]
\mbox{}\hfill
$
\displaystyle
W_n^{}\colon (x,y) \to\  x\;\!Y_n^{}(x,y) - yX_n^{}(x,y) \quad \forall (x,y)\in \R^2.
\hfill
$
\\[1.5ex]
\indent
Если 
\\[1ex]
\mbox{}\hfill
$
W_n^{}(x,y)\not\equiv (x^2+y^2)P(x,y)$ на $\R^2,
\hfill
$ 
\\[2ex]
где $P$ --- некоторый полином, то дифференциальная система (3.3) имеет вид 
\\[2.5ex]
$
\displaystyle
\dfrac{du}{d\tau} =
\dfrac{v^2 - u^2}{4}\ \sum\limits_{j = 0}^{n}\,
(u^2 + v^2)^{n - j}\  X_j^{}(4u,4v)  -
 \dfrac{uv}{2}\  \sum\limits_{j = 0}^{n}\,
(u^2 + v^2)^{n - j}\ Y_j^{}(4u,4v) \equiv U_{_0} (u,v),
\hfill
$
\\[0.25ex]
\mbox{}\hfill                   (4.1)
\\[0.25ex]
$
\displaystyle
\dfrac{dv}{d\tau} =
{} - \dfrac{uv}{2}\  \sum\limits_{j = 0}^{n}\,
(u^2+ v^2)^{n - j}\  X_j^{}(4u,4v)  +
\dfrac{u^2 - v^2}{4} \ \sum\limits_{j = 0}^{n}\,
(u^2 + v^2)^{n - j}\ Y_j^{}(4u,4v)\equiv V_{_0} (u,v),
\hfill
$
\\[3ex]
где  $(u^2 + v^2)^n\, d\tau = dt.$
\vspace{0.5ex}

Пусть 
\\[1ex]
\mbox{}\hfill
$
W_n^{}(x,y)=(x^2+y^2)P(x,y)
\quad 
\forall (x,y)\in \R^2,
\hfill
$ 
\\[2ex]
где $P$ --- некоторый полином 
\vspace{0.35ex}
(не исключается случай $P(x,y)=0\;\; \forall (x,y)\in \R^2$\!). Тогда, если имеют место тождества
\\[2ex]
\mbox{}\hfill  
$
\displaystyle
{}-2y\bigl(xY_{n-r+1}^{}(x,y) - yX_{n-r+1}^{}(x,y)\bigr)  -
 (x^2 + y^2) X_{n-r+1}^{}(x,y) =
(x^2 + y^2)^{k - r+ 1}\ K_r^{}(x,y),
\hfill
$
\\[2.5ex]
\mbox{}\hfill
$
\displaystyle
2x\bigl(xY_{n-r+1}^{}(x,y) - yX_{n-r+1}^{}(x,y)\bigr) -
 (x^2 + y^2) Y_{n-r+1}^{}(x,y) 
=
(x^2 + y^2)^{k - r + 1}\ Q_r^{}(x,y)
$
\hfill (4.2)
\\[2.5ex]
\mbox{}\hfill
$
\displaystyle
\forall (x,y)\in\R^2, \ \ r=1,\ldots, k,
\hfill
$
\\[2ex]
где натуральное число $k$ такое, что $2k\leq n + 2,$ а
\vspace{0.5ex}
$K_r^{}$ и $Q_r^{},\, r=1,\ldots, k,$ суть  некоторые полиномы,
то дифференциальная система (3.3) имеет вид 

\newpage

\mbox{}
\\[-2ex]
\mbox{}\qquad
$
\displaystyle
\dfrac{du}{d\theta_k^{}} =
\dfrac{v^2 - u^2}{4} \  \sum\limits_{j = 0}^{n-k}\,
(u^2 +  v^2)^{n - j - k} \, X_j^{}(4u,4v)
 -
\dfrac{uv}{2}\  \sum\limits_{j = 0}^{n - k}\,
(u^2 + v^2)^{n - j - k}\ Y_j^{}(4u,4v)  \,+
\hfill
$
\\[1.75ex]
\mbox{}\qquad\qquad\quad
$
\displaystyle
+\,
\sum \limits_{r = 1}^{k}\, 4^{2k-2r-1}\ K_{r}^{}(4u,4v)\equiv U_k^{} (u,v),
\hfill
$
\\[0.5ex]
\mbox{}\hfill (4.3)
\\[0.5ex]
\mbox{}\qquad
$
\displaystyle
\dfrac{dv}{d\theta_k^{}} = {}-\dfrac{uv}{2} \
\sum\limits_{j = 0}^{n - k}\,
(u^2 + v^2)^{n - j - k}\ X_j^{}(4u,4v)  
+
 \dfrac{u^2 - v^2}{4}\ 
\sum\limits_{j = 0}^{n - k}\,
(u^2 + v^2)^{n - j - k}\ Y_j^{}(4u,4v)  \,+
\hfill
$
\\[1.75ex]
\mbox{}\qquad\qquad\quad
$
\displaystyle
+\,
\sum\limits_{r=1}^{k} \, 4^{2k - 2r - 1}\ Q_{r}^{}(4u,4v)\equiv V_k^{} (u,v),
\hfill
$
\\[3ex]
где $(u^2 + v^2)^{n - k}\, d\theta_k^{}  = dt.$
\vspace{0.75ex}

{\bf Пример 4.1.} Рассмотрим автономную дифференциальную систему
\\[2ex]
\mbox{}\hfill
$
\displaystyle
\dfrac{dx}{dt}  = a_{_0}\equiv X(x,y),
\qquad 
\dfrac{dy}{dt}  = b_{_0}\equiv Y(x,y),
\qquad 
| a_{_0}|+| b_{_0}| \ne 0,
$
\hfill (4.4)
\\[2.5ex]
у которой 
\vspace{0.5ex}
$W_{_0}(x,y)=b_{_0} x-a_{_0} y\not\equiv (x^2+y^2)P(x,y)$ на $\R^2,$ где $P$ --- некоторый полином.

Стереографически сопряженной к дифференциальной  системе (4.4) является система 
\\[2.25ex]
\mbox{}\hfill
$
\displaystyle
\dfrac{du}{dt} =
{}-\dfrac{a_{_0}}{4}\; u^2-\dfrac{b_{_0}}{2}\;  uv+\dfrac{a_{_0}}{4}\; v^2,
 \qquad
\dfrac{dv}{dt} =
\dfrac{b_{_0}}{4}\; u^2-\dfrac{a_{_0}}{2}\;  uv-\dfrac{b_{_0}}{4}\; v^2,
\qquad
| a_{_0}|+| b_{_0}| \ne 0.
$
\hfill (4.5)
\\[3ex]
\indent
{\bf Пример 4.2.} Рассмотрим  стационарную линейную систему 
\\[2ex]
\mbox{}\hfill
$
\displaystyle
\dfrac{dx}{dt}  =  a_{_0}  +  a_{1}^{}x  +  a_{2}^{}y
\equiv X(x,y),
\qquad
\dfrac{dy}{dt} =  b_{_0}  +  b_{1}^{}x  +  b_{2}^{}y
\equiv Y(x,y),
$
\hfill (4.6)
\\[2ex]
где $|a_{1}^{}| + |a_{2}^{}| +|b_1^{}| + |b_2^{}|\ne 0.$
\vspace{0.75ex}

Стереографически сопряженной к  стационарной линейной  системе (4.6) пpи
\\[1.5ex]
\mbox{}\hfill
$
|a_1^{} - b_2^{}|  +| a_2^{} + b_1^{}|  \ne  0,
\quad
|a_{1}^{}| + |a_{2}^{}| +|b_1^{}| + |b_2^{}|\ne 0
\hfill
$
\\[1.5ex]
является дифференциальная система 
\\[2ex]
$
\displaystyle
\dfrac{du}{d\tau}  =
{} - a_1^{}u^3 -  (a_2^{} + 2b_1^{})u^2\;\!v
+  (a_1^{} - 2b_2^{})uv^2 +  a_2^{}v^3 \ -
\hfill 
$
\\[2ex]
\mbox{}\hspace{13.35em} 
$
-\ \dfrac{1}{4}\; a_{_0}u^4  -  2b_{_0}u^3\;\!v  -
2b_{_0}uv^3  +  \dfrac{1}{4}\; a_{_0}v^4\equiv U_{_0} (u,v),
$
\hfill (4.7)
\\[1.5ex]
$
\displaystyle
\dfrac{dv}{d\tau}  =
b_1^{}u^3  -  (2a_1^{} - b_2^{})u^2\;\!v  - (2a_2^{} + b_1^{})uv^2  -  b_2^{}v^3   +
b_{_0}u^4  -  2a_{_0}u^3\;\!v  - 2a_{_0}uv^3 +  b_{_0}v^4\equiv V_{_0} (u,v),
\hfill
$
\\[2.25ex]
где  $(u^2 + v^2)\,d\tau  =  dt.$
\vspace{0.5ex}

Стереографически сопряженной к  стационарной линейной системе (4.6) пpи
\\[1.5ex]
\mbox{}\hfill
$
b_2^{} = a_1^{},\quad 
b_1^{} ={} - a_2^{},
\quad
|a_{1}^{}| + |a_{2}^{}|\ne 0
\hfill
$
\\[1.5ex]
является дифференциальная система 
\\[2ex]
\mbox{}\hfill
$
\displaystyle
\dfrac{du}{dt}  =
{} - a_1^{}u  +  a_2^{}v  -  \dfrac{a_{_0}}{4}\; u^2  -
\dfrac{b_{_0}}{2}\; uv   +  \dfrac{a_{_0}}{4}\; v^2\equiv U_{1}^{} (u,v),
\hfill
$
\\[0.25ex]
\mbox{}\hfill                   (4.8)
\\[0.25ex]
\mbox{}\hfill
$
\displaystyle
\dfrac{dv}{dt}  =
{} -  a_2^{}u  -  a_1^{}v  +   \dfrac{b_{_0}}{4}\; u^2  -
\dfrac{a_{_0}}{2}\; uv   -  \dfrac{b_{_0}}{4}\; v^2\equiv V_{1}^{} (u,v).
\hfill 
$
\\[3.5ex]
\indent
{\bf Пример 4.3.} 
Рассмотрим автономную квадратичную систему 
\\[2ex]
\mbox{}\hfill
$
\displaystyle
\dfrac{dx}{dt}  =  a_{_0} +  a_1^{}x  +  a_2^{}y  +  a_3^{}x^2
 +  a_4^{}xy  +  a_5^{}y^2\equiv X(x,y),
\hfill
$
\\
\mbox{}\hfill                   (4.9)
\\
\mbox{}\hfill
$
\displaystyle
\dfrac{dy}{dt}  =  b_{_0} +  b_1^{}x  +  b_2^{}y  +
b_3^{}x^2  +  b_4^{}xy  +  b_5^{}y^2\equiv Y(x,y),
\;\ \hfill 
$
\\[2.75ex]
где  $|a_{3}^{}| + |a_{4}^{}| + |a_{5}^{}|+|b_3^{}| + |b_4^{}| + |b_{5}^{}|\ne 0.$
\vspace{1ex}

Пpи
\vspace{0.5ex}
$|a_5^{} - a_3^{} + b_4^{}|  +
|a_4^{} + b_3^{} - b_5^{}| \ne  0,\ 
|a_{3}^{}| + |a_{4}^{}| + 
|a_{5}^{}|+|b_3^{}| + |b_4^{}| + |b_{5}^{}|\ne 0$
стереографически сопряженной к   системе (4.9) 
является дифференциальная система 
\\[2.25ex]
\mbox{}\qquad
$
\dfrac{du}{d\tau}  =
{}- 4a_3^{}u^4 -  4(a_4^{} + 2b_3^{})u^3\;\!v  + 4(a_3^{} - 2b_4^{} - a_5^{})u^2\;\!v^2   +
 4(a_4^{} - 2b_5^{})uv^3 +4a_5^{}v^4
\, -
\hfill
$
\\[2ex]
\mbox{}\qquad\qquad\quad
$
-\, 
a_1^{}u^5  +  (2b_1^{} - a_2^{})u^4\;\!v  -
2b_2^{}u^3\;\!v^2  -
2b_1^{}u^2\;\!v^3  +  (a_1^{} - 2b_2^{})uv^4  +
 a_2^{}v^5
\,- 
\hfill
$
\\[2ex]
\mbox{}\qquad\qquad\quad
$
-\,  \dfrac{a_{_0}}{4}\; u^6  -
 \dfrac{b_{_0}}{2}\; u^5\;\!v  -
\dfrac{a_{_0}}{4}\; u^4\;\!v^2  -
b_{_0}u^3\;\!v^3  +  \dfrac{a_{_0}}{4}\; u^2\;\!v^4  -
\dfrac{b_{_0}}{2}\; uv^5  +  \dfrac{a_{_0}}{4}\; v^6 \equiv U_{_0} (u,v),
\hfill
$
\\[1.5ex]
\mbox{}\hfill (4.10)
\\[0.25ex]
\mbox{}\qquad
$
\dfrac{dv}{d\tau}  =
4b_3^{}u^4  +  4(b_4^{} - 2a_3^{})u^3\;\!v  -
4(b_3^{}+ 2a_4^{} + b_5^{})u^2\;\!v^2   -
4(b_4^{} + 2a_5^{})uv^3   - 4b_5^{}v^4 
\,+ 
\hfill
$
\\[2.25ex]
\mbox{}\qquad\qquad\quad
$
+\,    b_1^{}u^5
 +  (b_2^{} - 2a_1^{})u^4\;\!v  -
 2a_2^{}u^3\;\!v^2  -2a_1^{}u^2\;\!v^3  -   (b_1^{} + 2a_2^{})uv^4   
- b_2^{}v^5 
\,+
 \hfill
$
\\[2.5ex]
\mbox{}\qquad\qquad\quad
$
 +\,  \dfrac{b_{_0}}{4}\; u^6   -
\dfrac{a_{_0}}{2}\; u^5\;\!v  +
\dfrac{b_{_0}}{4}\; u^4\;\!v^2  -
a_{_0}u^3\;\!v^3  -  \dfrac{b_{_0}}{4}\;u^2\;\!v^4  -
\dfrac{a_{_0}}{2}\; uv^5  -  \dfrac{b_{_0}}{4}\; v^6\equiv V_{_0} (u,v),
\hfill
$
\\[2.25ex]
где $(u^2 + v^2)^2\,d\tau  =  dt.$
\vspace{0.75ex}

В случае
\vspace{1ex}
$a_4^{}  =  b_5^{}  -  b_3^{}, \  a_5^{}  =   a_3^{} -  b_4^{}, \
|b_2^{} - a_1^{}| +|b_1^{} + a_2^{}|+
|2a_3^{}- b_4^{}| +  |b_3^{} + b_5^{}| \ne 0,
$ 
$
|a_{3}^{}| + | b_5^{}  -  b_3^{}| + |a_3^{} -  b_4^{}|+|b_3^{}| + |b_4^{}| + |b_{5}^{}|\ne 0$
\vspace{0.75ex}
стереографически сопряженной к  системе (4.9) 
является дифференциальная  система 
\\[2ex]
\mbox{}\qquad\quad
$
\dfrac{du}{d\theta}  =
{}- 4a_3^{}u^2  -  4(b_3^{} + b_5^{})uv  +
4(a_3^{} - b_4^{})v^2  -  a_1^{}u^3  -
(a_2^{} + 2b_1^{})u^2\;\!v
\, +
\hfill
$
\\[1.75ex]
\mbox{}\qquad\qquad\qquad
$
+ \,
(a_1^{} - 2b_2^{})uv^2  + 2a_2^{}v^3   -
\dfrac{a_{_0}}{4}\; u^4  -  \dfrac{b_{_0}}{2}\; u^3\;\!v
-  \dfrac{b_{_0}}{2}\; uv^3  +  \dfrac{a_{_0} }{4}\; v^4\equiv U_{1}^{} (u,v),
\hfill
$
\\[1ex]
\mbox{}\hfill (4.11)
\\[-0.15ex]
\mbox{}\qquad\quad
$
\dfrac{dv}{d\theta}  =
{} -  4b_3^{}u^2  -
4(b_4^{} - 2a_3^{})uv  + 4b_5^{}v^2  +  b_1^{}u^3  +  (b_2^{} - 2a_1^{})u^2\;\!v 
\, -
\hfill
$
\\[1.75ex]
\mbox{}\qquad\qquad\qquad
$
- \, 
(b_1^{} + 2a_2^{})uv^2  -   b_2^{}v^3  +
\dfrac{b_{_0}}{4}\; u^4  -  \dfrac{a_{_0}}{2}\; u^3\;\!v  -
\dfrac{a_{_0}}{2}\; uv^3  -  \dfrac{b_{_0}}{4}\; v^4\equiv V_{1}^{} (u,v),
\hfill
$
\\[2.25ex]
где
\vspace{0.5ex}
$(u^2 + v^2)\,d\theta =  dt.$ 

Пpи
\vspace{0.75ex}
$a_4^{} = 2b_5^{}, \ a_5^{} = {}-a_3^{}, \ b_3^{} = {}- b_5^{},
\ b_4^{} = 2a_3^{}, \
 b_2^{} = a_1^{}, \ b_1^{} =  {}- a_1^{},\  
|a_{3}^{}| + | b_5^{}|\ne 0$
сте\-ре\-ог\-ра\-фи\-чес\-ки сопряженной к  системе (4.9) 
является дифференциальная система 
\\[2.5ex]
\mbox{}\hfill
$
\dfrac{du}{dt}  = {} - 4a_3^{}  - a_1^{}u  +  a_2^{}v  -  \dfrac{a_{_0}}{4}\; u^2  -
\dfrac{b_{_0}}{2}\; uv  +  \dfrac{a_{_0}}{4}\; v^2\equiv U_{2}^{} (u,v),
\hfill
$
\\[0.5ex]
\mbox{}\hfill (4.12)
\\[0.5ex]
\mbox{}\hfill
$
\dfrac{dv}{dt}  =
4b_5^{}  - a_2^{}u  -  a_1^{}v  +  \dfrac{b_{_0}}{4}\;u^2  -
\dfrac{a_{_0}}{2}\; uv  -  \dfrac{b_{_0}}{4}\; v^2\equiv V_{2}^{} (u,v).
\quad\;\, 
\hfill
$
\\

\newpage

\mbox{}
\\[-1.5ex]
\centerline{
{\bf  5. Стереографический атлас траекторий дифференциальной системы
}
}
\\[1.75ex]
\indent
Плоскостью 
\\[1ex]
\mbox{}\hfill
$
\{(x^{\star},y^{\star},z^{\star})\colon z^{\star}=1-\varepsilon_1^{}\},
\ \ \
0<\varepsilon_1^{}< 1,
\hfill
$ 
\\[1.5ex]
сферу (1.1) разделим на две части и возьмем часть
\\[1.5ex]
\mbox{}\hfill
$
S_1^2=\{(x^{\star},y^{\star},z^{\star})\colon\,
x^{\star}{}^{\,^{\scriptstyle 2}}+
y^{\star}{}^{\,^{\scriptstyle 2}}+
z^{\star}{}^{\,^{\scriptstyle 2}}=1, \ \,
{}-1 \leq z^{\star} \leq 1-\varepsilon_1^{}\},
\hfill
$
\\[1.75ex]
которая не содержит северный полюс $N(0,0,1).$
\vspace{0.35ex}
Число $\varepsilon_1^{}\in (0;1)$ выберем так, чтобы на части $S_1^2$
\vspace{0.15ex}
сферы (1.1) были расположены стереографические образы всех состояний равновесия
и изолированных замкнутых траекторий системы (D), 
\vspace{0.35ex}
лежащих в конечной части  расширенной фазовой плоскости $\overline{Oxy}.$
\vspace{0.25ex}
За круг  $K(x,y)$ примем  круг, лежащий на фазовой плоскости $Oxy$
\vspace{0.35ex}
с центром в начале координат  $O$ и являющийся прообразом части $S_1^2$ сферы  (1.1) 
\vspace{0.35ex}
при стереографическом проецировании с центром проекции в
северном полюсе $N(0,0,1)$ (рис. 5.1).
\\[3.75ex]
\mbox{}\hfill
{\unitlength=1mm
\begin{picture}(106,57)
\put(0,0){\includegraphics[width=106mm,height=57mm]{r05-01.eps}}

\put(53.9,20.7){\makebox(0,0)[cc]{\scriptsize $O$}}
\put(50.7,24.7){\makebox(0,0)[cc]{\scriptsize $S$}}
\put(68.5,21.5){\makebox(0,0)[cc]{\scriptsize $x$}}
\put(73,21){\makebox(0,0)[cc]{\scriptsize $x$}}
\put(43,4.9){\makebox(0,0)[cc]{\scriptsize $y$}}
\put(53.6,29){\makebox(0,0)[cc]{\scriptsize $y$}}
\put(63,6.6){\makebox(0,0)[cc]{\scriptsize $M$}}

\put(50.5,38){\makebox(0,0)[cc]{\scriptsize $O^{\star}$}}
\put(56.9,32.9){\makebox(0,0)[cc]{\scriptsize $P$}}
\put(70.6,34.1){\makebox(0,0)[cc]{\scriptsize $x^{\star}$}}
\put(61,50.3){\makebox(0,0)[cc]{\scriptsize $y^{\star}$}}

\put(50.6,51.5){\makebox(0,0)[cc]{\scriptsize $N$}}
\put(54.7,56.5){\makebox(0,0)[cc]{\scriptsize $z^{\star}$}}

\put(52.5,-6.5){\makebox(0,0)[cc]{\rm Рис. 5.1}}
\end{picture}}
\hfill\mbox{}
\\[7.5ex]
\indent
Плоскостью 
\\[1ex]
\mbox{}\hfill
$
\{(x^{\star},y^{\star},z^{\star})\colon\, 
z^{\star}=\varepsilon_2^{}-1\},
\ \ \
0<\varepsilon_2^{}< 1,
\hfill
$ 
\\[1.5ex]
сферу (1.1) разделим на две части и возьмем часть 
\\[1.5ex]
\mbox{}\hfill
$
S_2^2=\{(x^{\star},y^{\star},z^{\star})\colon\,
x^{\star}{}^{\,^{\scriptstyle 2}}+
y^{\star}{}^{\,^{\scriptstyle 2}}+
z^{\star}{}^{\,^{\scriptstyle 2}}=1, \ \,
\varepsilon_2^{}-1 \leq z^{\star} \leq 1\},
\hfill
$
\\[1.75ex]
которая не содержит южный  полюс $S(0,0,{}-1).$
\vspace{0.35ex}
Число $\varepsilon_2^{}\in (0;1)$ выберем так, чтобы на части $S_2^2$
\vspace{0.15ex}
сферы (1.1) были расположены стереографические образы всех состояний равновесия
и изолированных замкнутых траекторий системы (3.3), лежащих
\vspace{0.35ex}
в конечной части расширенной  фазовой плоскости $\overline{O^{\ast} uv}.$
\vspace{0.35ex}
За круг  $K(u,v)$ примем  круг, лежащий на фазовой плоскости $O^{\ast} uv$
\vspace{0.35ex}
с центром в начале координат  $O^{\ast}$ и являющийся пробразом части $S_2^2$ сферы
\vspace{0.35ex}
 (1.1) при стереографическом проецировании плоскости $O^{\ast} uv$ с центром проекции в
южном полюсе $S(0,0,{}-1)$ (рис. 5.2).
\vspace{0.5ex}

Упорядоченную пару $(K(x,y), K(u,v))$ 
\vspace{0.35ex}
кругов $K(x,y)$ и $K(u,v)$ c нанесенными на них траекториями систем (D) и (3.3) соответственно 
\vspace{0.35ex}
назовем {\it стереографическим атласом траекторий} системы (D).
\vspace{0.35ex}
Тогда (согласно теореме 2.1)
упорядоченная пара $(K(u,v), K(x,y))$ ---  стереографический атлас траекторий системы (3.3).

\newpage

\mbox{}
\\[-1ex]
\mbox{}\hfill
{\unitlength=1mm
\begin{picture}(80,58)
\put(0,0){\includegraphics[width=79.79mm,height=58.58mm]{r05-02.eps}}

\put(41,19.7){\makebox(0,0)[cc]{\scriptsize $O$}}
\put(38.4,23.7){\makebox(0,0)[cc]{\scriptsize $S$}}
\put(56,20.5){\makebox(0,0)[cc]{\scriptsize $x$}}
\put(60.7,20){\makebox(0,0)[cc]{\scriptsize $x$}}
\put(30.8,3.9){\makebox(0,0)[cc]{\scriptsize $y$}}
\put(40.7,27.5){\makebox(0,0)[cc]{\scriptsize $y$}}
\put(49,2.6){\makebox(0,0)[cc]{\scriptsize $M$}}

\put(37.4,37.1){\makebox(0,0)[cc]{\scriptsize $O^{\star}$}}
\put(44.9,32.5){\makebox(0,0)[cc]{\scriptsize $P$}}
\put(58.6,33.8){\makebox(0,0)[cc]{\scriptsize $x^{\star}$}}
\put(42.5,46.7){\makebox(0,0)[cc]{\scriptsize $y^{\star}$}}

\put(37.6,47){\makebox(0,0)[cc]{\scriptsize $N$}}
\put(42,54.5){\makebox(0,0)[cc]{\scriptsize $z^{\star}$}}

\put(40,-6){\makebox(0,0)[cc]{\rm Рис. 5.2}}
\end{picture}}
\hfill\mbox{}
\\[7ex]
\indent
Соответствие между  кругами $K(x,y)$ и $K(u,v)$
\vspace{0.15ex}
показано на рис. 5.3.
Числами $1,\ldots, 40$ отражены соответствия между полуокрестностями точек, 
\vspace{0.15ex}
лежащих на координатных осях и концентрических окружностях с центром в начале координат.
\\[3.5ex]
\mbox{}\hfill
{\unitlength=1mm
\begin{picture}(65,65)
\put(0,0){\includegraphics[width=65mm,height=65mm]{r05-03.eps}}

\put(38.9,34.5){\makebox(0,0)[cc]{\scriptsize $1$}}
\put(38.9,30.8){\makebox(0,0)[cc]{\scriptsize $2$}}
\put(49.9,34.5){\makebox(0,0)[cc]{\scriptsize $3$}}
\put(49.9,30.8){\makebox(0,0)[cc]{\scriptsize $4$}}
\put(59.9,34.5){\makebox(0,0)[cc]{\scriptsize $5$}}
\put(59.9,30.8){\makebox(0,0)[cc]{\scriptsize $6$}}
\put(39.9,39.9){\makebox(0,0)[cc]{\scriptsize $7$}}
\put(42.9,42.8){\makebox(0,0)[cc]{\scriptsize $8$}}
\put(46.9,46.9){\makebox(0,0)[cc]{\scriptsize $9$}}
\put(49.5,49.5){\makebox(0,0)[cc]{\scriptsize $10$}}
\put(34.2,38.7){\makebox(0,0)[cc]{\scriptsize $11$}}
\put(30.5,38.7){\makebox(0,0)[cc]{\scriptsize $12$}}
\put(34.2,49.9){\makebox(0,0)[cc]{\scriptsize $13$}}
\put(30.5,49.9){\makebox(0,0)[cc]{\scriptsize $14$}}
\put(34.2,59.7){\makebox(0,0)[cc]{\scriptsize $15$}}
\put(30.5,59.7){\makebox(0,0)[cc]{\scriptsize $16$}}
\put(25.1,39.9){\makebox(0,0)[cc]{\scriptsize $17$}}
\put(22.7,42.8){\makebox(0,0)[cc]{\scriptsize $18$}}
\put(18.1,47.1){\makebox(0,0)[cc]{\scriptsize $19$}}
\put(15.5,49.9){\makebox(0,0)[cc]{\scriptsize $20$}}
\put(26.3,34.5){\makebox(0,0)[cc]{\scriptsize $21$}}
\put(26.3,31){\makebox(0,0)[cc]{\scriptsize $22$}}
\put(15.2,34.5){\makebox(0,0)[cc]{\scriptsize $23$}}
\put(15.2,31){\makebox(0,0)[cc]{\scriptsize $24$}}
\put(5.2,34.5){\makebox(0,0)[cc]{\scriptsize $25$}}
\put(5.2 ,31){\makebox(0,0)[cc]{\scriptsize $26$}}
\put(25.1,24.9){\makebox(0,0)[cc]{\scriptsize $27$}}
\put(22.5,22.4){\makebox(0,0)[cc]{\scriptsize $28$}}
\put(17.8,18.2){\makebox(0,0)[cc]{\scriptsize $29$}}
\put(15.4,15.5){\makebox(0,0)[cc]{\scriptsize $30$}}
\put(34.2,26.3){\makebox(0,0)[cc]{\scriptsize $32$}}
\put(30.5,26.3){\makebox(0,0)[cc]{\scriptsize $31$}}
\put(34.2,15.1){\makebox(0,0)[cc]{\scriptsize $34$}}
\put(30.5,15.1){\makebox(0,0)[cc]{\scriptsize $33$}}
\put(34.2,5){\makebox(0,0)[cc]{\scriptsize $36$}}
\put(30.5,5){\makebox(0,0)[cc]{\scriptsize $35$}}
\put(40.1,25.1){\makebox(0,0)[cc]{\scriptsize $37$}}
\put(42.7,22.5){\makebox(0,0)[cc]{\scriptsize $38$}}
\put(47.2,18.1){\makebox(0,0)[cc]{\scriptsize $39$}}
\put(49.6,15.4){\makebox(0,0)[cc]{\scriptsize $40$}}
\put(64,31){\makebox(0,0)[cc]{\scriptsize $x$}}
\put(30.5,64){\makebox(0,0)[cc]{\scriptsize $y$}}

\end{picture}
}
\qquad\ \ 
{\unitlength=1mm
\begin{picture}(65,65)
\put(0,0){\includegraphics[width=65mm,height=65mm]{r05-03.eps}}

\put(38.7,34.5){\makebox(0,0)[cc]{\scriptsize $5$}}
\put(38.7,30.8){\makebox(0,0)[cc]{\scriptsize $6$}}
\put(49.9,34.5){\makebox(0,0)[cc]{\scriptsize $3$}}
\put(49.9,30.8){\makebox(0,0)[cc]{\scriptsize $4$}}
\put(59.9,34.5){\makebox(0,0)[cc]{\scriptsize $1$}}
\put(59.9,30.8){\makebox(0,0)[cc]{\scriptsize $2$}}
\put(40,40){\makebox(0,0)[cc]{\scriptsize $10$}}
\put(42.9,42.8){\makebox(0,0)[cc]{\scriptsize $9$}}
\put(46.9,46.9){\makebox(0,0)[cc]{\scriptsize $8$}}
\put(49.5,49.5){\makebox(0,0)[cc]{\scriptsize $7$}}
\put(34.2,38.7){\makebox(0,0)[cc]{\scriptsize $15$}}
\put(30.5,38.7){\makebox(0,0)[cc]{\scriptsize $16$}}
\put(34.2,49.9){\makebox(0,0)[cc]{\scriptsize $13$}}
\put(30.5,49.9){\makebox(0,0)[cc]{\scriptsize $14$}}
\put(34.2,59.7){\makebox(0,0)[cc]{\scriptsize $11$}}
\put(30.5,59.7){\makebox(0,0)[cc]{\scriptsize $12$}}
\put(25.2,40){\makebox(0,0)[cc]{\scriptsize $20$}}
\put(22.7,42.8){\makebox(0,0)[cc]{\scriptsize $19$}}
\put(18.1,47.1){\makebox(0,0)[cc]{\scriptsize $18$}}
\put(15.6,49.7){\makebox(0,0)[cc]{\scriptsize $17$}}
\put(26.3,34.5){\makebox(0,0)[cc]{\scriptsize $25$}}
\put(26.3,31){\makebox(0,0)[cc]{\scriptsize $26$}}
\put(15.2,34.5){\makebox(0,0)[cc]{\scriptsize $23$}}
\put(15.2,31){\makebox(0,0)[cc]{\scriptsize $24$}}
\put(5.2,34.5){\makebox(0,0)[cc]{\scriptsize $21$}}
\put(5.2 ,31){\makebox(0,0)[cc]{\scriptsize $22$}}
\put(25,25.1){\makebox(0,0)[cc]{\scriptsize $30$}}
\put(22.5,22.4){\makebox(0,0)[cc]{\scriptsize $29$}}
\put(17.8,18.2){\makebox(0,0)[cc]{\scriptsize $28$}}
\put(15.4,15.4){\makebox(0,0)[cc]{\scriptsize $27$}}
\put(34.2,26.3){\makebox(0,0)[cc]{\scriptsize $36$}}
\put(30.5,26.3){\makebox(0,0)[cc]{\scriptsize $35$}}
\put(34.2,15.1){\makebox(0,0)[cc]{\scriptsize $34$}}
\put(30.5,15.1){\makebox(0,0)[cc]{\scriptsize $33$}}
\put(34.2,5){\makebox(0,0)[cc]{\scriptsize $32$}}
\put(30.5,5){\makebox(0,0)[cc]{\scriptsize $31$}}
\put(40.1,25.1){\makebox(0,0)[cc]{\scriptsize $40$}}
\put(42.7,22.5){\makebox(0,0)[cc]{\scriptsize $39$}}
\put(47.2,18.1){\makebox(0,0)[cc]{\scriptsize $38$}}
\put(49.7,15.4){\makebox(0,0)[cc]{\scriptsize $37$}}
\put(64,31){\makebox(0,0)[cc]{\scriptsize $u$}}
\put(30.5,64){\makebox(0,0)[cc]{\scriptsize $v$}}

\end{picture}
}
\hfill\mbox{}
\\[1ex]
\mbox{}\hfill
Рис. 5.3
\hfill\mbox{}
\\[3ex]
\indent
{\bf Пример 5.1.}
Траекториями системы (4.4) являются параллельные прямые 
\\[1.5ex]
\mbox{}\hfill
$
b_{_0}x-a_{_0}y = C,
\quad 
C\in\R.
\hfill
$
\\[1.5ex]
\indent
Траекториями системы (4.5) являются состояние равновесия  $O^{\ast} (0,0)$ и
примыкающие к нему кривые 
\\[1.5ex]
\mbox{}\hfill
$
\dfrac{b_{_0}u-a_{_0}v}{u^2+v^2}=C^{\ast},
\quad 
|u|+|v|\ne 0,
\quad 
C^{\ast}\in\R,
\quad 
C^{\ast}=4C.
\hfill
$
\\[1.75ex]
\indent
При $a_{_0}=1,\ b_{_0}=0$ на рис. 5.4 
\vspace{0.25ex}
изображено поведение траекторий на сфере (1.1) стереографически
взаимосопряженных систем (4.4) и (4.5).
\vspace{0.25ex}
Круги, образующие стереографические атласы траекторий систем (4.4) и (4.5), построены на рис. 5.5.
\vspace{0.35ex}

Заметим, что состояние равновесия  $O^{\ast} (0,0)$ 
\vspace{0.25ex}
однородной квадратичной системы (4.5) является сложным и состоит из двух 
\vspace{0.25ex}
эллиптических секторов Бендиксона, ограниченных траекториями-лучами прямой $v=0.$

\newpage

\mbox{}
\\[-1ex]
\parbox{50mm}{
\mbox{}\hfill
{\unitlength=1mm
\begin{picture}(42,42)
\put(0,0){\includegraphics[width=42mm,height=42mm]{r05-04.eps}}
\end{picture}}
\hfill\mbox{}
\\[2ex]
\mbox{}\hfill
Рис. 5.4
\hfill\mbox{}
}
\parbox{104mm}{
\mbox{}\hfill
{\unitlength=1mm
\begin{picture}(42,42)
\put(0,0){\includegraphics[width=42mm,height=42mm]{r05-05a.eps}}
\put(41,19){\makebox(0,0)[cc]{\scriptsize $x$}}
\put(19,41){\makebox(0,0)[cc]{\scriptsize $y$}}
\end{picture}}
\qquad
{\unitlength=1mm
\begin{picture}(42,42)
\put(0,0){\includegraphics[width=42mm,height=42mm]{r05-05b.eps}}
\put(41,19){\makebox(0,0)[cc]{\scriptsize $u$}}
\put(19,41){\makebox(0,0)[cc]{\scriptsize $v$}}
\end{picture}}
\hfill\mbox{}
\\[2ex]
\mbox{}\hfill
Рис. 5.5
\hfill\mbox{}
}
\\[4ex]
\indent
{\bf Пример 5.2.}
Траекториями линейной системы
\\[2ex]
\mbox{}\hfill                        
$
\dfrac{dx}{dt}=x,
\qquad
\dfrac{dy}{dt}=y
$
\hfill (5.1)
\\[2.25ex]
являются $O\!$-лучи семейства прямых $C_1^{}y+C_2^{}x=0,\ C_1^{},C_2^{}\in \R,$
\vspace{0.35ex}
и  состояние равновесия  $O (0,0)$ --- неустойчивый дикритический узел.
\vspace{0.35ex}

Траекториями стереографически сопряженной системы
\\[2ex]
\mbox{}\hfill                        
$
\dfrac{du}{dt}={}-u,
\qquad
\dfrac{dv}{dt}={}-v
$
\hfill (5.2)
\\[2.25ex]
являются $O\!$-лучи семейства прямых $C_1^{}v+C_2^{}u=0,\ C_1^{},C_2^{}\in \R,$
\vspace{0.35ex}
и  состояние равновесия  $O^{\ast} (0,0)$ --- устойчивый дикритический узел.
\vspace{0.35ex}

На рис. 5.6 изображено поведение траекторий на сфере (1.1) стереографически
взаимосопряженных дифференциальных систем (5.1) и (5.2):
стереографическими образами траекторий систем (5.1) и (5.2)
на сфере (1.1) являются полуокружности меридианов сферы (1.1), 
примыкающие к северному и южному полюсам, а также северный и южный полюсы сферы (1.1).
На рис. 5.7 построены круги, образующие стереографические атласы траекторий 
дифференциальных систем (5.1) и (5.2).
\\[4ex]
\parbox{50mm}{
\mbox{}\hfill
{\unitlength=1mm
\begin{picture}(42,42)
\put(0,0){\includegraphics[width=42mm,height=42mm]{r05-06.eps}}
\end{picture}}
\hfill\mbox{}
\\[2ex]
\mbox{}\hfill
Рис. 5.6
\hfill\mbox{}
}
\parbox{104mm}{
\mbox{}\hfill
{\unitlength=1mm
\begin{picture}(42,42)
\put(0,0){\includegraphics[width=42mm,height=42mm]{r05-07a.eps}}
\put(41,19){\makebox(0,0)[cc]{\scriptsize $x$}}
\put(19,41){\makebox(0,0)[cc]{\scriptsize $y$}}
\end{picture}}
\qquad
{\unitlength=1mm
\begin{picture}(42,42)
\put(0,0){\includegraphics[width=42mm,height=42mm]{r05-07b.eps}}
\put(41,19){\makebox(0,0)[cc]{\scriptsize $u$}}
\put(19,41){\makebox(0,0)[cc]{\scriptsize $v$}}
\end{picture}}
\hfill\mbox{}
\\[2ex]
\mbox{}\hfill
Рис. 5.7
\hfill\mbox{}
}
\\[5ex]
\indent
{\bf Пример 5.3.}
Траекториями линейной системы
\\[2ex]
\mbox{}\hfill                        
$
\dfrac{dx}{dt}=y,
\qquad
\dfrac{dy}{dt}={}-x
$
\hfill (5.3)
\\[2.25ex]
являются концентрические окружности $x^2+y^2=C,\ C\in (0;{}+\infty),$
\vspace{0.35ex}
и  состояние равновесия  $O (0,0)$ --- центр.
\vspace{0.5ex}
Направление движения вдоль траекторий системы (5.3) определяют касательные векторы 
$\vec{a}(x,y)=(y,{}-x)\;\; \forall (x,y)\in \R^2\backslash\{(0,0)\}.$

\newpage

Траекториями стереографически сопряженной системы
\\[2ex]
\mbox{}\hfill                        
$
\dfrac{du}{dt}=v,
\qquad
\dfrac{dv}{dt}={}-u
$
\hfill (5.4)
\\[2.25ex]
являются концентрические окружности 
\vspace{0.35ex}
$u^2+v^2=C^\ast,\ C^\ast \in (0;{}+\infty), \ C^\ast=4C,$
и  состояние равновесия  $O^{\ast} (0,0)$ --- центр.
\vspace{0.35ex}

Стереографическими образами 
\vspace{0.15ex}
траекторий систем (5.3) и (5.4)
на сфере (1.1) являются параллели,  северный и южный полюсы сферы (1.1) (рис. 5.8).
\vspace{0.15ex}
На рис. 5.9 построены
круги, образующие стереографические атласы траекторий систем (5.3) и (5.4).
\\[3ex]
\parbox{50mm}{
\mbox{}\hfill
{\unitlength=1mm
\begin{picture}(42,42)
\put(0,0){\includegraphics[width=42mm,height=42mm]{r05-08.eps}}
\end{picture}}
\hfill\mbox{}
\\[2ex]
\mbox{}\hfill
Рис. 5.8
\hfill\mbox{}
}
\parbox{104mm}{
\mbox{}\hfill
{\unitlength=1mm
\begin{picture}(42,42)
\put(0,0){\includegraphics[width=42mm,height=42mm]{r05-09a.eps}}
\put(41,19){\makebox(0,0)[cc]{\scriptsize $x$}}
\put(19,41){\makebox(0,0)[cc]{\scriptsize $y$}}
\end{picture}}
\qquad
{\unitlength=1mm
\begin{picture}(42,42)
\put(0,0){\includegraphics[width=42mm,height=42mm]{r05-09b.eps}}
\put(41,19){\makebox(0,0)[cc]{\scriptsize $u$}}
\put(19,41){\makebox(0,0)[cc]{\scriptsize $v$}}
\end{picture}}
\hfill\mbox{}
\\[2ex]
\mbox{}\hfill
Рис. 5.9
\hfill\mbox{}
}
\\[5ex]
\indent
{\bf Пример 5.4.}
Траекториями линейной системы
\\[2ex]
\mbox{}\hfill                        
$
\dfrac{dx}{dt}=x-y,
\qquad
\dfrac{dy}{dt}=x+y
$
\hfill (5.5)
\\[2ex]
являются логарифмические спирали  
\\[1.75ex]
\mbox{}\hfill
$
(x^2+y^2)\exp \Bigl( {}-\arctg \dfrac{y}{x}\Bigr)=C,
\quad 
C\in (0;{}+\infty),
\hfill
$
\\[1.75ex]
и  состояние равновесия  $O (0,0)$ --- неустойчивый грубый фокус.
\vspace{0.35ex}

Траекториями стереографически сопряженной системы
\\[2ex]
\mbox{}\hfill                        
$
\dfrac{du}{dt}={}-u-v,
\qquad
\dfrac{dv}{dt}=u-v
$
\hfill (5.6)
\\[2ex]
являются логарифмические спирали  
\\[1.75ex]
\mbox{}\hfill
$
(u^2+v^2)\exp \Bigl( {}-\arctg \dfrac{v}{u}\Bigr)=C^\ast,
\quad 
C^\ast\in (0;{}+\infty),
\quad 
C^\ast=4C,
\hfill
$
\\[1.75ex]
и  состояние равновесия  $O^{\ast} (0,0)$ --- устойчивый грубый фокус.
\\[4ex]
\parbox{50mm}{
\mbox{}\hfill
{\unitlength=1mm
\begin{picture}(42,42)
\put(0,0){\includegraphics[width=42mm,height=42mm]{r05-10.eps}}
\end{picture}}
\hfill\mbox{}
\\[2ex]
\mbox{}\hfill
Рис. 5.10
\hfill\mbox{}
}
\parbox{104mm}{
\mbox{}\hfill
{\unitlength=1mm
\begin{picture}(42,42)
\put(0,0){\includegraphics[width=42mm,height=42mm]{r05-11a.eps}}
\put(41,19){\makebox(0,0)[cc]{\scriptsize $x$}}
\put(19,41){\makebox(0,0)[cc]{\scriptsize $y$}}
\end{picture}}
\qquad
{\unitlength=1mm
\begin{picture}(42,42)
\put(0,0){\includegraphics[width=42mm,height=42mm]{r05-11b.eps}}
\put(41,19){\makebox(0,0)[cc]{\scriptsize $u$}}
\put(19,41){\makebox(0,0)[cc]{\scriptsize $v$}}
\end{picture}}
\hfill\mbox{}
\\[2ex]
\mbox{}\hfill
Рис. 5.11
\hfill\mbox{}
}
\\[-3ex]

\newpage

На рис. 5.10 изображены траектории  на сфере (1.1) 
\vspace{0.15ex}
стереографически вза\-и\-мо\-со\-п\-ря\-жен\-ных систем (5.5) и (5.6).
\vspace{0.15ex}
На рис. 5.11 построены
круги, образующие стереографические атласы траекторий систем (5.5) и (5.6).
\vspace{0.75ex}

{\bf Пример 5.5.}
Траекториями линейной системы
\\[2ex]
\mbox{}\hfill                        
$
\dfrac{dx}{dt}=x,
\qquad
\dfrac{dy}{dt}={}-y
$
\hfill (5.7)
\\[2.25ex]
являются $O\!$-кривые семейства  $xy=C,\ C\in \R,$
\vspace{0.25ex}
и  состояние равновесия  $O (0,0)$ --- седло, сепаратрисами которого
являются координатные $O\!$-лучи.
\vspace{0.35ex}

Траекториями стереографически сопряженной системы
\\[2ex]
\mbox{}\hfill                        
$
\dfrac{du}{d\tau}={}-u^3+3uv^2,
\qquad
\dfrac{dv}{d\tau}={}-3u^2\;\!v+v^3, 
\qquad 
(u^2+v^2)\, d\tau=dt,
$
\hfill (5.8)
\\[2ex]
являются $O^{\ast}\!$-кривые семейства  
\\[1.75ex]
\mbox{}\hfill
$
\dfrac{uv}{(u^2+v^2)^2 }=C^{\ast},
\quad 
C^{\ast}\in \R,
\quad 
C^{\ast}=16C,
\hfill
$
\\[1.75ex]
и сложное  состояние равновесия  $O^{\ast} (0,0),$ 
\vspace{0.15ex}
которое состоит из четырех эллиптических секторов Бендексона, ограниченных координатными 
$O^{\ast}\!$-лучами, и образованных лемнискатами Бернулли.

На рис. 5.12
изображено поведение траекторий на сфере (1.1) стереографически
взаимосопряженных систем (5.7) и (5.8).
На рис. 5.13 построены
круги, образующие стереографические атласы траекторий систем (5.7) и (5.8).
\\[3.75ex]
\mbox{}\hfill
{\unitlength=1mm
\begin{picture}(42,42)
\put(0,0){\includegraphics[width=42mm,height=42mm]{r05-12a.eps}}
\end{picture}}
\qquad
\qquad
\qquad
{\unitlength=1mm
\begin{picture}(42,42)
\put(0,0){\includegraphics[width=42mm,height=42mm]{r05-12b.eps}}
\end{picture}}
\hfill\mbox{}
\\[2ex]
\mbox{}\hfill
Рис. 5.12
\hfill\mbox{}
\\[5.75ex]
\mbox{}\hfill
{\unitlength=1mm
\begin{picture}(42,42)
\put(0,0){\includegraphics[width=42mm,height=42mm]{r05-13a.eps}}
\put(41,19){\makebox(0,0)[cc]{\scriptsize $x$}}
\put(19,41){\makebox(0,0)[cc]{\scriptsize $y$}}
\end{picture}}
\qquad
\qquad
{\unitlength=1mm
\begin{picture}(42,42)
\put(0,0){\includegraphics[width=42mm,height=42mm]{r05-13b.eps}}
\put(41,19){\makebox(0,0)[cc]{\scriptsize $u$}}
\put(19,41){\makebox(0,0)[cc]{\scriptsize $v$}}
\end{picture}}
\hfill\mbox{}
\\[2ex]
\mbox{}\hfill
Рис. 5.13
\hfill\mbox{}
\\[3.75ex]
\indent
{\bf Пример 5.6.}
Траекториями линейной системы
\\[2ex]
\mbox{}\hfill                        
$
\dfrac{dx}{dt}=x,
\qquad
\dfrac{dy}{dt}=2y
$
\hfill (5.9)
\\[2.5ex]
являются $O\!$-кривые семейства  
\\[1.5ex]
\mbox{}\hfill
$
C_1^{}y+C_2^{}x^2=0,
\quad  
C_1^{},C_2^{}\in \R,
\hfill
$
\\[1.5ex]
и  состояние равновесия  $O (0,0)$ --- простой неустойчивый узел.
\vspace{0.35ex}

Траекториями стереографически сопряженной системы
\\[2ex]
\mbox{}\hfill                        
$
\dfrac{du}{d\tau}={}-u^3-3uv^2,
\qquad
\dfrac{dv}{d\tau}={}-2v^3, 
\qquad 
(u^2+v^2)\, d\tau=dt,
$
\hfill (5.10)
\\[2.5ex]
являются $O^{\ast}\!$-кривые семейства
\\[1.5ex]
\mbox{}\hfill
$
C_1^{\ast}v(u^2+v^2)+C_2^{\ast}\, u^2=0,
\quad  
C_1^{\ast},C_2^{\ast}\in \R,
\quad
C_1^{\ast}=C_1^{},
\quad 
C_2^{\ast}=4C_2^{},
\hfill
$
\\[1.5ex]
и сложное  состояние равновесия  $O^{\ast} (0,0)$ --- устойчивый узел.
\vspace{0.35ex}

На рис. 5.14 изображены 
\vspace{0.15ex}
траектории  на сфере (1.1) стереографически взаимосопряженных систем (5.9) и (5.10).
\vspace{0.15ex}
На рис. 5.15 построены
круги, образующие стереографические атласы траекторий систем (5.9) и (5.10).
\\[4ex]
\parbox{50mm}{
\mbox{}\hfill
{\unitlength=1mm
\begin{picture}(42,42)
\put(0,0){\includegraphics[width=42mm,height=42mm]{r05-14.eps}}
\end{picture}}
\hfill\mbox{}
\\[2ex]
\mbox{}\hfill
Рис. 5.14
\hfill\mbox{}
}
\parbox{104mm}{
\mbox{}\hfill
{\unitlength=1mm
\begin{picture}(42,42)
\put(0,0){\includegraphics[width=42mm,height=42mm]{r05-15a.eps}}
\put(41,19){\makebox(0,0)[cc]{\scriptsize $x$}}
\put(19,41){\makebox(0,0)[cc]{\scriptsize $y$}}
\end{picture}}
\qquad
{\unitlength=1mm
\begin{picture}(42,42)
\put(0,0){\includegraphics[width=42mm,height=42mm]{r05-15b.eps}}
\put(41,19){\makebox(0,0)[cc]{\scriptsize $u$}}
\put(19,41){\makebox(0,0)[cc]{\scriptsize $v$}}
\end{picture}}
\hfill\mbox{}
\\[2ex]
\mbox{}\hfill
Рис. 5.15
\hfill\mbox{}
}
\\[5ex]
\indent
{\bf Пример 5.7.}
Траекториями линейной системы
\\[2ex]
\mbox{}\hfill                        
$
\dfrac{dx}{dt}=x+y,
\qquad
\dfrac{dy}{dt}=y
$
\hfill (5.11)
\\[2ex]
являются $O\!$-кривые семейства  
\\[1.5ex]
\mbox{}\hfill
$
y \exp \Bigl({}-\dfrac{x}{y}\Bigr)=C,
\quad 
C\in \R,
\hfill
$
\\[1.5ex]
и  состояние равновесия  $O (0,0)$ --- неустойчивый вырожденный узел.
\vspace{0.35ex}

Траекториями стереографически сопряженной системы
\\[2ex]
\mbox{}\hfill                        
$
\dfrac{du}{d\tau}={}-u^3-u^2\;\!v-uv^2+v^3,
\quad
\dfrac{dv}{d\tau}={}-u^2\;\!v -2uv^2-v^3, \ \ (u^2+v^2)\, d\tau=dt,
$
\hfill (5.12)
\\[2ex]
являются $O^{\ast}\!$-кривые семейства  
\\[1.5ex]
\mbox{}\hfill
$
\dfrac{v}{u^2+v^2 }\, \exp\Bigl({}- \dfrac{u}{v}\Bigr)=C^{\ast},
\quad
C^{\ast}\in\R, 
\quad 
C^{\ast}=4C,
\hfill
$
\\[1.5ex]
и сложное  состояние равновесия  $O^{\ast} (0,0)$ --- устойчивый узел.
\vspace{0.35ex}

На рис. 5.16
\vspace{0.15ex}
изображены траектории  на сфере (1.1) стереографически
взаимосопряженных систем (5.11) и (5.12).
\vspace{0.15ex}
На рис. 5.17 построены
круги, образующие стереографические атласы траекторий систем (5.11) и (5.12).

\newpage

\mbox{}
\\[-1ex]
\parbox{50mm}{
\mbox{}\hfill
{\unitlength=1mm
\begin{picture}(42,42)
\put(0,0){\includegraphics[width=42mm,height=42mm]{r05-16.eps}}
\end{picture}}
\hfill\mbox{}
\\[2ex]
\mbox{}\hfill
Рис. 5.16
\hfill\mbox{}
}
\parbox{104mm}{
\mbox{}\hfill
{\unitlength=1mm
\begin{picture}(42,42)
\put(0,0){\includegraphics[width=42mm,height=42mm]{r05-17a.eps}}
\put(41,19){\makebox(0,0)[cc]{\scriptsize $x$}}
\put(19,41){\makebox(0,0)[cc]{\scriptsize $y$}}
\end{picture}}
\quad
{\unitlength=1mm
\begin{picture}(42,42)
\put(0,0){\includegraphics[width=42mm,height=42mm]{r05-17b.eps}}
\put(41,19){\makebox(0,0)[cc]{\scriptsize $u$}}
\put(19,41){\makebox(0,0)[cc]{\scriptsize $v$}}
\end{picture}}
\hfill\mbox{}
\\[2ex]
\mbox{}\hfill
Рис. 5.17
\hfill\mbox{}
}
\\[7.5ex]
\centerline{
{\bf\large \S\;\!3. Траектории стереографически сопряженных}} 
\\
\centerline{{\bf\large дифференциальных систем}}
\\[2.25ex]
\centerline{
{\bf  6. Регулярные точки и состояния равновесия }}
\\
\centerline{
{\bf  стереографически сопряженных дифференциальных систем}}
\\[1.5ex]
\indent
При диффеоморфном отображении (2.6) 
\vspace{0.15ex}
образом точки $M(x,y)$ фазовой плоскости $Oxy$ системы (D), 
отличной от начала координат $O(0,0),$ является точка 
\\[1.5ex]
\mbox{}\hfill
$
M^* \Bigl(\;\!\dfrac{4x}{x^{2}+y^{2}}\,,\
\dfrac{4y}{x^{2}+y^{2}}\Bigr)
\hfill
$
\\[1.5ex]
фазовой плоскости $O^*uv$ системы (3.3), отличная от начала координат $O^*(0,0).$
\vspace{0.25ex}

Пусть 
\vspace{0.25ex}
образом начала координат $O(0,0)$ расширенной плоскости $\overline{Oxy}$ 
является бесконечно удаленная точка $M^{\ast}_{\infty}$ 
\vspace{0.25ex}
стереографически сопряженной  расширенной  плоскости $\overline{O^{\ast}uv}\,,$
а образом бесконечно удаленной точки 
\vspace{0.15ex}
$M^{}_{\infty}$ расширенной плоскости $\overline{Oxy}$ 
является начало координат $O^{\ast}(0,0)$
\vspace{0.25ex}
стереографически сопряженной расширенной  плоскости $\overline{O^{\ast}uv}\,.$ 
Такое продолжение отображения $\varphi_{21}^{}$  
\vspace{0.35ex}
есть биективное отображение $\overline{\varphi_{21}^{}}$ 
расширенной фазовой плоскости  $\overline{Oxy}$ системы (D) 
\vspace{0.35ex}
на расширенную  фазовую плоскость $\overline{O^{\ast}uv}$ системы (3.3).
\vspace{0.35ex}
Точки $M$ и $M^{\ast},\ O$ и $M^{\ast}_{\infty}\,,\ M^{}_{\infty}$ и $O^{\ast}$ являются стереографически 
взаимосопряженными.
Также при отображении $\overline{\varphi_{21}^{}}$ 
\vspace{0.5ex}
кривая $l$ на расширенной плоскости $\overline{Oxy}$ и 
ее образ $l^\ast$ на расширенной плоскости $\overline{O^{\ast}uv}$ будут стереографически взаимосопряженными кривыми.

Согласно основному свойству стереографической проекции (пункт 1) и тому, 
что отображение (2.6) есть суперпозиция (2.5) стереографических отображений, 
\vspace{0.5ex}
получаем 

{\bf Свойство 6.1.}
{\it
Угол между кривыми равен углу между стереографически сопряженными к ним кривыми.}
\vspace{0.5ex}

Точка фазовой плоскости $Oxy$ 
\vspace{0.15ex}
является регулярной точкой системы (D), если она не является состоянием равновесия системы (D). 
\vspace{0.35ex}
Бесконечно удаленную точку $M_{\infty}^{}$ расширенной фазовой плоскости $\overline{Oxy}$ 
\vspace{0.25ex}
будем считать {\it регулярной бесконечно удаленной точкой} системы (D),  
\vspace{0.25ex}
если начало координат $O^{\ast}(0,0)$ фазовой плоскости 
$O^{\ast}uv$ не является состоянием равновесия системы (3.3). 
\vspace{0.35ex}
Если точка $O^{\ast}(0,0)$ есть состояние равновесия системы (3.3), то бесконечно удаленную точку $M^{}_{\infty}$ 
\vspace{0.25ex}
расширенной фазовой плоскости $\overline{Oxy}$  
примем за бесконечно удаленное состояние равновесия системы (D) такого же вида.
\vspace{0.5ex}

Таким образом, 
\vspace{0.15ex}
каждая точка расширенной фазовой плоскости $\overline{Oxy}$ --- регулярная  или состояние равновесия 
дифференциальной системы (D). 
\vspace{0.35ex}
На  расширенной фазовой плоскости $\overline{Oxy}$ каждая траектория 
\vspace{0.15ex}
дифференциальной системы (D) является состоянием равновесия или состоит из регулярных точек. 

\newpage

На основании диффеоморфности отображения (2.6) плоскости $Oxy,$ проколотой в начале координат $O(0,0),$ 
на стереографически сопряженную плоскость $O^{\ast}uv,$ проколотую в начале координат $O^{\ast}(0,0),$ 
\vspace{0.35ex}
и учитывая свойство 6.1, 
получаем, что на расширенных фазовых плоскостях $\overline{Oxy}$ и $\overline{O^{\ast}uv}$ траектории стереографически 
взаимосопряженных систем (D) и (3.3) такие, что имеет место
\vspace{0.5ex}

{\bf Свойство 6.2.}
\vspace{0.25ex}
{\it
Стереографически взаимосопряженные точки  расширенных фазовых плоскостей $\overline{Oxy}$ и $\overline{O^{\ast}uv}$ 
\vspace{0.15ex}
являются одновременно или регулярными  или состояниями равновесия одного и того же вида систем 
{\rm (D)} и {\rm (3.3)}.
}
\vspace{0.5ex}

Из того, что отображение (2.6) есть суперпозиция (2.5), следуют  
\vspace{0.5ex}

{\bf Свойство 6.3 (6.4).}
{\it
Образом  замкнутой на плоскости $Oxy$ кривой, проходящей 
\linebreak
{\rm(}не проходящей{\rm)} через начало координат $O(0,0),$  является  незамкнутая  
{\rm(}замкнутая{\rm)} кривая на стереографически сопряженной
плоскости $O^{\ast}uv,$ не проходящая через начало координат $O^{\ast}(0,0).$
}
\vspace{0.5ex}

Например, вычислениями, используя преобразование Бендиксона (3.1), получаем 
\vspace{0.5ex}

{\bf Свойство 6.5.}
{\it
Образом\;\!{\rm:}

а{\rm)}
проходящей через  начало координат $O(0,0)$ плоскости $Oxy$ окружности 
\\[1.5ex]
\mbox{}\hfill
$
(x+a)^2+(y+b)^2=a^2+b^2,
\quad 
|a|+|b|\ne 0;
\hfill
$
\\[1.5ex]
\indent
б{\rm)}
не проходящей через  начало координат $O(0,0)$ плоскости $Oxy$ окружности 
\\[1.5ex]
\mbox{}\hfill
$
(x+a)^2+(y+b)^2=r^2,
\quad 
r>0,
\quad 
r^2\ne a^2+b^2;
\hfill
$
\\[1.5ex]
\indent
в{\rm)}
окружности 
\\[0.5ex]
\mbox{}\hfill
$
x^2+y^2=r^2,
\quad 
r>0,
\hfill
$
\\[1.5ex]
с центром в начале координат $O(0,0)$ плоскости $Oxy;$
\vspace{0.5ex}

г{\rm)}
точки $A(a,b),\ |a|+|b|\ne 0,$ если она лежит внутри {\rm(}вне{\rm)}  окружности 
\\[1.5ex]
\mbox{}\hfill
$
x^2+y^2=r^2,
\quad 
r>0;
\hfill
$
\\[1.5ex]
\indent
д{\rm)}
проходящей через  начало координат $O(0,0)$ плоскости $Oxy$ прямой 
\\[1.5ex]
\mbox{}\hfill
$
Ax+By=0,
\quad 
|A|+|B|\ne 0;
\hfill
$
\\[1.5ex]
\indent
е{\rm)}
координатной прямой $x=0\ (y=0)$ плоскости $Oxy;$
\vspace{0.5ex}

ж{\rm)}
не проходящей через  начало координат $O(0,0)$ плоскости $Oxy$ прямой 
\\[1.5ex]
\mbox{}\hfill
$
Ax+By+C=0,
\quad 
|A|+|B|\ne 0,
\quad 
C\ne 0;
\hfill
$
\\[1.5ex]
на стереографически сопряженной плоскости $O^\ast uv$ является\;\!{\rm:}
\vspace{0.35ex}

а{\rm)}
прямая $au+bv+2=0;$
\vspace{0.5ex}

б{\rm)}
 окружность  
$\Bigl( u+\dfrac{4a}{a^2+b^2-r^2}\Bigr)^{\!2}+\Bigl( v+\dfrac{4b}{a^2+b^2-r^2}\Bigr)^{\!2}=\dfrac{4r^2}{(a^2+b^2-r^2)^2}\,;$
\vspace{0.75ex}

в{\rm)}
окружность  $u^2+v^2=\dfrac{4}{r^2}\,;$
\vspace{0.35ex}

г{\rm)}
стереографически сопряженная точка 
$A^\ast \biggl( \dfrac{4a}{a^2+b^2}\,,\, \dfrac{4b}{a^2+b^2}\biggr),$
лежащая вне {\rm(}внутри{\rm)}  
ок\-руж\-но\-сти $u^2+v^2=\dfrac{4}{r^2}\,;$
\vspace{0.75ex}

д{\rm)}
проходящая через начало координат $O^\ast(0,0)$ 
прямая $Au+Bv=0;$
\vspace{0.35ex}

е{\rm)}
координатная прямая $u=0\ (v=0);$
\vspace{0.35ex}

ж{\rm)}
проходящая через начало координат $O^\ast(0,0)$  окружность  
\\[1.5ex]
\mbox{}\hfill
$
\Bigl( u+\dfrac{2A}{C}\Bigr)^{\!2}+\Bigl(v+\dfrac{2B}{C}\Bigr)^{\!2}=\dfrac{4(A^2+B^2)}{C^2}\,.
\hfill
$
\\[-1.5ex]
}

\newpage

В примере 5.2 (рис. 5.6 и 5.7) рассмотрены стереографически взаимосопряженные дифференциальных системы (5.1) и (5.2) 
с траекториями, лежащими на прямых, проходящих через  начала координат $O(0,0)$ и $O^\ast(0,0)$ 
фазовых плоскостей $Oxy$ и $O^\ast uv$ соответственно (свойство 6.5, случай д). 
\vspace{0.5ex}

{\bf Пример 6.1.}
Заменой $x$ на $x-1,\ y$ на $y-1$ систему (5.1) приводим к системе 
\\[1.5ex]
\mbox{}\hfill                        
$
\dfrac{dx}{dt}=x-1,
\qquad
\dfrac{dy}{dt}=y-1.
$
\hfill (6.1)
\\[1.75ex]
\indent
Траекториями 
\vspace{0.25ex}
линейной дифференциальной системы (6.1) являются неустойчивый дикритический узел $A(1,1)$
и $A\!$-лучи пучка прямых 
\\[1ex]
\mbox{}\hfill
$
\dfrac{y-1}{x-1}=C,
\quad 
{}-\infty\leq C\leq {}+\infty.
\hfill
$ 
\\[1.5ex]
\indent
Траекториями стереографически сопряженной системы 
\\[2ex]
\mbox{}\hfill                        
$
\dfrac{du}{dt}={}-u+\dfrac{1}{4}\; u^2+\dfrac{1}{2}\; uv-\dfrac{1}{4}\; v^2,
\qquad
\dfrac{dv}{dt}={}-v-\dfrac{1}{4}\; u^2+\dfrac{1}{2}\; uv+\dfrac{1}{4}\; v^2
$
\hfill (6.2)
\\[2.25ex]
являются 
\vspace{0.5ex}
дикритические узлы $O^\ast (0,0)$ (устойчивый) и $A^\ast (2,2)$ (неустойчивый),
отрезок $O^\ast\, A^\ast$ без концов,
$O^\ast\!$-луч и $A^\ast\!$-луч прямой $v=u,$ 
дуги окружностей 
\\[1.75ex]
\mbox{}\hfill
$
\Bigl( u-\dfrac{2C}{C-1}\Bigr)^{\!2}+\Bigl( v+\dfrac{2}{C-1}\Bigr)^{\!2}=
\dfrac{4(C^2+1)}{(C-1)^2}\,,
\quad 
C\in \R\backslash\{1\},
\hfill
$
\\[1.75ex]
примыкающие к состояниям равновесия  $O^\ast$ и $A^\ast$
\vspace{0.15ex}
(центры окружностей лежат на прямой $u+v-2=0).$
\vspace{0.15ex}

Круги, образующие стереографические атласы траекторий систем (6.1) и (6.2),
построены на рис. 6.1.
\\[3.5ex]
\mbox{}\hfill
{\unitlength=1mm
\begin{picture}(45,45)
\put(0,0){\includegraphics[width=45mm,height=45mm]{r06-01a.eps}}
\put(44,20.5){\makebox(0,0)[cc]{\scriptsize $x$}}
\put(20.5,44){\makebox(0,0)[cc]{\scriptsize $y$}}
\end{picture}}
\qquad\qquad
{\unitlength=1mm
\begin{picture}(45,45)
\put(0,0){\includegraphics[width=45mm,height=45mm]{r06-01b.eps}}
\put(44,20.5){\makebox(0,0)[cc]{\scriptsize $u$}}
\put(20.5,44){\makebox(0,0)[cc]{\scriptsize $v$}}
\end{picture}}
\hfill\mbox{}
\\[1.75ex]
\mbox{}\hfill
Рис. 6.1
\hfill\mbox{}
\\[2.5ex]
\indent
В примере 5.3 (рис. 5.8 и 5.9) рассмотрены стереографически взаимосопряженные системы (5.3) и (5.4),
траекториями которых являются концентрические окружности, стягивающиеся в центры 
$O$ и $O^\ast$ соответственно (свойство 6.5, случай в).
\vspace{0.5ex}

{\bf Пример 6.2.}
Заменой $x$ на $x-1,\ y$ на $y-1$ систему (5.3) приводим к системе 
\\[1.5ex]
\mbox{}\hfill                        
$
\dfrac{dx}{dt}=y-1,
\qquad
\dfrac{dy}{dt}={}-x+1.
$
\hfill (6.3)
\\[1.5ex]
\indent
Траекториями  системы (6.3) являются центр  $A(1,1)$ и концентрические окружности  
\\[1.5ex]
\mbox{}\hfill
$
(x-1)^2\,+\,(y-1)^2=C,
\quad 
C\in (0;{}+\infty).
\hfill
$ 
\\[1.5ex]
\indent
Траекториями стереографически сопряженной системы 
\\[1.5ex]
\mbox{}\hfill                        
$
\dfrac{du}{dt}=v+\dfrac{1}{4}\; u^2-\dfrac{1}{2}\; uv-\dfrac{1}{4}\; v^2,
\qquad
\dfrac{dv}{dt}={}-u+\dfrac{1}{4}\; u^2+\dfrac{1}{2}\; uv-\dfrac{1}{4}\; v^2
$
\hfill (6.4)
\\[3ex]
являются прямая $u+v-2=0,$ окружности  
\\[1.5ex]
\mbox{}\hfill
$
\Bigl( u-\dfrac{4}{2-C}\Bigr)^{\!2}+\Bigl( v-\dfrac{4}{2-C}\Bigr)^{\!2}=\dfrac{16\;\!C}{(2-C)^2}\,,
\quad 
C\in (0;2)\sqcup (2;{}+\infty),
\hfill
$
\\[1.5ex]
центры $O^\ast (0,0)$ и $A^\ast (2,2).$
\vspace{0.25ex}

Круги, образующие стереографические атласы траекторий систем (6.3) и (6.4),
построены на рис. 6.2.
\\[2.75ex]
\mbox{}\hfill
{\unitlength=1mm
\begin{picture}(45,45)
\put(0,0){\includegraphics[width=45mm,height=45mm]{r06-02a.eps}}
\put(44,20.5){\makebox(0,0)[cc]{\scriptsize $x$}}
\put(20.5,44){\makebox(0,0)[cc]{\scriptsize $y$}}
\end{picture}}
\qquad\qquad
{\unitlength=1mm
\begin{picture}(45,45)
\put(0,0){\includegraphics[width=45mm,height=45mm]{r06-02b.eps}}
\put(44,20.5){\makebox(0,0)[cc]{\scriptsize $u$}}
\put(20.5,44){\makebox(0,0)[cc]{\scriptsize $v$}}
\end{picture}}
\hfill\mbox{}
\\[0ex]
\mbox{}\hfill
Рис. 6.2
\hfill\mbox{}
\\[3ex]
\indent
{\bf  Пример 6.3.}
Траекториями системы Якоби
\\[2ex]
\mbox{}\hfill        
$
\dfrac{dx}{dt}=1+x-y+x(x+y-1),
\qquad 
\dfrac{dy}{dt}=y(x+y-1)
$
\hfill(6.5)
\\[2ex]
являются кривые семейства 
\\[1.5ex]
\mbox{}\hfill
$
\dfrac{x^2+(y-1)^2}{y^2}\,\exp \Bigl({}-2\arctg \dfrac{y-1}{x}\Bigr)=C,
\quad
0\leq C\leq {}+\infty,
\hfill
$
\\[1.5ex]
среди которых прямая $y=0$ и неустойчивый фокус $A(0,1).$
\\[2.75ex]
\mbox{}\hfill
{\unitlength=1mm
\begin{picture}(45,45)
\put(0,0){\includegraphics[width=45mm,height=45mm]{r06-03a.eps}}
\put(44,20.5){\makebox(0,0)[cc]{\scriptsize $x$}}
\put(20.5,44){\makebox(0,0)[cc]{\scriptsize $y$}}
\end{picture}}
\qquad\qquad
{\unitlength=1mm
\begin{picture}(45,45)
\put(0,0){\includegraphics[width=45mm,height=45mm]{r06-03b.eps}}
\put(44,20.5){\makebox(0,0)[cc]{\scriptsize $u$}}
\put(20.5,44){\makebox(0,0)[cc]{\scriptsize $v$}}
\end{picture}}
\hfill\mbox{}
\\[0ex]
\mbox{}\hfill
Рис. 6.3
\hfill\mbox{}
\\[3ex]
\indent
Траекториями стереографически сопряженной системы 
\\[1.75ex]
\mbox{}\hfill        
$
\dfrac{du}{d\theta}={}-4u^2-4uv+u^2v+2uv^2-v^3-\dfrac{1}{4}\, u^4+\dfrac{1}{4}\, v^4,
\hfill
$
\\
\mbox{}\hfill (6.6)
\\
\mbox{}\hfill
$
\dfrac{dv}{d\theta}=v\Bigl({}-4u-4v-u^2+2uv+v^2-\dfrac{1}{2}\, u^3-\dfrac{1}{2}\,uv^2\Bigr),
\hfill
$
\\[2.5ex]
где $(u^2+v^2)\;\!d\theta=dt,$
являются кривые семейства 
\\[2ex]
\mbox{}\hfill
$
\dfrac{16u^2+(u^2+v^2-4v)^2}{v^2}\,\exp \Bigl({}-2\arctg \dfrac{4v-u^2-v^2}{4u}\Bigr)=C^\ast,
\ \ 
0\leq C^\ast\leq {}+\infty,\ C^\ast=16C,
\hfill
$
\\[2.25ex]
среди которых $O^\ast\!$-лучи прямой $v=0,$
неустойчивый фокус $A^\ast(0,2)$
и сложное состояние равновесия $O^\ast(0,0),$
состоящее из гиперболического, эллиптического и сопровождающих его двух параболических секторов Бендиксона.

Круги, образующие стереографические атласы траекторий систем (6.5) и (6.6),
построены на рис. 6.3.
\vspace{0.5ex}

{\bf  Пример 6.4.}
Заменой $x$ на $x,\ y$ на $y+1$ систему (6.5) приводим к системе Якоби
\\[2ex]
\mbox{}\hfill        
$
\dfrac{dx}{dt}=x-y+x(x+y),
\qquad 
\dfrac{dy}{dt}=(y+1)(x+y).
$
\hfill(6.7)
\\[2.5ex]
\indent
Кривые семейства 
\\[1.5ex]
\mbox{}\hfill
$
\dfrac{x^2+y^2}{(y-1)^2}\,\exp \Bigl({}-2\arctg \dfrac{y}{x}\Bigr)=C,
\quad
0\leq C\leq {}+\infty,
\hfill
$
\\[1.5ex]
являются траекториями системы (6.7),
\vspace{0.15ex}
среди которых траектория-прямая $y={}-1$ и неустойчивый фокус $O(0,0).$
\vspace{0.25ex}

Траекториями стереографически сопряженной системы 
\\[2ex]
\mbox{}\hfill        
$
\dfrac{du}{d\theta}={}-(4u^2+4uv+u^3+u^2v+uv^2+v^3),
\quad 
\dfrac{dv}{d\theta}={}-4uv-4v^2+u^3-u^2v+uv^2-v^3,
$
\hfill(6.8)
\\[2.5ex]
где $(u^2+v^2)\;\!d\theta=dt,$ являются кривые семейства 
\\[2ex]
\mbox{}\hfill
$
\dfrac{u^2+v^2}{(u^2+v^2-4v)^2}\,\exp \Bigl({}-2\arctg \dfrac{v}{u}\Bigr)=C^\ast,
\quad
0\leq C^\ast\leq {}+\infty,
\quad 
16C^\ast=C,
\hfill
$
\\[2ex]
среди которых сложное состояние равновесия $O^\ast(0,0),$ 
\vspace{0.15ex}
состоящее из гиперболического, эллиптического и двух сопровождающих его параболических секторов Бендиксона.
\vspace{0.25ex}

Круги, образующие стереографические атласы траекторий систем (6.7) и (6.8),
построены на рис. 6.4.
\\[3.75ex]
\mbox{}\hfill
{\unitlength=1mm
\begin{picture}(45,45)
\put(0,0){\includegraphics[width=45mm,height=45mm]{r06-04a.eps}}
\put(44,20.5){\makebox(0,0)[cc]{\scriptsize $x$}}
\put(20.5,44){\makebox(0,0)[cc]{\scriptsize $y$}}
\end{picture}}
\qquad\qquad
{\unitlength=1mm
\begin{picture}(45,45)
\put(0,0){\includegraphics[width=45mm,height=45mm]{r06-04b.eps}}
\put(44,20.5){\makebox(0,0)[cc]{\scriptsize $u$}}
\put(20.5,44){\makebox(0,0)[cc]{\scriptsize $v$}}
\end{picture}}
\hfill\mbox{}
\\[2ex]
\mbox{}\hfill
Рис. 6.4
\hfill\mbox{}
\\[4.75ex]
\centerline{
{\bf  7. Стереографические циклы
}
}
\\[1.35ex]
\indent
На основании свойства 6.4 и диффеоморфности отбражения (2.6) получаем
\vspace{0.35ex}

{\bf Свойство 7.1 (7.2).}
{\it
Образом не проходящего через начало координат $O(0,0)$ фазовой плоскости $Oxy$ цикла 
{\rm(}предельного цикла{\rm)} дифференциальной системы {\rm(D)}  на  фазовой плоскости 
$O^{\ast}uv$ является не проходящий через начало координат $O^{\ast}(0,0)$ цикл 
{\rm(}предельный цикл{\rm)} дифференциальной системы {\rm(3.3).}
}
\vspace{0.5ex}

На основании свойства 6.3 и диффеоморфности отбражения (2.6) получаем
\vspace{0.35ex}

{\bf Свойство 7.3 (7.4).}
{\it
Образом проходящего через начало координат $O(0,0)$ фазовой плоскости $Oxy$ цикла 
{\rm(}предельного цикла{\rm)} дифференциальной системы {\rm(D)}  на  фазовой плоскости 
$O^{\ast}uv$ является не проходящая через начало координат $O^{\ast}(0,0)$ незамкнутая траектория  
дифференциальной системы {\rm(3.3).}
}
\vspace{0.5ex}

{\bf Определение 7.1.}
Траекторию системы (D), стереографическим образом которой на сфере (1.1) является замкнутая кривая, 
каждая точка которой является образом регулярной точки (конечной или бесконечно удаленной) системы (D), 
назовем  {\it стереографическим циклом} системы (D).
\vspace{0.5ex}

{\bf Определение 7.2.}
{\it Предельным стереографическим циклом} системы (D) назовем 
такой ее стереографический цикл, у стереографического образа которого на сфере (1.1) 
существует окрестность, в которой нет стереографического образа другого стереографического цикла системы (D).
\vspace{0.5ex}

{\bf Свойство 7.5 (7.6).}
{\it
Цикл {\rm(}предельный цикл{\rm)} системы {\rm(D)}  является ее стереографическим циклом 
{\rm(}предельным стереографическим циклом{\rm)}.
}
\vspace{0.5ex}

{\bf Определение 7.3.}
\vspace{0.15ex}
Стереографический цикл (предельный стереографический цикл) системы (D), 
проходящий через бесконечно удаленную точку $M^{}_{\infty}$
\vspace{0.35ex}
расширенной фазовой плоскости $\overline{Oxy},$ назовем {\it разомкнутым}.
\vspace{0.5ex}

{\bf Свойство 7.7 (7.8).}
{\it
Образом проходящего через начало координат $O(0,0)$ фазовой плоскости $Oxy$ 
стереографического цикла {\rm(}предельного стереографического цикла{\rm)} системы {\rm(D)} на фазовой плоскости $O^{\ast}uv$ 
является разомкнутый стереографический цикл {\rm(}разомкнутый  предельный стереографический цикл{\rm)}
системы} (3.3).
\vspace{0.5ex}

{\bf Свойство 7.9 (7.10).}
{\it
Образом проходящего через начало координат $O(0,0)$ фазовой плоскости $Oxy$ 
цикла {\rm(}предельного цикла{\rm)} системы {\rm(D)} на фазовой плоскости $O^{\ast}uv$ 
является разомкнутый стереографический цикл {\rm(}разомкнутый  предельный стереографический цикл{\rm)}
системы} (3.3).
\vspace{0.5ex}

Так, прямая-траектория $u+v-2=0$ является разокнутым стереографическим циклом 
дифференциальной системы (6.4).
\vspace{0.75ex}

{\bf  Пример 7.1.}
Траекториями системы [8, с. 88]
\\[2ex]
\mbox{}\hfill        
$
\dfrac{dx}{dt}=x(x^2+y^2-1)-y(x^2+y^2+1),
\qquad
\dfrac{dy}{dt}=x(x^2+y^2+1)+y(x^2+y^2-1)
$
\hfill (7.1)
\\[2ex]
являются [9; 1] кривые, заданные уравнением
\\[2ex]
\mbox{}\hfill
$
\dfrac{x^2+y^2}{(x^2+y^2-1)^2}\,\exp \Bigl(2\arctg \dfrac{y}{x}\Bigr)=C,
\quad
0\leq C\leq {}+\infty,
\hfill
$
\\[2ex]
среди которых состояние равновесия $O(0,0)$ (устойчивый фокус) 
и неустойчивый предельный цикл $x^2+y^2=1.$
\\[4ex]
\parbox{50mm}{
\mbox{}\hfill
{\unitlength=1mm
\begin{picture}(45,45)
\put(0,0){\includegraphics[width=45mm,height=45mm]{r07-01.eps}}
\end{picture}}
\hfill\mbox{}
\\[2ex]
\mbox{}\hfill
Рис. 7.1
\hfill\mbox{}
}
\parbox{104mm}{
\mbox{}\hfill
{\unitlength=1mm
\begin{picture}(45,45)
\put(0,0){\includegraphics[width=45mm,height=45mm]{r07-02a.eps}}
\put(44,20.5){\makebox(0,0)[cc]{\scriptsize $x$}}
\put(20.5,44){\makebox(0,0)[cc]{\scriptsize $y$}}
\end{picture}}
\quad
{\unitlength=1mm
\begin{picture}(45,45)
\put(0,0){\includegraphics[width=45mm,height=45mm]{r07-02b.eps}}
\put(44,20.5){\makebox(0,0)[cc]{\scriptsize $u$}}
\put(20.5,44){\makebox(0,0)[cc]{\scriptsize $v$}}
\end{picture}}
\hfill\mbox{}
\\[2ex]
\mbox{}\hfill
Рис. 7.2
\hfill\mbox{}
}
\\[4.25ex]
\indent
Траекториями стереографически сопряженной системы 
\\[2ex]
\mbox{}\hfill        
$
\dfrac{du}{d\theta}=u(u^2+v^2-16)-v(u^2+v^2+16),
\quad\
\dfrac{dv}{d\theta}=u(u^2+v^2+16)+v(u^2+v^2-16),
$
\hfill (7.2)
\\[3ex]
где $(u^2+v^2)\,d\theta=dt,$
являются кривые, заданные уравнением
\\[2ex]
\mbox{}\hfill
$
\dfrac{u^2+v^2}{(u^2+v^2-16)^2}\,\exp \Bigl(2\arctg \dfrac{v}{u}\Bigr)=C^\ast,
\quad
0\leq C^\ast\leq {}+\infty,
\quad 
16C^\ast=C,
\hfill
$
\\[2ex]
среди которых состояние равновесия $O^\ast(0,0)$ (устойчивый фокус) 
\vspace{0.25ex}
и не\-устойчи\-вый предельный цикл $u^2+v^2=16.$
\vspace{0.25ex}

На рис. 7.1 
\vspace{0.1ex}
изображены траектории  на сфере (1.1) стереографически 
взаимосопряженных дифференциальных систем (7.1) и (7.2).
\vspace{0.1ex}
Круги, образующие стереографические атласы траекторий систем (7.1) и (7.2),
построены на рис. 7.2.
\vspace{0.75ex}

{\bf  Пример 7.2.}
Траекториями системы Дарбу [10]
\\[2ex]
\mbox{}\hfill                        
$
\dfrac{dx}{dt}={}-y-x(x^2+y^2-1),
\qquad
\dfrac{dy}{dt}=x-y(x^2+y^2-1)
$
\hfill (7.3)
\\[2ex]
являются [9; 1] кривые, заданные уравнением 
\\[2ex]
\mbox{}\hfill
$
\dfrac{x^2+y^2}{1-x^2-y^2}\, \exp\Bigl({}-2\arctg\dfrac{y}{x}\Bigr)=C,
\quad 
{}-\infty\leq C\leq {}+\infty,
\hfill
$
\\[2ex]
среди которых состояние равновесия $O(0,0)$
\vspace{0.25ex}
(неустойчивый фокус) и устойчивый предельный цикл $x^2+y^2=1.$
\vspace{0.25ex}

Траекториями стереографически сопряженной системы 
\\[2ex]
\mbox{}\hfill                        
$
\dfrac{du}{d\theta}=
16u-u^3-u^2v-uv^2-v^3,   
\qquad
\dfrac{dv}{d\theta}=16v+u^3-u^2v+uv^2-v^3,
$
\hfill (7.4)
\\[2ex]
где $(u^2+v^2)\;\!d\theta=dt,$ являются кривые, заданные уравнением 
\\[2ex]
\mbox{}\hfill
$
\dfrac{1}{u^2+v^2-16}\, \exp\Bigl({}-2\arctg\dfrac{v}{u}\Bigr)=C^\ast,
\quad
{}-\infty  \leq C^\ast \leq {}+\infty,
\quad
16C^\ast=C,
\hfill
$
\\[2ex]
среди которых состояние равновесия $O^{\ast}(0,0)$ 
\vspace{0.25ex}
(неустойчивый дикритический узел) и устойчивый предельный цикл $u^2+v^2=16.$
\vspace{0.25ex}

На рис. 7.3
\vspace{0.1ex}
изображены траектории  на сфере (1.1) стереографически 
взаимосопряженных дифференциальных систем (7.3) и (7.4).
\vspace{0.1ex}
Круги, образующие стереографические атласы траекторий систем (7.3) и (7.4),
построены на рис. 7.4.
\\[4ex]
\parbox{50mm}{
\mbox{}\hfill
{\unitlength=1mm
\begin{picture}(45,45)
\put(0,0){\includegraphics[width=45mm,height=45mm]{r07-03.eps}}
\end{picture}}
\hfill\mbox{}
\\[2ex]
\mbox{}\hfill
Рис. 7.3
\hfill\mbox{}
}
\parbox{104mm}{
\mbox{}\hfill
{\unitlength=1mm
\begin{picture}(45,45)
\put(0,0){\includegraphics[width=45mm,height=45mm]{r07-04a.eps}}
\put(44,20.5){\makebox(0,0)[cc]{\scriptsize $x$}}
\put(20.5,44){\makebox(0,0)[cc]{\scriptsize $y$}}
\end{picture}}
\quad
{\unitlength=1mm
\begin{picture}(45,45)
\put(0,0){\includegraphics[width=45mm,height=45mm]{r07-04b.eps}}
\put(44,20.5){\makebox(0,0)[cc]{\scriptsize $u$}}
\put(20.5,44){\makebox(0,0)[cc]{\scriptsize $v$}}
\end{picture}}
\hfill\mbox{}
\\[2ex]
\mbox{}\hfill
Рис. 7.4
\hfill\mbox{}
}
\\[6.25ex]
\indent
{\bf  Пример 7.3.}
Заменой $x$ на $x-1,\ y$ на $y$ систему (7.3) приводим к системе Дарбу 
\\[2.2ex]
\mbox{}\hfill                        
$
\dfrac{dx}{dt}={}-2x-y+3x^2+y^2-x(x^2+y^2),
\qquad
\dfrac{dy}{dt}={}-1+x+2xy-y(x^2+y^2).
$
\hfill (7.5)
\\[-3.25ex]

\newpage

Кривые семейства 
\\[2ex]
\mbox{}\hfill
$
\dfrac{(x-1)^2+y^2}{1-(x-1)^2-y^2}\, \exp\Bigl({}-2\arctg\dfrac{y}{x-1}\Bigr)=C,
\quad 
{}-\infty\leq C\leq {}+\infty,
\hfill
$
\\[2ex]
являются траекториями системы (7.5). При этом состояние равновесия 
\vspace{0.25ex}
$A(1,0)$ --- неустойчивый фокус, а окружность $(x-1)^2+y^2=1,$ 
\vspace{0.25ex}
проходящая через начало координат $O(0,0),$ суть устойчивый предельный цикл системы (7.5).
\vspace{0.25ex}

Траекториями стереографически сопряженной системы 
\\[2ex]
\mbox{}\hfill                        
$
\dfrac{du}{d\theta}=
16u-12u^2+4v^2+2u^3-u^2v-2uv^2-v^3+\dfrac{1}{2}\,u^3v+\dfrac{1}{2}\,uv^3,   
\hfill
$
\\[0.5ex]
\mbox{}\hfill (7.6)
\\[0.5ex]
\mbox{}\hfill
$
\dfrac{dv}{d\theta}=
16v-16uv+u^3+4u^2v+uv^2-\dfrac{1}{4}\,u^4+\dfrac{1}{4}\,v^4,   
\hfill
$
\\[2.25ex]
где $(u^2+v^2)\;\!d\theta=dt,$
являются кривые семейства
\\[2ex]
\mbox{}\hfill
$
\dfrac{(u-4)^2+v^2}{u-2}\, \exp\Bigl(2\;\!\arctg\dfrac{4v}{(u-2)^2+v^2-4}\Bigr)=C^\ast,
\quad 
{}-\infty\leq C^\ast\leq {}+\infty,
\quad
 C^\ast=8C,
\hfill
$
\\[2ex]
среди которых прямая $\!u\!=\!2,\!$ 
\vspace{0.15ex}
являющаяся разомкнутым предельным стереографичес\-ким циклом, 
неустойчивый дикритический узел $O^{\ast}(0,0)$ 
\vspace{0.25ex}
и  неустойчивый фокус  $A^{\ast}(4,0).$

Круги, образующие стереографические атласы траекторий 
\vspace{0.1ex}
дифференциальных сис\-тем (7.5) и (7.6),
построены на рис. 7.5.
\\[3.5ex]
\mbox{}\hfill
{\unitlength=1mm
\begin{picture}(45,45)
\put(0,0){\includegraphics[width=45mm,height=45mm]{r07-05a.eps}}
\put(44,20.5){\makebox(0,0)[cc]{\scriptsize $x$}}
\put(20.5,44){\makebox(0,0)[cc]{\scriptsize $y$}}
\end{picture}}
\qquad\qquad
{\unitlength=1mm
\begin{picture}(45,45)
\put(0,0){\includegraphics[width=45mm,height=45mm]{r07-05b.eps}}
\put(44,20.5){\makebox(0,0)[cc]{\scriptsize $u$}}
\put(20.5,44){\makebox(0,0)[cc]{\scriptsize $v$}}
\end{picture}}
\hfill\mbox{}
\\[1.5ex]
\mbox{}\hfill
Рис. 7.5
\hfill\mbox{}
\\[4ex]
\indent
{\bf  Пример 7.4.}
Заменой $u$ на $u+2,\ v$ на $v$ систему (7.6) приводим к системе 
\\[2ex]
\mbox{}\hfill                        
$
\dfrac{du}{d\theta}=
{}-8u+2uv+2u^3+2u^2v-2uv^2+\dfrac{1}{2}\,u^3v+\dfrac{1}{2}\,uv^3,   
\hfill
$
\\[0.5ex]
\mbox{}\hfill (7.7)
\\[0.5ex]
\mbox{}\hfill
$
\dfrac{dv}{d\theta}=
4+4u+2v^2-u^3+4u^2v+uv^2-\dfrac{1}{4}\,u^4+\dfrac{1}{4}\,v^4.   
\hfill
$
\\[2ex]
\indent
Кривые семейства 
\\[2ex]
\mbox{}\hfill
$
\dfrac{(u-2)^2+v^2}{u}\, \exp\Bigl(2\arctg\dfrac{4v}{u^2+v^2-4}\Bigr)=C^{\ast},
\quad 
{}-\infty\leq C^{\ast}\leq {}+\infty,
\hfill
$
\\[2ex]
являются траекториями системы (7.7), 
\vspace{0.25ex}
среди которых прямая $u=0,$ являющаяся разомкнутым предельным стереографическим циклом, 
\vspace{0.25ex}
неустойчивый дикритический узел $A^{}_1({}-2,0)$ и  неустойчивый фокус  $A^{}_2(2,0).$

\newpage

Траекториями стереографически сопряженной системы 
\\[2ex]
\mbox{}\hfill                        
$
\dfrac{dx}{dt}=
{}-32x-8xy+8x^3-8x^2y-8xy^2-2x^3y-2xy^3,   
\hfill
$
\\[0.5ex]
\mbox{}\hfill (7.8)
\\[0.5ex]
\mbox{}\hfill
$
\dfrac{dy}{dt}=
{}-16-16x-8y^2+4x^3+16x^2y-4xy^2+x^4-y^4,   
\hfill
$
\\[2.25ex]
где $(x^2+y^2)\;\!dt=d\theta,$
являются кривые семейства
\\[2ex]
\mbox{}\hfill
$
\dfrac{(x-2)^2+y^2}{x}\, \exp\Bigl(2\arctg\dfrac{4y}{4-x^2-y^2}\Bigr)=C^\ast,
\quad
{}-\infty\leq C^\ast\leq {}+\infty,
\hfill
$
\\[2ex]
среди которых прямая $\!x\!=\!0,\!$ 
\vspace{0.25ex}
являющаяся разомкнутым предельным стереографическим циклом, 
неустойчивый дикритический узел $A^{\ast}_1({}-2,0)$ 
\vspace{0.25ex}
и  неустойчивый фокус $A^{\ast}_2(2,0).$

Круги, образующие стереографические атласы траекторий 
\vspace{0.1ex}
дифференциальных систем (7.7) и (7.8),
построены на рис. 7.6.
\\[3.25ex]
\mbox{}\hfill
{\unitlength=1mm
\begin{picture}(45,45)
\put(0,0){\includegraphics[width=45mm,height=45mm]{r07-06a.eps}}
\put(44,20.5){\makebox(0,0)[cc]{\scriptsize $x$}}
\put(20.5,44){\makebox(0,0)[cc]{\scriptsize $y$}}
\end{picture}}
\qquad\qquad
{\unitlength=1mm
\begin{picture}(45,45)
\put(0,0){\includegraphics[width=45mm,height=45mm]{r07-06b.eps}}
\put(44,20.5){\makebox(0,0)[cc]{\scriptsize $u$}}
\put(20.5,44){\makebox(0,0)[cc]{\scriptsize $v$}}
\end{picture}}
\hfill\mbox{}
\\[1ex]
\mbox{}\hfill
Рис. 7.6
\hfill\mbox{}
\\[4.25ex]
\centerline{
{\bf  8. Симметpичность фазового поля направлений}}
\\
\centerline{{\bf стереографически сопряженых дифференциальных систем}}
\\[1.5ex]
\indent
На основании аналитических условий симметричности фазового поля направлений 
дифференциальной системы [11; 1] получаем критерии симметричности 
для стереографически сопряженых дифференциальных систем.
\vspace{0.5ex}

{\bf Свойство 8.1.}
{\it 
Равносильными являются следующие утверждения\;\!{\rm:}
\vspace{0.15ex}

{\rm 1.} 
\vspace{0.15ex}
Фазовое поле направлений системы {\rm (D)} симметpично
относительно начала кооpдинат фазовой плоскости $Oxy;$ 
\vspace{0.25ex}

{\rm 2.} 
Выполняется тождество
\\[1.5ex]
\mbox{}\hfill
$
X(x,y)\;\!Y({}-x,{}-y)  - X({}-x,{}-y)\;\!Y(x,y)  = 0
\quad
\forall (x,y)\in \R^2;
\hfill
$
\\[1.5ex]
\indent
{\rm 3.} 
\vspace{0.15ex}
Фазовое поле направлений системы {\rm (3.3)} симметpично
относительно начала кооpдинат фазовой плоскости $O^{\ast}uv;$ 
\vspace{0.25ex}

{\rm 4.} 
Выполняется тождество}
\\[1.5ex]
\mbox{}\hfill
$
U(u,v)\;\!V({}-u,{}-v)  - U({}-u,{}-v)\;\!V(u,v)  = 0
\quad
\forall (u,v)\in \R^2.
\hfill
$
\\[1.75ex]
\indent
Например, такой симметричностью обладают фазовые поля направлений 
стереографически сопряженных систем 
\vspace{0.25ex}
(4.4) и (4.5) при $a_{_0}=b_{_0},$ 
(5.1) и (5.2), 
(5.3) и (5.4), 
(5.5) и (5.6), 
(5.7) и (5.8), 
(5.9) и (5.10),
(5.11) и (5.12),
(7.1) и (7.2), 
(7.3) и (7.4).
\vspace{1ex}

{\bf Свойство 8.2.}
{\it 
Равносильными являются следующие утверждения\;\!{\rm:}
\vspace{0.15ex}

{\rm 1.} 
\vspace{0.15ex}
Фазовое поле направлений дифференциальной системы {\rm (D)} симметpично
относительно кооpдинатной оси $Ox;$

{\rm 2.} 
Выполняется тождество
\\[1.5ex]
\mbox{}\hfill
$
X(x,y)\;\!Y(x,{}-y)  + X(x,{}-y)\;\!Y(x,y)  = 0
\quad
\forall (x,y)\in \R^2;
\hfill
$
\\[1.5ex]
\indent
{\rm 3.} 
\vspace{0.15ex}
Фазовое поле направлений дифференциальной системы {\rm (3.3)} симметpично
относительно кооpдинатной оси\,$O^{\ast}\;\!\!u;$ 
\vspace{0.25ex}

{\rm 4.} 
Выполняется тождество}
\\[1.5ex]
\mbox{}\hfill
$
U(u,v)\;\!V(u,{}-v)  + U(u,{}-v)\;\!V(u,v)  = 0
\quad
\forall (u,v)\in \R^2.
\hfill
$
\\[1.75ex]
\indent
Например, такой симметричностью обладают фазовые поля направлений 
стереографически сопряженных дифференциальных систем 
\vspace{0.25ex}
(4.4) и (4.5) при $a_{_0}b_{_0}=0,$ 
(5.1) и (5.2), 
(5.3) и (5.4), 
(5.7) и (5.8), 
(5.9) и (5.10).
\vspace{1ex}

{\bf Свойство 8.3.}
\vspace{0.15ex}
{\it 
Равносильными являются следующие утверждения\;\!{\rm:}

{\rm 1.} 
\vspace{0.15ex}
Фазовое поле направлений дифференциальной системы {\rm (D)} симметpично
относительно кооpдинатной оси $Oy;$
\vspace{0.25ex}

{\rm 2.} 
Выполняется тождество
\\[1.5ex]
\mbox{}\hfill
$
X(x,y)\;\!Y({}-x,y)  + X({}-x,y)\;\!Y(x,y)  = 0
\quad
\forall (x,y)\in \R^2;
\hfill
$
\\[1.5ex]
\indent
{\rm 3.} 
\vspace{0.15ex}
Фазовое поле направлений дифференциальной системы {\rm (3.3)} симметpично
относительно кооpдинатной оси\,$O^{\ast}v;$ 
\vspace{0.25ex}

{\rm 4.} 
Выполняется тождество}
\\[1.5ex]
\mbox{}\hfill
$
U(u,v)\;\!V({}-u,v)  + U({}-u,v)\;\!V(u,v)  = 0
\quad
\forall (u,v)\in \R^2.
\hfill
$
\\[1.75ex]
\indent
Например, такой симметричностью обладают фазовые поля направлений 
стереографически сопряженных дифференциальных систем 
\vspace{0.25ex}
(4.4) и (4.5) при $a_{_0}b_{_0}=0,$ 
(5.1) и (5.2), 
(5.3) и (5.4), 
(5.7) и (5.8), 
(5.9) и (5.10).
\vspace{1ex}

{\bf Свойство 8.4.}
\vspace{0.15ex}
{\it 
Равносильными являются следующие утверждения\;\!{\rm:}

{\rm 1.} 
\vspace{0.15ex}
Фазовое поле направлений дифференциальной системы {\rm (D)} симметpично
относительно пpямой $y = x;$
\vspace{0.25ex}

{\rm 2.} 
Выполняется тождество
\\[1.5ex]
\mbox{}\hfill
$
X(x,y)\;\!X(y,x)  -  Y(x,y)\;\!Y(y,x)  =  0
\quad 
\forall (x,y)\in \R^2;
\hfill
$
\\[1.5ex]
\indent
{\rm 3.} 
\vspace{0.15ex}
Фазовое поле направлений дифференциальной системы {\rm (3.3)} симметpично
относительно пpямой $v = u;$
\vspace{0.25ex}

{\rm 4.} 
Выполняется тождество}
\\[1.5ex]
\mbox{}\hfill
$
U(u,v)\;\!U(v,u)  -  V(u,v)\;\!V(v,u)  =  0
\quad 
\forall (u,v)\in \R^2.
\hfill
$
\\[1.75ex]
\indent
Например, такой симметричностью обладают фазовые поля направлений 
стереографически сопряженных дифференциальных систем 
\vspace{0.25ex}
(4.4) и (4.5) при $|a_{_0}|=|b_{_0}|,$ 
(5.1) и (5.2), 
(5.3) и (5.4), 
(5.7) и (5.8), 
(6.1) и (6.2), 
(6.3) и (6.4).
\vspace{1ex}

{\bf Свойство 8.5.}
\vspace{0.15ex}
{\it 
Равносильными являются следующие утверждения\;\!{\rm:}

{\rm 1.} 
\vspace{0.15ex}
Фазовое поле направлений дифференциальной системы {\rm (D)} симметpично
относительно пpямой $y = {}-x;$
\vspace{0.25ex}

{\rm 2.} 
Выполняется тождество
\\[1.5ex]
\mbox{}\hfill
$
X({}-x,{}-y)\;\!X(y,x)  -  Y({}-x,{}-y)\;\!Y(y,x)  =  0
\quad 
\forall (x,y)\in \R^2;
\hfill
$
\\[1.5ex]
\indent
{\rm 3.} 
\vspace{0.15ex}
Фазовое поле направлений дифференциальной системы {\rm (3.3)} симметpично
относительно пpямой $v = {}-u;$
\vspace{0.25ex}

{\rm 4.} 
Выполняется тождество}
\\[1.5ex]
\mbox{}\hfill
$
U({}-u,{}-v)\;\!U(v,u)  -  V({}-u,{}-v)\;\!V(v,u)  =  0
\quad 
\forall (u,v)\in \R^2.
\hfill
$
\\[1.75ex]
\indent
Например, такой симметричностью обладают фазовые поля направлений 
стереографически сопряженных дифференциальных систем 
\vspace{0.25ex}
(4.4) и (4.5) при $|a_{_0}|=|b_{_0}|,$ 
(5.1) и (5.2), 
(5.3) и (5.4), 
(5.7) и (5.8).
\\[3.5ex]
\centerline{
{\bf  9. Бесконечно удаленное состояние равновесия
}
}
\\[1.5ex]
\indent
В соответствии со свойством 6.2 бесконечно удаленная точка $M_\infty^{}$ 
\vspace{0.25ex}
расширенной фазовой плоскости $\overline{Oxy}$ является состоянием равновесия системы (D), 
\vspace{0.15ex}
если и только если точка $O^\ast(0,0)$
\vspace{0.25ex}
является состоянием равновесия системы (3.3), причем у 
состояний равновесия $M_\infty^{}$ и $O^\ast$ один и тот же вид.
\vspace{0.15ex}

Поведение траекторий системы (D) 
\vspace{0.25ex}
в окрестности бесконечно удаленной точки $M_\infty^{}$ 
расширенной фазовой плоскости $\overline{Oxy}$ 
\vspace{0.25ex}
определяется поведением ее траекторий в  окрестности бесконечно удаленной прямой 
\vspace{0.25ex}
проективной фазовой плоскости $\P\R(x,y)$ [12; 1].

Пусть $L$ --- бесконечно удаленное состояние равновесия системы (D) на  
\vspace{0.15ex}
проективной фазовой плоскости $\P\R(x,y).$ 
\vspace{0.15ex}
Секторы Бендиксона состояния равновесия $L$ 
разделим на два вида: внутренние и внешние. 
\vspace{0.15ex}
Точки бесконечно удаленной прямой проективной фазовой плоскости $\P\R(x,y),$ 
лежащие в проколотой окрестности состояния равновесия $L,$ принадлежат внешним секторам Бендиксона 
и не принадлежат внутренним секторам Бендиксона. При этом имеет место
\vspace{0.5ex}

{\bf Свойство 9.1.}
{\it
Внутреннему сектору Бендиксона состояния равновесия $L$
соответствует сектор Бендиксона состояния равновесия $M_\infty^{}$
такого же вида с таким же направлением движения вдоль траекторий.
}
\vspace{0.5ex}

Внешние секторы Бендиксона бесконечно удаленных состояний равновесия
проективной фазовой плоскости $\P\R(x,y)$ стягиваются, образуя секторы Бендиксона
\vspace{0.25ex}
бесконечно удаленного состояния равновесия $M_\infty^{}$
расширенной фазовой плоскости $\overline{Oxy}\,.$ 
\vspace{0.5ex}
При этом существенное значение имеет то, что [12; 1]
у проективно неособой системы (D) 
\vspace{0.25ex}
бесконечно удаленная прямая проективной фазовой плоскости $\P\R(x,y)$
состоит из траекторий. Это позволяет указать следующие закономерности. 
\vspace{0.25ex}

Пусть $L_1^{}$ и $L_2^{}$ --- соседние 
\vspace{0.25ex}
бесконечно удаленные состояния равновесия 
на окружности проективного круга  $\P\K(x,y)$ [12; 1] 
системы (D). 
\vspace{0.75ex}

{\bf Свойство 9.2.}
\vspace{0.15ex}
{\it
Если у проективно неособой системы {\rm(D)} смежные внешние секторы Бендиксона состояний равновесия 
$L_1^{}$ и $L_2^{}\colon$ 
\vspace{0.25ex}
а{\rm)} параболические{\rm;} 
б{\rm)} гиперболические{\rm;} 
\linebreak
в{\rm)} один --- параболический, другой --- гиперболический{\rm;} 
\vspace{0.15ex}
г{\rm)} один --- гиперболический, другой --- эллиптический, 
то 
они стягиваются в{\rm:} 
\vspace{0.25ex}
а{\rm)} эллиптический{\rm;} 
б{\rm)} гиперболический{\rm;} 
\linebreak
в{\rm)} параболический{\rm;} 
г{\rm)} эллиптический  
сектор Бендиксона состояния равновесия $M_\infty^{}.$
}
\\[3.75ex]
\mbox{}\hfill
{\unitlength=1mm
\begin{picture}(45,45)
\put(0,0){\includegraphics[width=45mm,height=45mm]{r09-01a.eps}}
\put(44,20.5){\makebox(0,0)[cc]{\scriptsize $x$}}
\put(20.5,44){\makebox(0,0)[cc]{\scriptsize $y$}}
\end{picture}}
\qquad\qquad
{\unitlength=1mm
\begin{picture}(45,45)
\put(0,0){\includegraphics[width=45mm,height=45mm]{r09-01b.eps}}
\put(44,20.5){\makebox(0,0)[cc]{\scriptsize $u$}}
\put(20.5,44){\makebox(0,0)[cc]{\scriptsize $v$}}
\end{picture}}
\hfill\mbox{}
\\[0ex]
\mbox{}\hfill
Рис. 9.1
\hfill\mbox{}
\\[-4ex]

\newpage

{\bf  Пример 9.1.}
На  проективной фазовой плоскости $\P\R(x,y)$ у дифференциальной системы [5, c. 84 --- 85; 1 ]
\\[2ex]
\mbox{}\hfill        
$
\dfrac{dx}{dt}=1-x^2-y^2,
\qquad 
\dfrac{dy}{dt}=xy-1
$
\hfill(9.1)
\\[2.35ex]
одно состояние равновесия --- узел, которое лежит на <<концах>> оси $Ox.$
\vspace{0.25ex}

Состояние равновесия $O^\ast(0,0)$ стереографически сопряженной системы
\\[2ex]
\mbox{}\hfill        
$
\dfrac{du}{d\theta}=4u^4-8u^2\;\!v^2-4v^4-\dfrac{1}{4}\,u^6+\dfrac{1}{2}\,u^5\;\!v-
\dfrac{1}{4}\,u^4\;\!v^2+u^3\;\!v^3+\dfrac{1}{4}\, u^2\;\!v^4+\dfrac{1}{2}\,uv^5+\dfrac{1}{4}\, v^6,
\hfill
$
\\[0.75ex]
\mbox{}\hfill (9.2)
\\[0.5ex]
\mbox{}\hfill
$
\dfrac{dv}{d\theta}=
12u^3\;\!v+4uv^3-\dfrac{1}{4}\,u^6-\dfrac{1}{2}\,u^5\;\!v-\dfrac{1}{4}\,u^4\;\!v^2-
u^3\;\!v^3+\dfrac{1}{4}\, u^2\;\!v^4-\dfrac{1}{2}\,uv^5+\dfrac{1}{4}\, v^6,
\hfill
$
\\[2.35ex]
где $(u^2+v^2)^2\;\!d\theta=dt,$
\vspace{0.25ex}
состоит из двух эллиптических секторов, разделенных двумя параболическими 
секторами (свойство 9.2, случай а).
\vspace{0.25ex}
Учитывая качественное исследование [1] поведения траекторий системы (9.1),
на рис. 9.1 построены круги, образующие стереографические атласы траекторий систем (9.1) и (9.2).
\vspace{1ex}

{\bf  Пример 9.2.}
Дифференциальная система [13, с. 61 --- 65]
\\[2ex]
\mbox{}\hfill        
$
2\,\dfrac{dx}{dt}=2y+i\;\!(x-i\;\!y)^q-i\;\!(x+i\;\!y)^q,
\qquad 
2\,\dfrac{dy}{dt}={}-2x+(x-i\;\!y)^q+ (x+i\;\!y)^q
$
\hfill (9.3)
\\[2.25ex]
при $i=\sqrt{{}-1}\,,\, q=4$ и $q=5$  
\vspace{0.35ex}
является проективно неособой [12], 
а все бесконечно удаленные состояния равновесия на проективной фазовой плоскости $\P\R(x,y)$ 
\vspace{0.35ex}
являются узлами (см.  рис. 2.12 из [13, с. 65]  или рис. 8.9 и 8.10 из [12]).
\vspace{0.5ex}

При $q=4$  система (9.3) имеет  вид
\\[2ex]
\mbox{}\hfill        
$
\dfrac{dx}{dt}=y+4x^3y-4xy^3,
\qquad 
\dfrac{dy}{dt}={}-x+x^4-6x^2y^2+y^4,
$
\hfill (9.4)
\\[2.5ex]
а стереографически сопряженной к ней является 
система
\\[2.25ex]
\mbox{}\qquad\qquad        
$
\dfrac{du}{d\theta}=
v({}-384u^5+1280u^3v^2-384uv^4+
u^8+4u^6v^2+6u^4v^4+4u^2v^6+v^8),
\hfill
$
\\[0.5ex]
\mbox{}\hfill (9.5)
\\[0.5ex]
\mbox{}\qquad\qquad
$
\dfrac{dv}{d\theta}=
64u^6-960u^4v^2+960u^2v^4-64v^6
-u^9-4u^7v^2-6u^5v^4-4u^3v^6-uv^8,
\hfill
$
\\[2.5ex]
где $(u^2+v^2)^4\;\!d\theta=dt.$
\vspace{0.5ex}

При $q=5$ система (9.3) имеет  вид
\\[2ex]
\mbox{}\hfill        
$
\dfrac{dx}{dt}=y+5x^4\;\!y-10x^2\;\!y^3+y^5,
\qquad 
\dfrac{dy}{dt}={}-x+x^5-10x^3\;\!y^2+5xy^4,
$
\hfill (9.6)
\\[2.5ex]
а стереографически сопряженной к ней является система 
\\[2.25ex]
\mbox{}        
$
\dfrac{du}{d\theta}=
v({}-1792u^6+8960u^4v^2-5376u^2v^4+256v^6+
u^{10}+5u^8v^2+10u^6v^4+10u^4v^6+5u^2v^8+v^{10}),
\ \ 
\hfill
$
\\[0.5ex]
\mbox{}\hfill (9.7)
\\[0.5ex]
\mbox{}
$
\dfrac{dv}{d\theta}=
u(256u^6-5376u^4v^2+8960u^2v^4-1792v^6
-u^{10}-5u^8v^2-10u^6v^4-10u^4v^6-5u^2v^8-v^{10}),
\ \ 
\hfill
$
\\[2.5ex]
где $(u^2+v^2)^5\;\!d\theta=dt.$
\vspace{0.5ex}

Согласно свойству 9.2 состояние равновесия $O^\ast(0,0)$ 
\vspace{0.25ex}
как дифференциальной системы (9.5), так и дифференциальной системы (9.7)
\vspace{0.25ex}
состоит из эллиптических секторов Бендиксона.
Круги, образующие стереографические атласы траекторий дифференциальных систем (9.4) и (9.5),
построены на рис. 9.2,
а круги, образующие стереографические атласы траекторий дифференциальных систем (9.6) и (9.7),
построены на рис. 9.3.
\\[3.75ex]
\mbox{}\hfill
{\unitlength=1mm
\begin{picture}(45,45)
\put(0,0){\includegraphics[width=45mm,height=45mm]{r09-02a.eps}}
\put(44,20.5){\makebox(0,0)[cc]{\scriptsize $x$}}
\put(20.5,44){\makebox(0,0)[cc]{\scriptsize $y$}}
\end{picture}}
\qquad\qquad
{\unitlength=1mm
\begin{picture}(45,45)
\put(0,0){\includegraphics[width=45mm,height=45mm]{r09-02b.eps}}
\put(44,20.5){\makebox(0,0)[cc]{\scriptsize $u$}}
\put(20.5,44){\makebox(0,0)[cc]{\scriptsize $v$}}
\end{picture}}
\hfill\mbox{}
\\[1ex]
\mbox{}\hfill
Рис. 9.2
\hfill\mbox{}
\\[6.75ex]
\mbox{}\hfill
{\unitlength=1mm
\begin{picture}(45,45)
\put(0,0){\includegraphics[width=45mm,height=45mm]{r09-03a.eps}}
\put(44,20.5){\makebox(0,0)[cc]{\scriptsize $x$}}
\put(20.5,44){\makebox(0,0)[cc]{\scriptsize $y$}}
\end{picture}}
\qquad\qquad
{\unitlength=1mm
\begin{picture}(45,45)
\put(0,0){\includegraphics[width=45mm,height=45mm]{r09-03b.eps}}
\put(44,20.5){\makebox(0,0)[cc]{\scriptsize $u$}}
\put(20.5,44){\makebox(0,0)[cc]{\scriptsize $v$}}
\end{picture}}
\hfill\mbox{}
\\[2ex]
\mbox{}\hfill
Рис. 9.3
\hfill\mbox{}
\\[4ex]
\indent
{\bf  Пример 9.3.}
\vspace{0.25ex}
На  проективной фазовой плоскости $\P\R(x,y)$ у 
дифференциальной системы [5, c. 85 --- 87; 14, c. 209 --- 212]
\\[2.25ex]
\mbox{}\hfill        
$
\dfrac{dx}{dt}={}-1+x^2+y^2,
\qquad 
\dfrac{dy}{dt}={}-5+5xy
$
\hfill(9.8)
\\[2.5ex]
три состояния равновесия:
\vspace{0.25ex}
седло,  лежащее на <<концах>> прямой $y=0,$ 
и два устойчивых узла,
одно из которых  лежит на <<концах>> прямой $y={}-2x,$
\vspace{0.25ex}
а другое  лежит на <<концах>> прямой $y=2x.$
\vspace{0.35ex}

У стереографически сопряженной системы
\\[2.75ex]
\mbox{}\hfill        
$
\dfrac{du}{d\theta}={}-4u^4-40u^2v^2+4v^4+\dfrac{1}{4}\,u^6+\dfrac{5}{2}\,u^5v+\dfrac{1}{4}\,u^4v^2+
5u^3v^3-\dfrac{1}{4}\, u^2v^4+\dfrac{5}{2}\,uv^5-\dfrac{1}{4}\, v^6,
\hfill
$
\\[0.5ex]
\mbox{}\hfill (9.9)
\\[0.5ex]
\mbox{}\hfill
$
\dfrac{dv}{d\theta}=
12u^3v-28uv^3-\dfrac{5}{4}\,u^6+\dfrac{1}{2}\,u^5v-\dfrac{5}{4}\,u^4v^2+
u^3v^3+\dfrac{5}{4}\, u^2v^4+\dfrac{1}{2}\,uv^5+\dfrac{5}{4}\, v^6,
\hfill
$
\\[2.75ex]
где $(u^2+v^2)^2\;\!d\theta=dt,$
\vspace{0.25ex}
состояние равновесия $O^\ast(0,0),$ 
состоящее из двух эллиптических и двух параболических 
секторов Бендиксона (свойство 9.2, случаи а и в).
\vspace{0.25ex}
Учитывая качественное исследование [1]
поведения траекторий системы (9.8),
\vspace{0.25ex}
на рис. 9.4 построены круги, образующие стереографические атласы траекторий систем (9.8) и (9.9).

\newpage

\mbox{}
\\[-1ex]
\mbox{}\hfill
{\unitlength=1mm
\begin{picture}(45,45)
\put(0,0){\includegraphics[width=45mm,height=45mm]{r09-04a.eps}}
\put(44,20.5){\makebox(0,0)[cc]{\scriptsize $x$}}
\put(20.5,44){\makebox(0,0)[cc]{\scriptsize $y$}}
\end{picture}}
\qquad\qquad
{\unitlength=1mm
\begin{picture}(45,45)
\put(0,0){\includegraphics[width=45mm,height=45mm]{r09-04b.eps}}
\put(44,20.5){\makebox(0,0)[cc]{\scriptsize $u$}}
\put(20.5,44){\makebox(0,0)[cc]{\scriptsize $v$}}
\end{picture}}
\hfill\mbox{}
\\[0ex]
\mbox{}\hfill
Рис. 9.4
\hfill\mbox{}
\\[4ex]
\indent
{\bf Свойство 9.3.}
\vspace{0.15ex}
{\it
Если граничная окружность проективного круга $\P\K(x,y)$
пересекается каждой траекторией проективно особой системы {\rm(D)}
\vspace{0.15ex}
ортогонально, то 
бесконечно удаленное состояние равновесия $M_\infty^{}$
является дикритическим узлом.
}
\vspace{0.75ex}

Такой, например,  является дифференциальная система Дарбу (7.3), проективный атлас траекторий которой
построен на рис. 16.3 из [1], а стереографический атлас траекторий построен на рис. 7.4.
\vspace{0.75ex}

{\bf  Пример 9.4.}
\vspace{0.35ex}
В [11; 1] качественно исследовано поведение траекторий на  проективной фазовой плоскости $\P\R(x,y)$ 
и построены проективные атласы траекторий систем
\\[2ex]
\mbox{}\hfill        
$
\dfrac{dx}{dt}={}-y+x^3,
\qquad 
\dfrac{dy}{dt}=x(1+xy),
$
\hfill(9.10)
\\[3.5ex]
\mbox{}\hfill                                       
$
\begin{array}{l}
\dfrac{dx}{dt}=x(x^2+y^2-1)(x^2+y^2-9)-y(x^2+y^2-2x-8),
\\[3ex]
\dfrac{dy}{dt}=y(x^2+y^2-1)(x^2+y^2-9)+x(x^2+y^2-2x-8),
\end{array}
$
\hfill (9.11)
\\[3.5ex]
\mbox{}\hfill                                       
$
\begin{array}{l}
\dfrac{dx}{dt}=x(2x^2+2y^2+1)\Bigl((x^2+y^2)^2+x^2-y^2+\dfrac{1}{10}\Bigr)-y(2x^2+2y^2-1),
\\[3ex]
\dfrac{dy}{dt}=y(2x^2+2y^2-1)\Bigl((x^2+y^2)^2+x^2-y^2+\dfrac{1}{10}\Bigr)+x(2x^2+2y^2+1).
\end{array}
$
\hfill (9.12)
\\[3ex]
\indent
Граничную окружность проективного круга $\P\K(x,y)$
 траектории каждой из этих систем пересекают ортогонально.

Cтереографически сопряженными к  системам (9.10), (9.11), (9.12) соответственно являются системы
\\[1ex]
\mbox{}\hfill        
$
\dfrac{du}{d\theta}={}-(16u^3+u^4v+2u^2v^3+v^5),
\quad 
\dfrac{dv}{d\theta}=u({}-16uv+u^4+2u^2v^2+v^4),
$
\hfill(9.13)
\\[2.35ex]
где $(u^2+v^2)^2\;\!d\theta=dt,$
\\[2.5ex]
\mbox{}\hfill                                       
$
\begin{array}{l}
\dfrac{du}{d\theta}=
{}-256u +160u^3-16u^2v+160uv^2-16v^3
+8u^3v+8uv^3\ -
\\[2.5ex]
\mbox{}\qquad \ \ 
-\ 9u^5 +8u^4v-18u^3v^2+16u^2v^3-9uv^4+8v^5,
\\[3ex]
\dfrac{dv}{d\theta}=
{}-256v +16u^3+160u^2v +16uv^2+160v^3
-8u^4-8u^2v^2\ -
\\[2.5ex]
\mbox{}\qquad \ \
- \ 8u^5-9u^4v-16u^3v^2-18u^2v^3-8uv^4-9v^5,
\end{array}
$
\hfill (9.14)
\\[2.35ex]
где $(u^2+v^2)^2\;\!d\theta=dt,$
\\[2.35ex]
\mbox{}                                       
$
\dfrac{du}{d\theta}=
{}-8192u-768u^3+1280uv^2
-\dfrac{96}{5}\,u^5-32u^4v+\dfrac{288}{5}\,u^3v^2-64u^2v^3 -\
\hfill
$
\\[2.5ex]
\mbox{}\qquad
$
-\ \dfrac{256}{5}\,uv^4-32v^5 - \dfrac{1}{10}\,u^7-3u^6v+
\dfrac{1}{10}\,u^5v^2-5u^4v^3+\dfrac{1}{2}\,u^3v^4-u^2v^5+\dfrac{3}{10}\,uv^6+v^7,
\hfill
$
\\[0.5ex]
\mbox{}\hfill (9.15)
\\[0.5ex]
\mbox{}
$
\dfrac{dv}{d\theta}=
{}-8192v-1280u^2v+768v^3
+32u^5-\dfrac{256}{5}\,u^4v+64u^3v^2+\dfrac{288}{5}\,u^2v^3\ +
\hfill
$
\\[2.5ex]
\mbox{}\qquad
$
+\ 32uv^4-\dfrac{96}{5}\,v^5 +
u^7-\dfrac{3}{10}\,u^6v-u^5v^2-\dfrac{1}{2}\,u^4v^3-5u^3v^4-\dfrac{1}{10}\,u^2v^5
-3uv^6+\dfrac{1}{10}\,v^7,
\hfill
$
\\[2.35ex]
где $(u^2+v^2)^3\;\!d\theta=dt.$
\vspace{0.35ex}

Круги, образующие стереографические атласы 
\vspace{0.25ex}
траекторий дифференциальных систем (9.10) и (9.13),
построены на рис. 9.5,   
\vspace{0.25ex}
дифференциальных систем (9.11) и (9.14) --- на рис. 9.6,
дифференциальных систем (9.12) и (9.15) --- на рис. 9.7.      
\\[3.75ex]
\mbox{}\hfill
{\unitlength=1mm
\begin{picture}(45,45)
\put(0,0){\includegraphics[width=45mm,height=45mm]{r09-05a.eps}}
\put(44,20.5){\makebox(0,0)[cc]{\scriptsize $x$}}
\put(20.5,44){\makebox(0,0)[cc]{\scriptsize $y$}}
\end{picture}}
\qquad\qquad
{\unitlength=1mm
\begin{picture}(45,45)
\put(0,0){\includegraphics[width=45mm,height=45mm]{r09-05b.eps}}
\put(44,20.5){\makebox(0,0)[cc]{\scriptsize $u$}}
\put(20.5,44){\makebox(0,0)[cc]{\scriptsize $v$}}
\end{picture}}
\hfill\mbox{}
\\[0ex]
\mbox{}\hfill
Рис. 9.5
\hfill\mbox{}
\\[6ex]
\mbox{}\hfill
{\unitlength=1mm
\begin{picture}(45,45)
\put(0,0){\includegraphics[width=45mm,height=45mm]{r09-06a.eps}}
\put(44,20.5){\makebox(0,0)[cc]{\scriptsize $x$}}
\put(20.5,44){\makebox(0,0)[cc]{\scriptsize $y$}}
\end{picture}}
\qquad\qquad
{\unitlength=1mm
\begin{picture}(45,45)
\put(0,0){\includegraphics[width=45mm,height=45mm]{r09-06b.eps}}
\put(44,20.5){\makebox(0,0)[cc]{\scriptsize $u$}}
\put(20.5,44){\makebox(0,0)[cc]{\scriptsize $v$}}
\end{picture}}
\hfill\mbox{}
\\[0ex]
\mbox{}\hfill
Рис. 9.6
\hfill\mbox{}
\\[6ex]
\mbox{}\hfill
{\unitlength=1mm
\begin{picture}(45,45)
\put(0,0){\includegraphics[width=45mm,height=45mm]{r09-07a.eps}}
\put(44,20.5){\makebox(0,0)[cc]{\scriptsize $x$}}
\put(20.5,44){\makebox(0,0)[cc]{\scriptsize $y$}}
\end{picture}}
\qquad\qquad
{\unitlength=1mm
\begin{picture}(45,45)
\put(0,0){\includegraphics[width=45mm,height=45mm]{r09-07b.eps}}
\put(44,20.5){\makebox(0,0)[cc]{\scriptsize $u$}}
\put(20.5,44){\makebox(0,0)[cc]{\scriptsize $v$}}
\end{picture}}
\hfill\mbox{}
\\[0ex]
\mbox{}\hfill
Рис. 9.7
\hfill\mbox{}
\\[-5ex]

\newpage

\mbox{}
\\[-2.25ex]

{\Large\bf Список литературы}
\vspace{1.75ex}

1. 
{\it Gorbuzov V.N.}
Projective atlas of trajectories of differential systems //  
Mathematics. Dy\-na\-mi\-cal Sys\-tems (arXiv: 1401.1000v1 [math.DS].  Cornell Univ., 
Ithaca, New York). -- 2014. -- P. 1 -- 61.
\vspace{0.75ex}

2. 
{\it Горбузов В.Н.}
Стереографический атлас траекторий дифференциальных систем второго  порядка //
Веснiк Гродзенскага дзяржа\u{y}нага \u{y}нiверсiтэта. Сер. 2. -- 2014. -- 
\linebreak
\No\, 1(170). -- С. 12 -- 20.
\vspace{0.75ex}

3. 
{\it Горбузов В.Н.}
Траектории стереографически сопряженных дифференциальных систем //
Веснiк Гродзенскага дзяржа\u{y}нага \u{y}нiверсiтэта. Сер. 2. -- 2014. -- \No\, 3(180). -- С. 27 -- 36.
\vspace{0.75ex}

4.
{\it Лаврентьев М.А., Шабат Б.В.} 
Методы теории функций комплексного переменного.  -- М.: Наука, 1987.   --  688 с.
\vspace{0.75ex}

5.
{\it Математическая} энциклопедия. Т. 5. -- М.: Сов. энцикл., 1984. -- 1248 с.
\vspace{0.75ex}

6.
{\it Мищенко А.С., Фоменко А.Т.}
Курс дифференциальной геометрии и топологии. -- М.: Изд-во Моск. ун-та, 1980. -- 439 с.
\vspace{0.75ex}

7.
{\it Андронов А.А., Леонтович Е.А., Гордон И.И., Майер А.Г.} 
Качественная теория динамических систем второго порядка.  -- М.: Наука, 1966.  -- 568 с.
\vspace{0.75ex}

8.
{\it Пуанкаре А.}
О кривых, определяемых дифференциальными уравнениями. -- М.;Л.: ГИТТЛ, 1947. -- 392 с.
\vspace{0.75ex}

9.
{\it Горбузов В.Н., Павлючик П.Б.} 
Линейные и разомкнутые предельные циклы дифференциальных систем  // 
Веснiк Гродзенскага дзяржа\u{y}нага \u{y}нiверсiтэта. Сер. 2. -- 2013. -- \No\, 3(159). -- С.  23 -- 32.
\vspace{0.75ex}

10.
{\it Горбузов В.Н., Самодуров А.А.} 
Уравнение Дарбу и  его аналоги.  -- Гродно:  ГрГУ, 1985.  -- 94 с.
\vspace{0.75ex}

11.
{\it Горбузов В.Н.} 
Траектории проективно приведенных дифференциальных систем // 
Веснiк Гродзенскага дзяржа\u{y}нага \u{y}нiверсiтэта. Сер. 2. -- 2012. -- \No\, 1(126). -- 
\linebreak
С. 39 -- 52.
\vspace{0.75ex}

12.
{\it Горбузов В.Н.} 
Проективный атлас траекторий дифференциальных систем второго порядка // 
Веснiк Гродзенскага дзяржа\u{y}нага \u{y}нiверсiтэта. Сер. 2. -- 2011. -- \No\, 2(111). -- С. 15 -- 26.
\vspace{0.75ex}

13.
{\it Сибирский К.С.} 
Алгебраические инварианты
дифференциальных уравнений и матриц.  -- Кишинёв:  Штиинца,
1976.  -- 268 с.
\vspace{0.75ex}

14.
{\it Лефшец С.} Геометрическая теория дифференциальных уравнений. -- М.: ИЛ, 1961. -- 387 с.

}
\end{document}